\documentclass[11pt,reqno]{amsart}

\usepackage[latin1]{inputenc}%(only for the pdftex engine)
\usepackage{microtype}
\usepackage{amssymb,amsfonts,amsmath,bm,algorithm,algorithmic,upgreek,epsfig}
\usepackage{color}
\usepackage[mathscr]{eucal}
\usepackage{hyperref}
\hypersetup{draft = false,colorlinks = true}
\usepackage{textcomp}

\newcommand{\bmu}{{\bm\mu}}
\newcommand{\bx}{{\bm x}}
\newcommand{\bn}{{\bm x}}

\newtheorem{theorem}{Theorem}

\newtheorem{remark}[theorem]{Remark}
\newtheorem{example}{Example}
\newtheorem{assumption}{Assumption}
\newtheorem{run}{Run}

\DeclareMathOperator*{\argmax}{arg\,max}
\DeclareMathOperator*{\argmin}{arg\,min}

  \title[POD model order reduction]{Model Order Reduction by\\ Proper Orthogonal Decomposition}
  \author[Gr\"a\ss{}le, Hinze \& Volkwein]{Carmen Gr\"a{\ss}le (Universit\"at Hamburg),\\ Michael Hinze (Universit\"at Koblenz-Landau), \\Stefan Volkwein (Universit\"at Konstanz)}

\begin{document}

\maketitle

{\bf Abstract.} We provide an introduction to POD-MOR with focus on (nonlinear) parametric PDEs and (nonlinear) time-dependent PDEs, and PDE constrained optimization with POD surrogate models as application. We cover the relation of POD and SVD, POD from the infinite-dimensional perspective, reduction of nonlinearities, certification with a priori and a posteriori error estimates, spatial and temporal adaptivity, input dependency of the POD surrogate model, POD basis update strategies in optimal control with surrogate models, and sketch related algorithmic frameworks. The perspective of the method is demonstrated with several numerical examples.\\
  
{\bf Key words.} POD model order reduction, (discrete) empirical interpolation, adaptivity, parametric PDEs, evolutionary PDEs, certification with error analysis.\\

{\bf Mathematics Subject Classification.} 35B30, 37M99, 41A05, 65K99, 93A15, 93C05

%%%%%%%%%%%%%%%%%%%%%%%%%%%%%%%%%%%%%%%%%%%%%%%%%%%
\section{Introduction}
%%%%%%%%%%%%%%%%%%%%%%%%%%%%%%%%%%%%%%%%%%%%%%%%%%%

Proper orthogonal decomposition (POD) is a method which comprises the essential information contained in data sets. Data sets may have their origin in various sources, like, e.g., (uncertain) measurements of geophysical processes, numerical simulations of (parameter-dependent) complex physical problems, or (dynamical) imaging. In order to illustrate the POD idea of information extraction let $\{y_1,\ldots,y_n\} \subset \mathbb R^m$ denote a vector cloud (which here serves as our data set), where we suppose at least one of the vectors $y_j$ is nonzero. Let us collect the vectors $y_j$ in the data matrix
\[
Y=[y_1\,|\ldots|\,y_n] \in \mathbb R^{m\times n}.
\]
Then we have $\mathsf r=\mathrm{rank}\,Y\in\{1,\ldots,\min(m,n)\}$. Our aim now is to find a vector $\bar\psi\in\mathbb R^m$ with length one which carries as much information of this vector cloud as possible. Of course, we here have to specify what {\em information} in this context means. For this purpose we equip $\mathbb R^m$ with some inner product $\langle\cdot\,,\cdot\rangle$ and induced norm $\|\cdot\|$. We define the information content of vector $y$ with respect to some unit vector $\psi$ by the quantity $\vert\langle y,\psi\rangle\vert$. Then we determine the special vector $\bar\psi\in\mathbb R^m$ by solving the maximization problem
\begin{equation}
\label{MaxProblem}
\bar\psi\in\argmax\bigg\{\sum_{j=1}^n \big|{\langle y_j,\psi\rangle}\big|^2\,\Big|\,\psi\in\mathbb R^m\text{ with }{\|\psi\|}=1\bigg\}.
\end{equation}
Notice that the solution to the maximization problem in \eqref{MaxProblem} is not unique. If $\bar\psi$ is a vector, where the maximum is attained, then $-\bar\psi$ is an optimal solution, too. Let us label the vector $\bar\psi$ by $\psi_1$. We now iterate this procedure; suppose that for $2\le\ell\le\mathsf r$ we have already computed such $\ell-1$ orthonormal vectors $\{\psi_i\}_{i=1}^{\ell-1}$, then seek a unit vector $\psi_\ell\in \mathbb R^m$ which is perpendicular to the $(\ell-1)$-dimensional subspace
\[
\mathscr V^{\ell-1}=\mathrm{span}\big\{\psi_1,\dots, \psi_{\ell-1}\big\}\subset\mathbb R^m,
\]
and which carries as much information of our vector cloud as possible, i.e., satisfies
\[
\psi_\ell= \argmax\bigg\{\sum_{j=1}^n \big|{\langle y_j,\psi\rangle}\big|^2\,\Big|\,\psi\in\mathbb R^m\text{ with }{\|\psi\|}=1\text{ and }\psi\perp \mathscr V^{\ell-1}\bigg\}.
\]
It is now straightforward to see that the vectors $\{\psi_i\}_{i=1}^\mathsf r$ are given by
\begin{equation}
\label{YtY1}
W^{1/2}\psi_i=\tilde\psi_i,\quad1\le i\le\mathsf r,
\end{equation}
where the $\tilde\psi_i$'s solve the eigenvalue problem (cf. \cite{GV17,HLBR12})
\[
\bar Y\bar Y^\top\tilde\psi_i=\lambda_i\tilde\psi_i,\quad i=1,\ldots,\mathsf r\text{ and }\lambda_1\ge\ldots\ge\lambda_\mathsf r>0,
\]
where $\bar Y=W^{1/2}Y\in \mathbb R^{m\times n}$ with the symmetric, positive definite (weighting) matrix 
\begin{equation}
\label{Wmatrix}
W=\big(\big({\langle e_i,e_j\rangle}\big)\big)_{1\le i,j\le m}. 
\end{equation}
In \eqref{Wmatrix} the vector $e_i$ denotes the $i$-th unit vector in $\mathbb{R}^m$. The modes $\{\psi_i\}_{i=1}^\ell$ obtained in this way are called {\em POD Modes} or {\em Principal Components} of our data cloud.
If now $m\gg n\ge\mathsf r$ it is advantageous to consider the eigenvalue problem
\[
\bar Y^\top \bar Y\phi_i = \lambda_i\phi_i,\quad i=1,\ldots,\mathsf{r}\text{ and }\lambda_1\ge\ldots\ge\lambda_\mathsf r>0,
\]
which admits the same eigenvalues $\lambda_i$ as before. The modes $\psi_i$ and $\phi_i$, $i=1,\ldots,\mathsf r$, are related by {\em singular value decomposition} (SVD):
\[
\psi_i=\frac{1}{\sigma_i}\,\bar Y\phi_i,\quad i=1,\ldots,\mathsf r,
\]
and $\sigma_i=\sqrt{\lambda_i}>0$ is the $i$-th singular value of the weighted data matrix $\bar Y$. Notice that in contrast to \eqref{YtY1} the square root matrix $W^{1/2}$ is not required.

It is now clear that a vector cloud also could be replaced by a function cloud $\{y(\bmu_j)\,|\,j=1,\ldots,n\}\subset X$ in some Hilbert space $(X,\langle\cdot\,,\cdot\rangle_X)$, where $\{\bmu_j\}_{j=1}^n$ are parameters which may refer to, e.g., time instances of a dynamic process, or stochastic variables, and the concept of information extraction by the above maximization problems directly carries over to this situation. As it is shown in the next section we can even extend this concept to general Hilbert spaces. This will be formalized in Section~\ref{P1.1} below. From the considerations above it also becomes clear that POD is closely related to SVD. This is outlined in Section~\ref{P1.2}. The POD method for abstract nonlinear evolution problems is explained in Section~\ref{P1.3}. The Hilbert space perspective also allows us to treat spatially discrete evolution equations, which include adaptive concepts for the spatial discretization. This is outlined in Section~\ref{POD-adapt}. The POD-Galerkin procedure is explained in Section~\ref{POD-Galerkin}, including a discussion of the treatment of nonlinearities. The certification of the POD method with a priori and a posteriori error bounds is outlined in Section~\ref{POD-cert}. The POD approach heavily relies on the choice of the snapshots. Related approaches are discussed in Section~\ref{POD-snap}. In Section~\ref{POD-opt} we briefly address the scope of the POD method in the context of optimal control of PDEs. Finally, in Section~\ref{POD-misc} we sketch further important research trends related to POD. Our analytical exposition is supported by several numerical experiments which give an impression of the power of the approach.

POD is one of the most successfully used model reduction techniques for nonlinear dynamical systems; see, e.g., \cite{Cha00,GV17,HLBR12,Pin08,Sir87} and the references therein. It is applied in a variety of fields including fluid dynamics, coherent structures \cite{AH01,AFS00} and inverse problems \cite{BJWW00}. Moreover in \cite{ABK01} POD is successfully applied to compute reduced-order controllers. The relationship between POD and balancing was considered in \cite{LMG02,Row05,WP02}. An error analysis for nonlinear dynamical systems in finite dimensions was carried out in \cite{RP02} and a missing point estimation in models described by POD was studied in \cite{AWWB08}.

%%%%%%%%%%%%%%%%%%%%%%%%%%%%%%%%%%%%%%%%%%%%%%%%%%%
\section{Proper Orthogonal Decomposition (POD)}
\label{P1}
%%%%%%%%%%%%%%%%%%%%%%%%%%%%%%%%%%%%%%%%%%%%%%%%%%%

In this section we introduce a discrete variant of the POD method, where we follow partially \cite[Section~1.2.1]{GV17}. For a continuous variant of the POD method and its relationship to the discrete one we refer the reader to \cite{KV02} and \cite[Sections~1.2.2 and 1.2.3]{GV17}.

%%%%%%%%%%%%%%%%%%%%%%%%%%%%%%%%%%%%%%%%%%%%%%%%%%%
\subsection{The POD method}
\label{P1.1}
%%%%%%%%%%%%%%%%%%%%%%%%%%%%%%%%%%%%%%%%%%%%%%%%%%%

Suppose that $K,n_1,\ldots,n_K$ are fixed natural numbers. Let the so-called {\em snapshot ensembles} $\{y_j^k\}_{j=1}^{n_k}\subset X$ be given for $1\le k\le K$, where $X$ is a separable real Hilbert space. For POD in complex Hilbert spaces we refer the reader to \cite{Vol01}. We set $n=n_1+\ldots+n_K$. To avoid a trivial case we suppose that at least one of the $y_j^k$'s is nonzero. Then we introduce the finite dimensional, linear {\em snapshot space}
\begin{equation}
\label{SIAM:Eq-I.1.1.1}
\mathscr V=\mathrm{span}\,\Big\{y_j^k\,\big|\,1\le j\le n_k\text{ and } 1\le k\le K\Big\}\subset X
\end{equation}
with finite dimension $d\le n$. We distinguish two cases:
\begin{itemize}
\item [1)] The separable Hilbert space $X$ has finite dimension $m$: Then $X$ is isomorphic to $\mathbb R^m$; see, e.g., \cite[p.~47]{RS80}. We define the finite index set $\mathbb I=\{1,\ldots,m\}$. Clearly, we have $1\le\mathsf r\le\min(n,m)$. Especially in case of $X=\mathbb R^m$, the snapshots $y_j^k=(y_{ij}^k)_{1\le i\le m}$ are vectors in $\mathbb R^m$ for $k=1,\ldots,K$.
\item [2)] $X$ is infinite-dimensional: Since $X$ is separable, each orthonormal basis of $X$ has countably many elements. In this case $X$ is isomorphic to the set $\ell_2$ of sequences $\{x_i\}_{i\in\mathbb N}$ of complex numbers which satisfy $\sum_{i=1}^\infty |x_i|^2<\infty$; see \cite[p.~47]{RS80}, for instance. The index set $\mathbb I$ is now the countable, but infinite set $\mathbb N$.
\end{itemize}

The {\em POD method} consists in choosing a complete orthonormal basis $\{\psi_i\}_{i\in \mathbb I}$ in $X$ such that for every $\ell\in\{1,\ldots,\mathsf r\}$ the information content of the given snapshots $y_j^k$ is maximized in the following sense:
\begin{equation}
\tag{$\mathbf P^\ell$}
\label{Pell}
\left.
\begin{aligned}
&\max\sum_{i=1}^\ell\sum_{k=1}^K\sum_{j=1}^{n_k}\alpha_j^k\,\big| {\langle y_j^k,\psi_i\rangle}_X\big|^2\\
&\hspace{0.5mm}\text{s.t. } \{\psi_i\}_{i=1}^\ell\subset X\text{ and }{\langle\psi_i,\psi_j\rangle}_X=\delta_{ij},~1 \le i,j \le \ell
\end{aligned}
\right\}
\end{equation}
with positive weighting parameters $\alpha_j^k$, $j=1,\dots,n_k$ and $k=1,\dots,K$. Here, the symbol $\delta_{ij}$ denotes the Kronecker symbol satisfying $\delta_{ii}=1$ and $\delta_{ij}=0$ for $i\neq j$.

An optimal solution $\{\Psi_i\}_{i=1}^\ell$ to \eqref{Pell} is called a {\em POD basis of rank $\ell$}. It is proved in \cite[Theorem~1.8]{GV17} that for every $\ell\in\{1,\ldots,\mathsf r\}$ a solution $\{\Psi_i\}_{i=1}^\ell$ to \eqref{Pell} is characterized by the eigenvalue problem
\begin{equation}
\label{EigPODPro}
\mathcal R\Psi_i=\lambda_i\Psi_i\quad\text{for }1\le i\le\ell,
\end{equation}
where $\lambda_1\ge\ldots\ge\lambda_\mathsf r>0$ denote the largest eigenvalues of the linear, bounded, nonnegative and self-adjoint operator $\mathcal R: X \to X$ given as
\begin{equation}
\label{OperatorR}
\mathcal R\Psi=\sum_{k=1}^K\sum_{j=1}^{n_k}\alpha_j^k\,{\langle\Psi,y_j^k\rangle}_X\,y_j^k\quad\text{for }\Psi\in X.
\end{equation}
Moreover, the operator $\mathcal R$ can be presented in the form
\begin{equation}
\label{R=YY*}
\mathcal R = \mathcal Y \mathcal Y^*
\end{equation}
with the mapping
\[
\mathcal Y: \mathbb R^n \to X,\quad\mathcal Y(\Phi) = \sum_{k=1}^K\sum_{j=1}^{n_k} \sqrt{\alpha_j^k} \phi_j^k y_j^k\quad\text{for }\Phi=\big(\phi_1^1,\ldots,\phi^K_{n_K}\big)\in\mathbb R^n,
\]
where $\mathcal Y^*:X \to \mathbb{R}^n$ denotes the Hilbert space adjoint of $\mathcal Y$, whose action is given by
\[
\mathcal Y^*(\Psi) = \left(\left\langle \Psi,\sqrt{\alpha_1^1} y_1^1\right\rangle_X, \dots, \left\langle \Psi,\sqrt{\alpha_{n_K}^K} y_{n_K}^K\right\rangle_X\right)^\top\quad\text{for }\Psi\in X.
\]
The operator $\mathcal K: \mathbb R^n\to\mathbb R^n$, $\mathcal K:=\mathcal Y^*\mathcal Y$ then admits the same nonzero eigenvalues $\lambda_1\ge\ldots\ge\lambda_\mathsf r>0$ with corresponding eigenvectors $\Phi_1, \dots,\Phi_\mathsf r$, and its action is given by
\begin{equation}
\label{OperatorK}
\mathcal K\Phi=\sum_{k=1}^K\sum_{j=1}^{n_k}\bigg(\sqrt{\alpha_1^1\alpha_j^k}\phi_j^k\,{\langle y_j^k,y_1^1\rangle}_X, \dots, \sqrt{\alpha_{n_K}^K\alpha_j^k}\phi_j^k\,{\langle y_j^k,y_{n_K}^K\rangle}_X\bigg)^\top
\end{equation}
with the vector $\Phi=(\phi_1^1,\ldots,\phi^K_{n_K})\in\mathbb R^n$. For the eigensystems of $\mathcal R$ and $\mathcal K$ there holds the relation
\begin{equation}
\label{R-K}
\Phi_i = \frac{1}{\sqrt{\lambda_i}}\,\mathcal Y^*\Psi_i, \text{ and } \Psi_i = \frac{1}{\sqrt{\lambda_i}}\,\mathcal Y\Phi_i, \quad \text{for }i=1,\dots,\mathsf r.
\end{equation}
Furthermore, we obtain
\[
\sum_{i=1}^\ell\sum_{k=1}^K\sum_{j=1}^{n_k}\alpha_j^k\,\big| {\langle y_j^k,\Psi_i\rangle}_X\big|^2=\sum_{i=1}^\ell\lambda_i,
\]
and for the POD projection error we get
\begin{equation}\label{projection_error}
\sum_{k=1}^K\sum_{j=1}^{n_k}\alpha_j^k \left\| y_j^k - \sum_{i=1}^\ell\sum_{k=1}^K\sum_{j=1}^{n_k} \langle y_j^k, \Psi_i \rangle_X \; \Psi_i\right\|_X^2 = \sum_{i=\ell+1}^\mathsf r \lambda_i.
\end{equation}
Thus, the decay rate of the positive eigenvalues $\{\lambda_i\}_{i=1}^\mathrm r$ plays an essential role for a successful application of the POD method. In general, one has to utilize a complete orthonormal basis $\{\Psi_i\}_{i\in\mathbb I}\subset X$ to represent elements in the snapshot space $\mathscr V$ by their Fourier sum. This leads to a high-dimensional or even infinite-dimensional approximation scheme. Nevertheless, if the term $\sum_{i=\ell+1}^\mathsf r\lambda_i$ is sufficiently small for a not too large $\ell$, elements in the subspace $\mathscr V$ can be approximated by a linear combination of the few basis elements $\{\Psi_i\}_{i=1}^\ell$. This offers the chance to reduce the number of terms in the Fourier series using the POD basis of rank $\ell$, as shown in the following examples. For this reason it is useful to define information content of the basis $\{\Psi_i\}_{i=1}^\ell$ in $\mathscr V$ by the quantity
\begin{equation}\label{info-content}
\mathcal E(\ell)=\frac{\sum_{i=1}^\ell\lambda_i}{\sum_{i=1}^\mathsf r\lambda_i}\in[0,1].
\end{equation}
It can e.g. be utilized to determine a basis of length $\ell\in\{1,\dots,\mathsf r\}$ containing $\approx99\%$ of the information contained in $\mathscr V$ by requiring $\mathcal E(\ell)\approx99\%$. Now it is shown in \cite[Section~1.2.1]{GV17} that
\[
\sum_{i=1}^\mathsf r\lambda_i = \sum_{k=1}^K\sum_{j=1}^{n_k}\alpha_j^k\,{\|y_j^k\|}_X^2
\]
holds true. This implies
\[
\mathcal E(\ell)=\frac{\sum_{i=1}^\ell\lambda_i}{\sum_{k=1}^K\sum_{j=1}^{n_k}\alpha_j^k\,{\|y_j^k\|}_X^2}\in[0,1],
\]
so that the quantity $\mathcal E(\ell)$ can be computed without knowing the eigenvalues $\lambda_{\ell+1},\ldots,\lambda_\mathsf r$.

%%%%%%%%%%%%%%%%%%%%%%%%%%%%%%%%%%%%%%%%%%%%%%%%%%%
\subsection{Singular value decomposition and POD}
\label{P1.2}
%%%%%%%%%%%%%%%%%%%%%%%%%%%%%%%%%%%%%%%%%%%%%%%%%%%

To investigate the relationship between singular value decomposition (SVD) and POD let us discuss the POD method for the specific case $X=\mathbb R^m$. Then we define the matrices
\begin{align*}
D^k&=\left(
\begin{array}{ccc}
\alpha_1^k&&0\\
&\ddots&\\
0&&\alpha_{n_k}^k
\end{array}
\right)\in\mathbb R^{n_k\times n_k}\quad\text{for }1\le k\le K,\\
D&=\left(
\begin{array}{ccc}
D^1&&0\\
&\ddots&\\
0&&D^K
\end{array}
\right)\in\mathbb R^{n\times n},\\
Y^k&=\big[y_1^k\,|\ldots|\,y_{n_k}^k\big]\in\mathbb R^{m\times n_k}\quad\text{for }1\le k\le K,\\
Y&=\big[Y^1\,|\ldots|\,Y^K\big]\in\mathbb R^{m\times n},\quad\bar Y=W^{1/2}YD^{1/2}\in\mathbb R^{m\times n},
\end{align*}
where we have introduced the weighting matrix $W\in\mathbb R^{m\times m}$ in \eqref{Wmatrix}.

\medskip

\begin{remark}
\em
Let us mention that $\bar Y=Y$ holds true provided all $\alpha_j^k$ are equal to one (i.e., $D$ is the identity matrix) and the inner product in $X$ is given by the Euclidean inner product (i.e., $W$ is the identity matrix).\hfill$\Diamond$
\end{remark}

Now \eqref{EigPODPro} is equivalent to the $m\times m$ eigenvalue problem%
\begin{equation}
\label{YYt}
\bar Y\bar Y^\top \bar\Psi_i=\lambda_i\bar\Psi_i\quad\text{for }1\le i\le\ell
\end{equation}
with $\Psi_i=W^{-1/2}\bar\Psi_i$ and the $n\times n$ eigenvalue problem%
\begin{equation}
\label{YtY}
\bar Y^\top \bar Y\bar\Phi_i=\lambda_i\bar\Phi_i\quad\text{for }1\le i\le\ell
\end{equation}
with $\Psi_i=YD^{1/2}\bar\Phi_i/\sqrt{\lambda_i}$. If $m\ll n$ holds, we solve \eqref{YYt}. However, we have to solve the linear system $W^{1/2}\Psi_i=\bar\Psi_i$ for any $i=1,\ldots,\ell$ in order to get the POD basis $\{\Psi_i\}_{i=1}^\ell$. Thus, if $n\le m$ holds, we will compute the solution $\{\bar\Phi_i\}_{i=1}^\ell$ to \eqref{YtY} and get the POD basis by the formula $\Psi_i=YD^{1/2}\bar\Phi_i/\sqrt{\lambda_i}$. In that case we also have $\bar Y^\top \bar Y=Y^\top WY$ so that we do not have to compute the square root matrix $W^{1/2}$. On the other hand, the diagonal matrix $D^{1/2}$ can be computed easily. The relationship between \eqref{YYt} and \eqref{YtY} is given by SVD: There exist real numbers $\sigma_1\ge\ldots\ge\sigma_\mathsf r>0$ and orthogonal matrices $\Uppsi\in\mathbb R^{m \times m}$, $\Upphi\in\mathbb R^{n\times n}$ with column vectors $\{\bar\Psi_i\}_{i=1}^m$,  $\{\bar\Phi_i\}_{i=1}^n$, respectively, such that
\begin{equation}
\label{SIAM:Eq-I.1.1.38}
\Uppsi^\top \bar Y\Upphi = \left(
\begin{array}{cc}
\Sigma^\mathsf r&0\\
0&0
\end{array}
\right)
=:\Upsigma\in\mathbb R^{m \times n},
\end{equation}
where $\Sigma^\mathsf r=\mathrm{diag}\,(\sigma_1,\ldots,\sigma_\mathsf r)\in\mathbb R^{\mathsf r\times \mathsf r}$ and the zeros in \eqref{SIAM:Eq-I.1.1.38} denote matrices of appropriate dimensions. Moreover, the vectors $\{\bar\Psi_i\}_{i=1}^\mathsf r$ and $\{\bar\Phi_i\}_{i=1}^\mathsf r$ are eigenvectors of $\bar Y\bar Y^\top $ and $\bar Y^\top \bar Y$, respectively, with eigenvalues $\lambda_i=(\sigma_i)^2>0$ for $i=1,\ldots,\mathsf r$. The vectors $\{\bar\Psi_i\}_{i=\mathsf r+1}^m$ and $\{\bar\Phi_i\}_{i=\mathsf r+1}^n$ (if $\mathsf r<m$ respectively $\mathsf r<n$) are eigenvectors of $\bar Y\bar Y^\top $ and $\bar Y^\top \bar Y$ with eigenvalue $0$. We summarize the computation of the POD basis in the pseudo code \textbf{function} \texttt{[}$\Psi$, $\Lambda$\texttt{]=}\,\texttt{POD(}$Y$, $W$, $D$, $\ell$, \texttt{flag)}.
\begin{center}
\medskip\noindent
\framebox[11.5cm]{\begin{minipage}{11.25cm}
\textbf{function} \texttt{[}$\Psi$, $\Lambda$\texttt{]=}\,\texttt{POD(}$Y$, $W$, $D$, $\ell$, \texttt{flag)}
\hrule
\vspace{0.5mm}
\begin{algorithmic}[1]
\REQUIRE Snapshots matrix $Y=[Y^1,\ldots,Y^K]$ with rank $\mathsf r$, weighting matrices $W$, $D$, number $\ell$ of POD functions and {\tt flag} for the solver;
\IF{{\tt flag} = 0}
   \STATE Set $\bar Y=W^{1/2}YD^{1/2}$;
   \STATE Compute singular value decomposition $[\Uppsi,\Upsigma,\Upphi]=\mathtt{svd}\,(\bar Y)$;
   \STATE Define $\bar\Psi_i$ as the $i$-th column of $\Uppsi$ and $\sigma_i=\Upsigma_{ii}$ for $1\le i\le\ell$;
   \STATE Set $\Psi_i=W^{-1/2}\bar\Psi_i$ and $\lambda_i=\sigma_i^2$ for $i=1,\ldots,\ell$;
\ELSIF{{\tt flag} = 1}
   \STATE Compute eigenvalue decomposition $[\Uppsi,\Uplambda]=\mathtt{eig}\,(\bar Y\bar Y^\top)$;
   \STATE Define $\bar\Psi_i$ as the $i$-th column of $\Uppsi$ and $\lambda_i=\Uplambda_{ii}$ for $1\le i\le\ell$;
   \STATE Set $\Psi_i=W^{-1/2}\bar\Psi_i$ for $i=1,\ldots,\ell$;
\ELSIF{{\tt flag} = 2}
   \STATE Compute eigenvalue decomposition $[\Upphi,\Uplambda]=\mathtt{eig}\,(\bar Y^\top \bar Y)$;
   \STATE Define $\bar\phi_i$ as the $i$-th column of $\Upphi$ and $\lambda_i=\Uplambda_{ii}$ for $1\le i\le\ell$;
   \STATE Set $\Psi_i=YD^{1/2}\bar\phi_i/\sqrt{\lambda_i}$ for $i=1,\ldots,\ell$;
\ENDIF
\RETURN $\Psi=[\Psi_1\,|\ldots|\,\Psi_\ell]$ and $\Uplambda=[\lambda_1\,|\ldots|\,\lambda_\ell]$;
\end{algorithmic}
\end{minipage}}
\end{center}

%%%%%%%%%%%%%%%%%%%%%%%%%%%%%%%%%%%%%%%%%%%%%%%%%%%
\subsection{The POD method for nonlinear evolution problems}
\label{P1.3}
%%%%%%%%%%%%%%%%%%%%%%%%%%%%%%%%%%%%%%%%%%%%%%%%%%%

In this subsection we explain the POD method for abstract nonlinear evolution problems. We focus on the numerical realization. For detailed theoretical investigations we refer the reader to \cite{GV17,HV05,HV08,KV02b,KV02}; for instance. 

%%%%%%%%%%%%%%%%%%%%%%%%%%%%%%%%%%%%%%%%%%%%%%%%%%%
\subsubsection{The nonlinear evolution problems}
\label{P1.3.1}
%%%%%%%%%%%%%%%%%%%%%%%%%%%%%%%%%%%%%%%%%%%%%%%%%%%

Let us formulate the nonlinear evolution problem. For that purpose we suppose the following hypotheses.

\begin{assumption}
\label{A1}
Suppose that $T>0$ holds, where $[0,T]$ is the considered finite time horizon.
\begin{enumerate}
\item [1)] $V$ and $H$ are real, separable Hilbert spaces and suppose that $V$ is dense in $H$ with compact embedding. By $\langle\cdot\,,\cdot\rangle_H$ and $\langle\cdot\,,\cdot\rangle_V$ we denote the inner products in $H$ and $V$, respectively. We identify $H$ with its dual (Hilbert) space $H'$ by the Riesz isomorphism so that we have the Gelfand triple
\[
V\hookrightarrow H\simeq H'\hookrightarrow V',
\]
where each embedding is continuous and dense. The last embedding is understood as follows: For every element $h \in H'$ and $v \in V$, we also have $v \in H$ by the embedding $V \hookrightarrow H$, so we can define $\langle h',v \rangle_{V',V} = \langle h',v \rangle_{H',H}$. 
\item [2)] For almost all $t\in[0,T]$ we define a time-dependent bilinear form $a(t;\cdot\,,\cdot):V\times V\to\mathbb R$ satisfying
\begin{subequations}
\label{BilForm}
\begin{align}
\label{BilForm-1}
\big|a(t;\varphi,\phi)\big|&\le\gamma\,{\|\varphi\|}_V{\|\phi\|}_V&&\text{for all }\varphi,\phi\in V,&&t\in[0,T]\text{ a.e.},\\
\label{BilForm-2}
a(t;\varphi,\varphi)&\ge\gamma_1\,{\|\varphi\|}_V^2-\gamma_2\,{\|\varphi\|}_H^2&&\text{for all }\varphi\in V,&&t\in[0,T]\text{ a.e.}
\end{align}
\end{subequations}
for time-independent constants $\gamma,\gamma_2 \ge 0$, $\gamma_1>0$ and where a.e. stands for almost everywhere.
\item [3)] Assume that $y_\circ\in V$, $f\in L^2(0,T;H)$ holds. Here we refer to {\em\cite[pp.~469-472]{DL00}} for vector-valued function spaces.
\end{enumerate}
\end{assumption}

Recall the function space
\[
W(0,T)=\big\{\varphi\in L^2(0,T;V)\,\big|\,\varphi_t\in L^2(0,T;V')\big\}
\]
which is a Hilbert space endowed with the standard inner product; cf. \cite[pp.~472-479]{DL00}. Furthermore, we have
\[
\frac{\mathrm d}{\mathrm dt}\,{\langle\varphi(t),\phi\rangle}_H={\langle\varphi_t(t),\phi\rangle}_{V',V}\quad\text{for }\varphi\in W(0,T),~\phi\in V
\]
in the sense of distributions in $[0,T]$. Here, $\langle\cdot\,,\cdot\rangle_{V',V}$ stands for the dual pairing between $V$ and its dual $V'$.

Now the evolution problem is given as follows: Find the state $y\in W(0,T)\cap C([0,T];V)$ such that
\begin{equation}
\label{State}
\begin{aligned}
\frac{\mathrm d}{\mathrm dt}\,{\langle y(t),\varphi\rangle}_H+a(t;y(t),\varphi)+{\langle\mathcal N(y(t)),\varphi\rangle}_{V',V}&={\langle f(t),\varphi\rangle}_H\\
&\quad\forall\varphi\in V,t\in(0,T]\text{ a.e.},\\
{\langle y(0),\varphi\rangle}_H&={\langle y_\circ,\varphi\rangle}_H\quad\forall\varphi\in H. 
\end{aligned}
\end{equation}

Throughout we assume that \eqref{State} admits a unique solution $y\in W(0,T)\cap C([0,T];V)$. Of course, this requires some properties for the nonlinear mapping $\mathcal N$ which we will not specify here. 

\begin{example}[Semilinear heat equation]
\label{slhe}
\em
Let $\Omega\subset\mathbb R^d$, $d\in\{2,3\}$ be a bounded open domain with Lipschitz-continuous boundary $\partial \Omega$ and let $T> 0$ be a fixed end time. We set $Q := (0,T) \times \Omega$ and $\Sigma := (0,T) \times \partial \Omega$ and $c \geq 0$. For a given forcing term $f \in L^2(Q)$ and initial condition $y_\circ\in L^2(\Omega)$, we consider the semilinear heat equation with homogeneous Dirichlet boundary condition:
\begin{equation}
\label{heat}
\left.
\begin{aligned}
y_t (t,\bx)-\Delta y(t,\bx)+c y^3(t,\bx)&=f(t,\bx)&&\text{in } Q,\\   
y(t,\bx)&=0&&\text{on } \Sigma, \\  
y(0,\bx)&=y_\circ(\bx)&&\text{in } \Omega.
\end{aligned}
\right\}
\end{equation}
The existence of a unique solution to \eqref{heat} is proved in \cite{RZ99}, for example. We can write \eqref{heat} as an abstract evolution problem of type \eqref{State} by deriving a variational formulation for \eqref{heat} with $V=H_0^1(\Omega)$ as the space of test functions, $H=L^2(\Omega)$ and integrating over the space $\Omega$. The bilinear form $a:V\times V\to\mathbb R$ is introduced by
\[
a(\varphi,\phi)=\int_\Omega\nabla\varphi\cdot\nabla\phi\,\mathrm d\bx\quad\text{for }\varphi,\phi\in V
\]
and the operator $\mathcal N:V\to V'$ is defined as $\mathcal N(\varphi)=c\varphi^3$ for $\varphi\in V$. For $c\equiv0$, the heat equation \eqref{heat} is linear.\hfill$\Diamond$
\end{example}
 
\begin{example}[Cahn-Hilliard equations]
\label{che}
\em
Let $\Omega$, $T$, $Q$ and $ \Sigma$ be defined as in Example~\ref{slhe}. The Cahn-Hilliard system was proposed in \cite{CH58} as a model for phase separation in binary alloys. Introducing the chemical potential $w$, the Cahn-Hilliard equations can be formulated in the common setting as a coupled system for the phase field $c$ and the chemical potential $w$:
\begin{equation}
\label{CHcoupled} 
\left.
\begin{aligned}
c_t(t,\bx)+y\cdot\nabla c(t,\bx)&=\mathsf m\Delta w(t,\bx)&&\text{in }Q,\\
w(t,\bx)&=-\sigma\varepsilon\Delta c(t,\bx)+\frac{\sigma}{\varepsilon} W'\big(c(t,\bx)\big)&&\text{in }Q,\\
\nabla c (t,\bx)\cdot\nu_\Omega& =\nabla w (t,\bx)\cdot\nu_\Omega=0&&\text{on }\Sigma,\\
c(0,x)&=c_\circ(\bx)&&\text{in }\Omega.
\end{aligned}
\right\}
\end{equation}
By $\nu_\Omega$ we denote the outward normal on $\partial\Omega$, $\mathsf m\geq 0$ is a constant mobility, $\sigma > 0$ denotes the surface tension and $0 < \varepsilon \ll 1$ represents the interface parameter. Note that the convective term $y \cdot \nabla c$ describes the transport with (constant) velocity $y$. The transport term represents the coupling to the Navier-Stokes equations in the context of multiphase flow, see e.g. \cite{HH77} and \cite{AGG12}. The phase field function $c$ describes the phase of a binary material with components $A$ and $B$ and takes the values $c \equiv -1$ in the pure $A$-phase and $c \equiv +1$ in the pure $B$-phase. The interfacial region is described by $c \in (-1,1)$ and admits a thickness of order $\mathcal{O}(\varepsilon)$, see e.g. Fig. \ref{fig:CH_adapt_sim}, left column, where the binary phases are colored in blue and red, respectively, and the interfacial region is depicted in white. The function $W(c)$ represents the free energy and is of double well-type. A typical choice for $W$ is the polynomial free energy function 
\begin{equation}
\label{poly}
W^p(c)=(1-c^2)^2/4 
\end{equation}
with two minima at $c=\pm 1$, which describe the energetically favorable states. It is infinitely often differentiable. Another choice for $W$ is the $C^1$ relaxed double obstacle free energy
\begin{equation}
\label{rdofe}
W^{\text{rel}}_s(c) = \frac{1}{2} (1-c^2) + \frac{s}{2}(\max(c-1,0)^2+\min(c+1,0)^2),
\end{equation}
with relaxation parameter $s\gg 0$, which is introduced in \cite{HHT11} as the Moreau-Yosida relaxation of the double obstacle free energy
\[
W^\infty(c)=
\begin{cases}
\frac{1}{2}(1-c^2),  & \text{if }c \in [-1,1],\\
+\infty, & \text{otherwise.}
\end{cases}
\]
The energies $W^p(c)$ and $W^{\text{rel}}_s(c)$ later will be used to compare the performance of POD on systems with smooth and less smooth nonlinearities. For more details on the choices for $W$ we refer to \cite{Abe07} and \cite{BE91}, for example. Concerning existence, uniqueness and regularity of a solution to \eqref{CHcoupled}, we refer to \cite{BE91}. In order to derive a variational form of type \eqref{State}, we rewrite \eqref{CHcoupled} as a single fourth-order parabolic equation for $c$ by
\begin{equation}
\label{CHsingle} 
\hspace{-2mm}\left.
\begin{aligned}
c_t(t,\bx)+y\cdot\nabla c(t,\bx)&=\mathsf m\Delta\bigg(-\sigma\varepsilon \Delta c(t,\bx)+ \frac{\sigma}{\varepsilon} W'(c(t,\bx))\bigg)&&\text{in }Q,\\
0=\nabla c(t,\bx)\cdot\nu_\Omega&=\nabla\bigg(-\sigma\varepsilon \Delta c(t,\bx)+\frac{\sigma}{\varepsilon} W'(c(t,\bx))\bigg) \cdot \nu_\Omega&&\text{on }\Sigma,\\
c(0,\bx)&=c_\circ(\bx)&&\text{in }\Omega.
\end{aligned}
\right\}
\end{equation}
We choose $V= \{ v \in H^1(\Omega): \frac{1}{|\Omega|} \int_\Omega v = 0 \}$ equipped with the inner product $(u,v)_V:= \int_\Omega \nabla u \nabla v$, so that the dual space of $V$ is given by $V'=\{f \in (H^1(\Omega))': \langle f,1 \rangle = 0\}$ such that $V \hookrightarrow H=V'$ and $\langle ., .\rangle$ denotes the duality pairing. We note that $(V,(.,.)_V)$ is a Hilbert space. We define the $V'-$inner product for $f,g\in V'$ as $(f,g)_{V'} := \int_\Omega \nabla (-\Delta)^{-1} f \cdot \nabla (-\Delta)^{-1} g$ where $(-\Delta)^{-1}$ denotes the inverse of the negative Laplacian with zero Neumann boundary data. Note that $(f,g)_{V'} = (f,(-\Delta)^{-1}g)_{L^2(\Omega)} = ((-\Delta)^{-1}f,g)_{L^2(\Omega)}$. We introduce the bilinear form $a: V \times V \to \mathbb{R}$ by 
 $$ a(u,v) = \sigma \varepsilon (\nabla u, \nabla v )_{L^2(\Omega)} + \frac{1}{\mathsf m} (y \cdot \nabla u, v)_{V'}$$
 and define the nonlinear operator $\mathcal{N}$ by $\mathcal{N}(c) = \frac{\sigma}{\varepsilon} W'(c).$
 The evolution problem can be written in the form
 \begin{equation}\label{CH-weak-abstract}
  \frac{1}{\mathsf m}(c_t(t),v)_{V'} + a(c(t),v) + \langle \mathcal{N}(c(t)),v \rangle = 0 \quad \forall v \in V \text{ and a.a. } t \in (0,T]. \footnote{We acknowledge a hint of Harald Garcke who pointed this form of the weak formulation of \eqref{CHsingle} to us.}
 \end{equation}
 We note that this fits our abstract setting formulated in \eqref{State} with the Gelfand triple $V\hookrightarrow H \equiv V' \hookrightarrow V'$. $\hfill\Diamond$
\end{example}

% \medskip

% Further examples include the two-dimensional incompressible Navier-Stokes equations, see e.g. \cite[p.~111]{Tem88}.

\bigskip

%%%%%%%%%%%%%%%%%%%%%%%%%%%%%%%%%%%%%%%%%%%%%%%%%%%
\subsubsection{Temporal discretization and POD method}
\label{P1.3.2}
%%%%%%%%%%%%%%%%%%%%%%%%%%%%%%%%%%%%%%%%%%%%%%%%%%%

Let $0=t_1<\ldots<t_{n_t}=T$ be a given time grid with step sizes $\Delta t_j=t_j-t_{j-1}$ for $j=2,\ldots,n_t$. Suppose that for any $j\in\{1,\ldots,n_t\}$ the element $y_j\in V\subset H$ is an approximation of $y(t_j)$ computed by applying a temporal integration method (e.g., the implicit Euler method) to \eqref{State}. Then we consider the snapshot ensemble
\[
\mathscr V=\mathrm{span}\,\big\{y_j\,\big|\,1\le j\le n\big\}\subset V\subset H
\]
with $n=n_t$ and $\mathsf r=\dim\mathscr V\le n$. In the context of \eqref{Pell} we choose $K=1$ and $n=n_1=n_t$. Moreover, $X$ can be either $V$ or $H$. For the weighting parameters we take the trapezoidal weights
\begin{equation}
\label{alpha}
\alpha_1=\frac{\Delta t_1}{2},\quad\alpha_j=\frac{\Delta t_j+\Delta t_{j-1}}{2}\text{ for }j=2,\ldots,n_t-1,\quad\alpha_{n_t}=\frac{\Delta t_{n_t}}{2}.
\end{equation}
Of course, other quadrature weights are also possible. Now, instead of \eqref{Pell} we consider the minimization problem
\begin{equation}
\label{PellDynSys}
\left.
\begin{aligned}
&\max\sum_{i=1}^\ell\sum_{j=1}^n\alpha_j\,\big| {\langle y_j,\Psi_i\rangle}_X\big|^2\\
&\hspace{0.5mm}\text{s.t. } \{\Psi_i\}_{i=1}^\ell\subset X\text{ and }{\langle\Psi_i,\Psi_j\rangle}_X=\delta_{ij},~1 \le i,j \le \ell
\end{aligned}
\right\}
\end{equation}
with either $X=V$ or $X=H$.

\medskip

\begin{remark}
\em
In \cite[Sections~1.2.2 and 1.3.2]{GV17} a continuous variant of the POD method is considered. In that case the trapezoidal approximation in \eqref{PellDynSys} is replaced by integrals over the time interval $[0,T]$. More precisely, we consider
\begin{equation}
\label{PellDynSysCont}
\left.
\begin{aligned}
&\max\sum_{i=1}^\ell\int_0^T\big| {\langle y(t),\Psi_i\rangle}_X\big|^2\,\mathrm dt\\
&\hspace{0.5mm}\text{s.t. } \{\Psi_i\}_{i=1}^\ell\subset X\text{ and }{\langle\Psi_i,\Psi_j\rangle}_X=\delta_{ij},~1 \le i,j \le \ell
\end{aligned}
\right\}
\end{equation}
with either $X=V$ or $X=H$. For the relationship between solutions to \eqref{PellDynSys} and \eqref{PellDynSysCont} we refer to \cite{KV02} and \cite[Section~1.2.3]{GV17}.\hfill$\Diamond$
\end{remark}

To compute the POD basis $\{\Psi_i\}_{i=1}^\ell$ of rank $\ell$ we have to evaluate the inner products $\langle y_j,\Psi_i\rangle_X$, where either $X=V$ or $X=H$ hold. In typical applications the space $X$ is usually infinite dimensional. Therefore, a discretization of $X$ is required in order to get a POD method that can be realized on a computer. This is the topic of the next subsection.

%%%%%%%%%%%%%%%%%%%%%%%%%%%%%%%%%%%%%%%%%%%%%%%%%%%
\subsubsection{Galerkin discretization}
\label{P1.3.3}
%%%%%%%%%%%%%%%%%%%%%%%%%%%%%%%%%%%%%%%%%%%%%%%%%%%

We discretize the state equation by applying any spatial approximation scheme. Let us consider here a Galerkin scheme for \eqref{State}. For this reason we are given linearly independent elements $\varphi_1,\ldots,\varphi_m\in V$ and define the $m$-dimensional subspace
\[
V^h=\mathrm{span}\,\big\{\varphi_1,\ldots,\varphi_m\big\}\subset V
\]
endowed with the $V$ topology. Then a Galerkin scheme for \eqref{State} is given as follows: find $y^h\in W(0,T)\cap C([0,T];V^h)$ satisfying
\begin{equation}
\label{StateGS}
\begin{aligned}
&\frac{\mathrm d}{\mathrm dt}\,{\langle y^h(t),\varphi^h\rangle}_H+a(t;y^h(t),\varphi^h)+{\langle\mathcal N(y^h(t)),\varphi\rangle}_{V',V}\\
&\qquad\qquad\qquad= {\langle f(t),\varphi^h\rangle}_H\quad\forall\varphi^h\in V^h,\,t\in(0,T]\text{ a.e.},\\  
&{\langle y^h(0),\varphi^h\rangle}_H={\langle y_\circ,\varphi^h\rangle}_H\quad\forall\varphi^h\in V^h. 
\end{aligned}
\end{equation}
Inserting the representation $y^h(t)=\sum_{i=1}^m\mathrm y^h_i(t)\varphi_i\in V^h$, $t\in[0,T]$, in \eqref{StateGS} and choosing $\varphi^h=\varphi_i$ for $i=1,\ldots,m$ we derive the following $m$-dimensional initial value problem
\begin{equation}
\label{StateGS-ODE}
\begin{aligned}
\mathrm M^h\dot{\mathrm y}^h(t)+\mathrm A^h(t)\mathrm y^h(t)+\mathrm N^h(\mathrm y^h(t))&=\mathrm F^h(t)\text{ for }t\in(0,T],\\ 
\mathrm M^h\mathrm y^h(0)&=\mathrm y_\circ^h,
\end{aligned}
\end{equation}
where we have used the matrices and vectors
\begin{align*}
&\mathrm M^h=\big(\big({\langle\varphi_j,\varphi_i\rangle}_H\big)\big)_{1\le i,j\le m},\hspace{8.5mm}\mathrm y^h(t)=\big(\mathrm y^h_i\big)_{1\le i\le m}\text{ for }t\in[0,T]\text{ a.e.,}\\
&\mathrm A^h(t)=\big(\big(a(t;\varphi_j,\varphi_i)\big)\big)_{1\le i,j\le m},\quad\mathrm y_\circ^h=\big({\langle y_\circ,\varphi_i\rangle}_H\big)_{1\le i\le m},\\
&\mathrm N^h(\mathrm v)=\bigg(\Big\langle\mathcal N\big({\textstyle\sum_{j=1}^m}\mathrm v_j\varphi_j\big),\varphi_i\Big\rangle_{V',V}\bigg)_{1\le i\le m}\text{ for }\mathrm v=\big(\mathrm v_j\big)_{1\le j\le m},\\&\mathrm F^h(t)=\big({\langle f(t),\varphi_i\rangle}_H\big)_{1\le i\le m}\text{ for }t\in[0,T].
\end{align*}

In the pseudo code \textbf{function} \texttt{[}$Y$\texttt{]=}\,\texttt{StateSol(}$\mathrm y_\circ^h$\texttt{)} we present a solution method for \eqref{StateGS-ODE} using the implicit Euler method.

\begin{center}
\medskip\noindent
\framebox[11.25cm]{\begin{minipage}{11cm}
\textbf{function} \texttt{[}$Y$\texttt{]=}\,\texttt{StateSol\big(}$\mathrm y_\circ^h$\texttt{\big)}
\hrule
\vspace{0.5mm}
\begin{algorithmic}[1]
\REQUIRE Initial condition $\mathrm y_\circ^h$;
\STATE Compute $\mathrm y_1^h\in\mathbb R^m$ solving $\mathrm M^h\mathrm y_1^h=\mathrm y_\circ^h$;
\FOR{$j=2$ \textbf{to} $n_t$}
   \STATE Set $\mathrm A^h_j=\mathrm A^h(t_j)\in\mathbb R^{m\times m}$ and $\mathrm F^h_j=\mathrm F^h(t_j)\in\mathbb R^m$;
   \STATE Solve (e.g., by applying Newton's method) for $\mathrm y_j^h\in\mathbb R^m$
   \[
   \big(\mathrm M^h+\Delta t_j\mathrm A_j^h\big)\mathrm y^h_j+ \Delta t_j \mathrm N^h(\mathrm y^h_j)=\mathrm M^h\mathrm y_{j-1}^h+ \Delta t_j \mathrm F^h_j;
   \]
\ENDFOR   
\RETURN matrix $Y=[\mathrm y_1^h\,|\ldots|\,\mathrm y_{n_t}^h]\in\mathbb R^{m\times n_t}$;
\end{algorithmic}
\end{minipage}}
\end{center}

In the next subsection we discuss how a POD basis $\{\Psi_j\}_{j=1}^\ell$ of rank $\ell\le\mathsf r$ can be computed from numerical approximations for the solution $y^h$ to \eqref{StateGS-ODE}.
 
%%%%%%%%%%%%%%%%%%%%%%%%%%%%%%%%%%%%%%%%%%%%%%%%%%%
\subsubsection{POD method for the fully discretized nonlinear evolution problem}
\label{P1.3.4}
%%%%%%%%%%%%%%%%%%%%%%%%%%%%%%%%%%%%%%%%%%%%%%%%%%%

Recall that we have introduced the temporal grid $\{t_j\}_{j=1}^{n_t}\subset[0,T]$ and set $n=n_t$. Let $y_1^h,\ldots,y_n^h\in V^h$ be numerical approximations to the solution $y^h(t)$ to \eqref{StateGS-ODE} at time instances $t=t_j$, $j=1\ldots,n_t$. Then, a coefficient matrix $Y\in\mathbb R^{m\times n}$ is defined by the elements $Y_{ij}$ given by
\[
y_j^h=\sum_{i=1}^mY_{ij}\varphi_i\in V^h\quad\text{for }1\le j\le n.
\]
The $j$-th column of $Y$ (denoted by $\mathrm y_j=Y_{\cdot,j}$) contains the Galerkin coefficients of the snapshot $y_j^h\in V^h$. We set $\mathsf r=\mathrm{rank}\,Y\le\min(m,n)$ and
\[
\mathscr V^h=\mathrm{span}\,\big\{y_j^h\,\big|\,1\le j\le n\big\}\subset V^h.
\]
Due to $\mathscr V^h\subset V^h$ we have $\Psi_j\in V^h$ for $1\le j\le\ell$. Therefore, there exists a coefficient matrix $\Uppsi\in\mathbb R^{m\times\ell}$ that is defined by the elements $\Uppsi_{ij}$ satisfying
\[
\Psi_j=\sum_{i=1}^m\Uppsi_{ij}\varphi_i\in V^h\quad\text{for }1\le j\le \ell,
\]
where the $j$-th column $\Uppsi_{\cdot,j}$ of the matrix $\Uppsi$  consists of  the Galerkin coefficients of the element $\Psi_j$. Notice that
\[
{\langle v^h,w^h\rangle}_H=(\mathrm v^h)^\top\mathrm M^h\mathrm w^h,\quad{\langle v^h,w^h\rangle}_H=(\mathrm v^h)^\top\mathrm S^h\mathrm w^h
\]
hold for $v^h=\sum_{i=1}^m\mathrm v_i^h\varphi_i$, $w^h=\sum_{i=1}^m\mathrm w_i^h\varphi_i\in V^h$ and for the symmetric, positive definite stiffness matrix
\[
\mathrm S^h=\big(\big({\langle\varphi_j,\varphi_i\rangle}_V\big)\big)_{1\le i,j\le m}.
\]
Then, we have for $X=H$
\[
{\langle y_j^h,\Psi_i\rangle}_X=\mathrm y_j^\top\mathrm M^h\Uppsi_{\cdot,i}=Y_{\cdot,j}^\top\mathrm M^h\Uppsi_{\cdot,i}\quad\text{for }1\le j\le n,~1\le i\le\ell,
\]
and for $X=V$
\[
{\langle y_j^h,\Psi_i\rangle}_X=\mathrm y_j^\top\mathrm S^h\Uppsi_{\cdot,i}=Y_{\cdot,j}^\top\mathrm S^h\Uppsi_{\cdot,i}\quad\text{for }1\le j\le n,~1\le i\le\ell.
\]
Thus, we can apply the approach presented in Section~\ref{P1.2} choosing $W=\mathrm M^h$ for $X=H$ and $W=\mathrm S^h$ for $X=V$. Moreover, we set $K=1$, $n_1=n_t=n$ and $\alpha^1_j=\alpha_j$ defined in \eqref{alpha}. Now a POD basis of rank $\ell$ for \eqref{StateGS-ODE} can be computed by the pseudo code \textbf{function} \texttt{[}$Y$,$\Psi$\texttt{]=}\,\texttt{PODState\big(}$\mathrm y_\circ^h$, $W$, $D$, $\ell$, \texttt{flag\big)}.

\begin{center}
\medskip\noindent
\framebox[11.25cm]{\begin{minipage}{11cm}
\textbf{function} \texttt{[}$Y$, $\Psi$\texttt{]=}\,\texttt{PODState\big(}$\mathrm y_\circ^h$, $W$, $D$, $\ell$, \texttt{flag\big)}
\hrule
\vspace{0.5mm}
\begin{algorithmic}[1]
\REQUIRE Initial condition $\mathrm y_\circ^h$, weighting matrices $W$, $D$, number $\ell$ of POD functions and {\tt flag} for the solver;
\STATE Call \texttt{[}$Y$\texttt{]=}\,\texttt{StateSol\big(}$\mathrm y_\circ^h$\texttt{\big)};
\STATE Call \texttt{[}$\Psi$,$\Lambda$\texttt{]=}\,\texttt{POD(}$Y$, $W$, $D$, $\ell$, \texttt{flag)};
\RETURN $Y=[\mathrm y_1^h\,|\ldots|\,\mathrm y_{n_t}^h]$ and $\Psi=[\Psi_1\,|\ldots|\,\Psi_\ell]$;
\end{algorithmic}
\end{minipage}}
\end{center}

In the next subsection we will discuss in detail how the POD method has to be applied in that case if we have -- instead of $V^h$ -- different spaces $V^{h_j}$ for each $j=1,\ldots,n$.

%%%%%%%%%%%%%%%%%%%%%%%%%%%%%%%%%%%%%%%%%%%%%%%%%%%
\subsection{The POD method with snapshots generated by spatially adaptive finite element methods}\label{POD-adapt}
%%%%%%%%%%%%%%%%%%%%%%%%%%%%%%%%%%%%%%%%%%%%%%%%%%%

In practical applications it often is desirable to provide POD models for time-dependent PDE systems, whose numerical treatment requires adaptive numerical techniques in space and/or time. Snapshots generated by those methods are not directly amenable to the POD procedure described in Section \ref{P1.3.4}, since the application of spatial adaptivity means, that the snapshots at each time instance may have different lengths due to their different spatial resolutions. In fact, there is not one single discrete Galerkin space $V^h$ for all snapshots generated by the fully discrete evolution, but at every time instance $t_j$ the adaptive procedure generates a discrete Galerkin space $V^{h_j} \subset X$, so that in this case $y_j^h \equiv y_j^{h_j}\in V^{h_j}$. For this reason, no snapshot matrix $Y$ can be formed with columns containing the basis coefficient vectors of the snapshots.

To obtain also a POD basis in this situation we inspect the operator $\mathcal K$ of \eqref{OperatorK} and observe that its action can be computed if the inner products $\langle y_j^k, y_i^l\rangle_X$ can be evaluated for all $1\le i \le n_l, 1\le j \le n_k$ and $1\le k,l\le K$. 

Let us next demonstrate how to compute a POD basis for snapshots residing in arbitrary finite element (FE) spaces. To begin with we drop the superindex $h$ and set $V_j:=V^{h_j}$. For each time instant $j=1, \dotsc , n$ of our time discrete PDE system the snapshots $\{ y_j\}_{j=1}^n$ are taken from different finite element spaces $V_j \subseteq X$ $(j=1,\dots,n)$, where $X$ denotes a common (real) Hilbert space. Let $V_j=\mathrm{span}\,\{\varphi_1^j,\dots,\varphi_{m_j}^j\}$. Then we have the expansions 
\begin{equation}
\label{yFE}
y_j = \sum_{i=1}^{m_j} \mathrm y_j^i\varphi_i^j\quad\text{for }j=1,\ldots,n
\end{equation}
with coefficient vectors
\[
\mathrm y_j=\big(\mathrm y_j^i\big)\in \mathbb R^{m_j}\quad\text{for }j=1,\ldots,n
\]
containing the finite element coefficients. The inner product of the associated functions can thus be computed as
\[
{\langle y_i, y_j \rangle}_X=\sum_{k=1}^{m_i} \sum_{l=1}^{m_j} \mathrm y_i^k \mathrm y_j^l\,{\langle\varphi_k^i,\varphi_l^j \rangle}_X\quad\text{for }i,j=1,\ldots,n,
\]
so that the evaluation of the action $\mathcal K \Phi$ only relies on the evaluation of the inner products $\langle\varphi_k^i,\varphi_l^j\rangle_X$ ($1\le i,j\le n$, $1\le k \le m_i$, $1\le l \le m_j$). In other words, once we are able to compute those inner products we are in the position to set up the eigensystem $\{(\lambda_i,\Phi_i)\}_{i=1}^\mathsf r$ of $\mathcal K$ from \eqref{OperatorK}. The POD modes $\{\Psi_i\}_{i=1}^\mathsf r$ can then be computed according to \eqref{R-K} by 
\[
\Psi_i = \frac{1}{\sqrt{\lambda_i}}\,\mathcal Y\Phi_i\quad\text{for }i=1,\dots,\mathsf r.
\]
Details on this procedure can be found in \cite{LU16,GH18}.

To illustrate how this procedure can be implemented we summarize Examples 6.1-6.3 from \cite{GH18}, which deal with nested and non-nested meshes. All coding was done in C++ with using FEniCS \cite{ABHJKLRW15, LMW12} for the solution of the differential equations and ALBERTA \cite{SS05} for dealing with hierarchical meshes. The numerical tests were run on a compute server with 512 GB RAM.

\begin{run}[{\cite[Example 6.1]{GH18}}]
\label{lhe-n}
\em
We consider the Example \ref{slhe} with homogeneous Dirichlet boundary condition and vanishing nonlinearity, i.e. we set $c\equiv 0$ so that the equation becomes linear. The spatial domain is chosen as $\Omega = (0,1)\times(0,1)\subset \mathbb R^2$, the time interval is $[0,T]=[0,1.57]$. Furthermore, we choose $X=L^2(\Omega)$. For the temporal discretization we introduce the uniform time grid by
\[
t_j = (j-1) \Delta t\quad\text{for }j=1,\ldots,n_t=1571
\]
with $\Delta t = 0.001$. For the spatial discretization we use $h$-adapted piecewise linear, continuous finite elements on hierarchical and nested meshes. Snapshots of the analytical solution at three different time points are shown in Figure \ref{fig:Heat_true}. Details on the construction of the analytical solution and the corresponding right hand side $f$ are given in \cite[Example 6.1]{GH18}.\\
\begin{figure}[H]
\centering
\includegraphics[width=0.3 \textwidth]{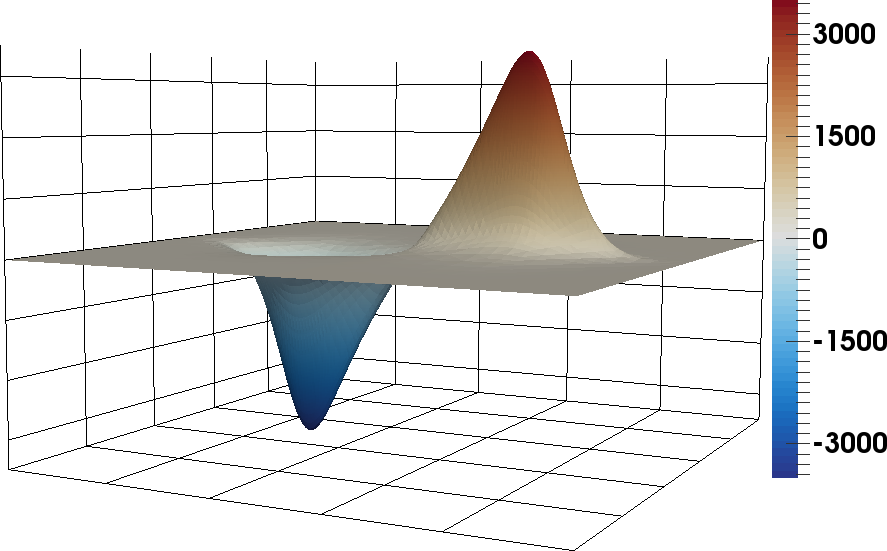} \hfill \includegraphics[width=0.3 \textwidth]{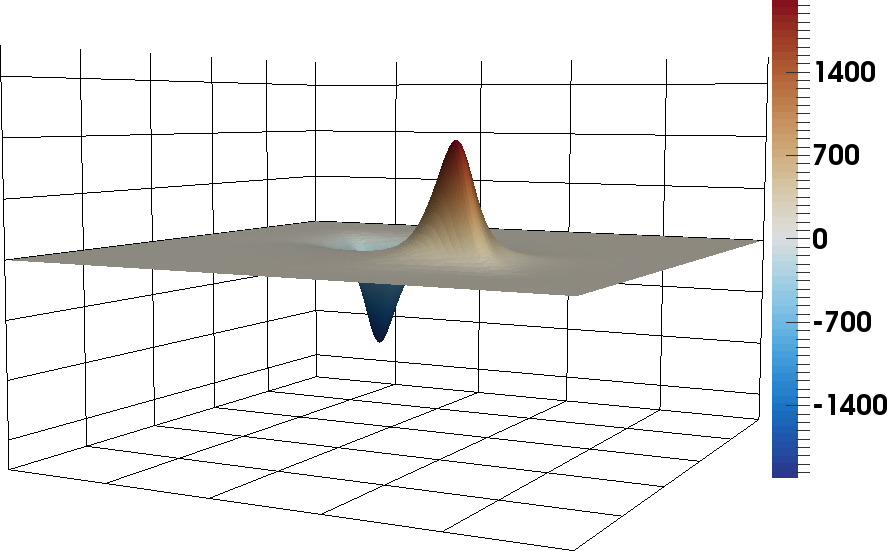}  \hfill \includegraphics[width=0.3 \textwidth]{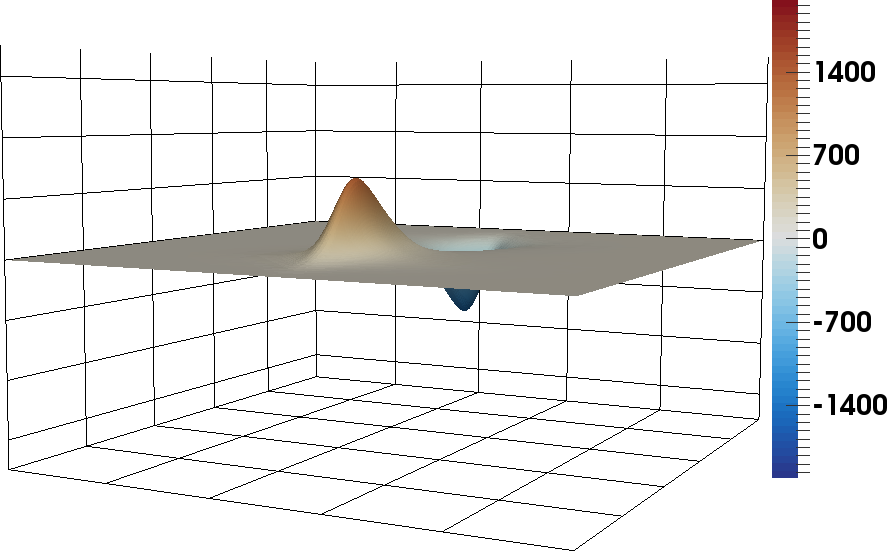} \\       
\includegraphics[width=0.33 \textwidth]{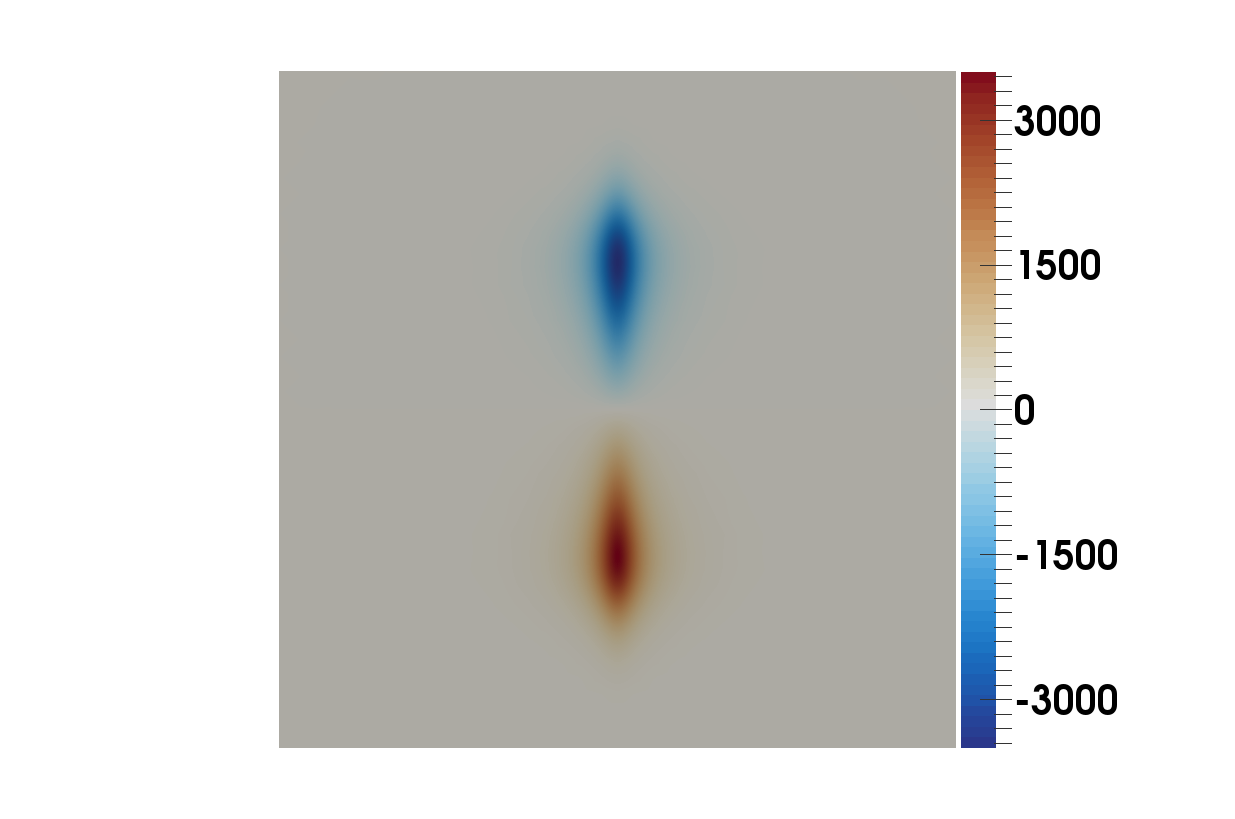}\includegraphics[width=0.33 \textwidth]{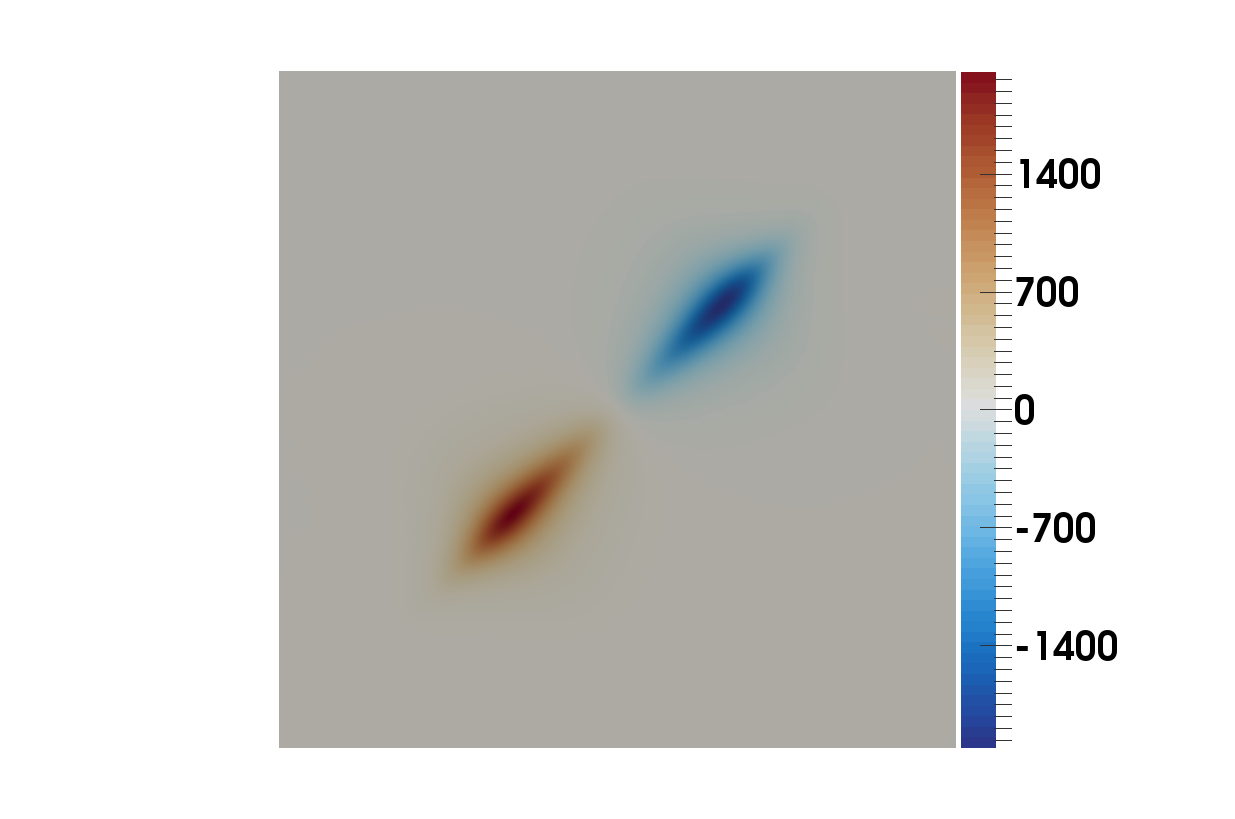}\includegraphics[width=0.33 \textwidth]{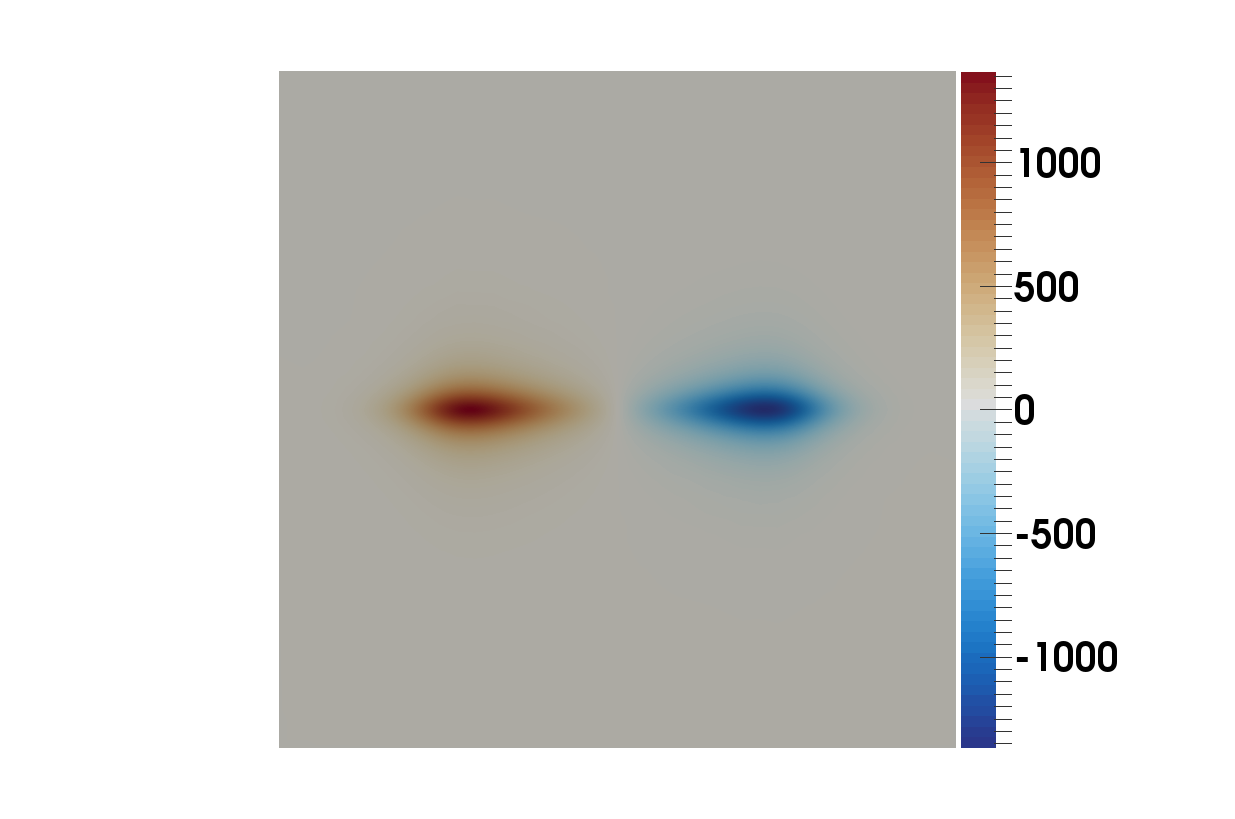}        
\caption{\em Run~{\rm\ref{lhe-n}}. Surface plot (top) and view from above (bottom) of the analytical solution of \eqref{heat} at $t=t_1$ (left), $t=T/2$ (middle) and $t=T$ (right).}
\label{fig:Heat_true}
\end{figure}

\noindent Due to the steep gradients in the neighborhood of the minimum and maximum, respectively, the use of an adaptive finite element discretization is justified. The resulting computational meshes as well as the corresponding finest mesh (reference mesh at the end of the simulation which is the union of all adaptive meshes generated during the simulation) are shown in Figure \ref{fig:Heat_meshes}.
\begin{figure}[H]
\centering
\includegraphics[width=0.24 \textwidth]{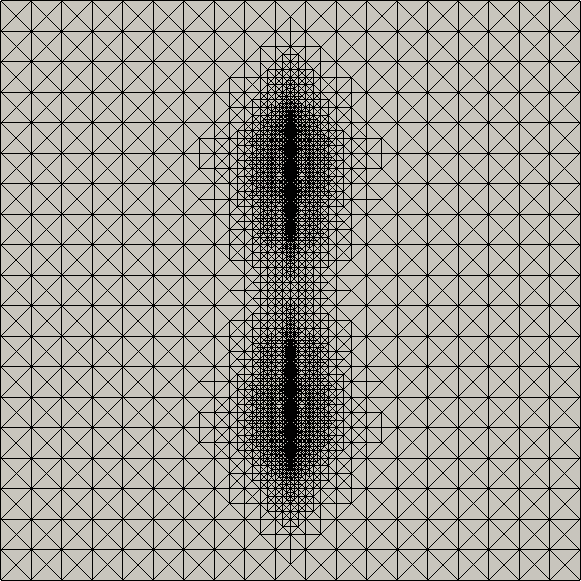} \hfill \includegraphics[width=0.24 \textwidth]{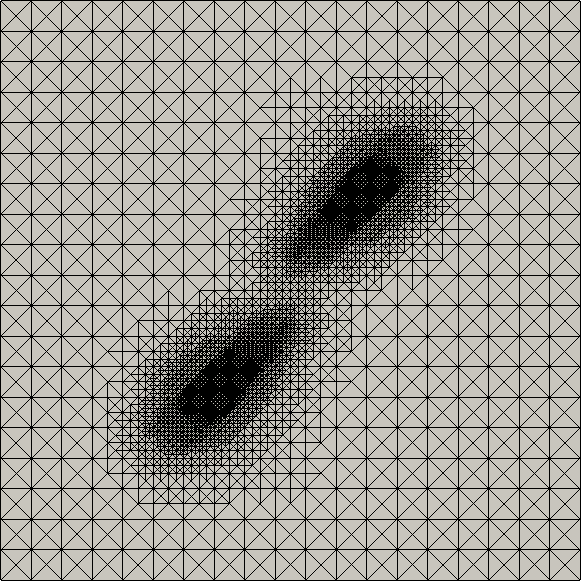} \hfill \includegraphics[width=0.24 \textwidth]{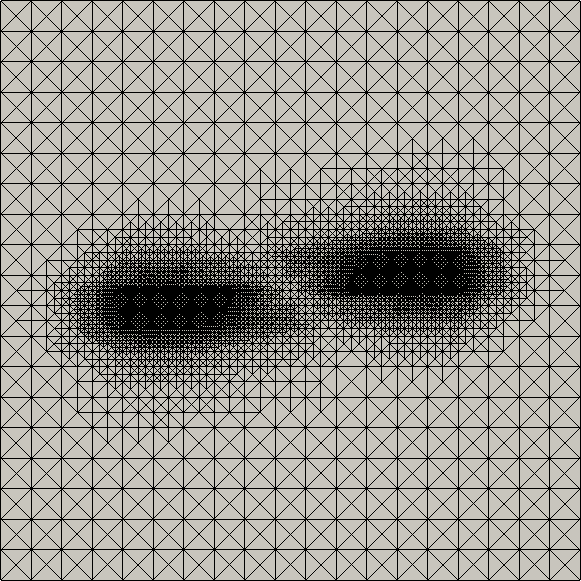}  \hfill \includegraphics[width=0.24 \textwidth]{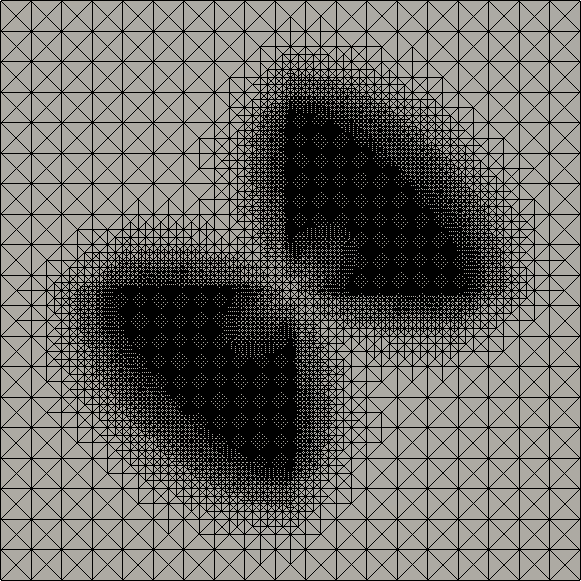} 
\caption{\em Run~{\rm\ref{lhe-n}}. Adaptive finite element meshes at $t=t_1$ (left), $t=T/2$ (middle left), $t=T$ (middle right) and finest mesh (right).}
\label{fig:Heat_meshes}
\end{figure}
\noindent The number of nodes of the adaptive meshes varies between 3637 and 7071 points. The finest mesh has 18628 degrees of freedom. A uniform mesh with grid size of order of the diameter of the smallest triangles in the adaptive grids ($h_{\min}=0.0047$) would have 93025 degrees of freedom. This clearly reveals the benefit of using adaptive meshes for snapshot generation which is also well reflected in the comparison of the computational times needed for the snapshot generation on the adaptive mesh taking 944 seconds compared to 8808 seconds on the uniform mesh, see Table \ref{tab:CPU_times_heat}) for the speedup factors obtained by spatial adaptation.
%
% \begin{table}[htbp]
% \centering
% \begin{tabular}{ l | c | c | c}
% & adaptive FE mesh  & uniform FE mesh & speedup factor\\\hline
% FE simulation & 944 sec & 8808 sec & 9.3 \\
% POD offline computations & 264 sec & 1300 sec & 4.9 \\
% POD simulation & \multicolumn{2}{c|}{ 187 sec}  & --\\\hline
% speedup factor & 5.0 &  47.1 & --\\
% \end{tabular}
% \vspace{0.4cm} \caption{\em Example~{\rm\ref{lhe-n}}: CPU times for FE and POD simulation using uniform finite element meshes and adaptive finite element snapshots which are interpolated onto the finest mesh, respectively, and using $\ell = 50$ POD modes.}
% \label{tab:CPU_times_heat}
% \end{table}
%
\noindent In Figure \ref{fig:Heat_ev}, the resulting normalized eigenspectrum of the correlation matrix $\mathrm K$ is shown for snapshots obtained by uniform spatial discretization (``uniform FE mesh''), for snapshots obtained by interpolation on the finest mesh (``adaptive FE mesh''), and for snapshots without interpolation (``infPOD''), where $\mathrm K$ is associated to the operator $\mathcal K$ from \eqref{OperatorK}, see also \eqref{matrix-k}.
\begin{figure}[H]
\centering
\includegraphics[width=0.32 \textwidth]{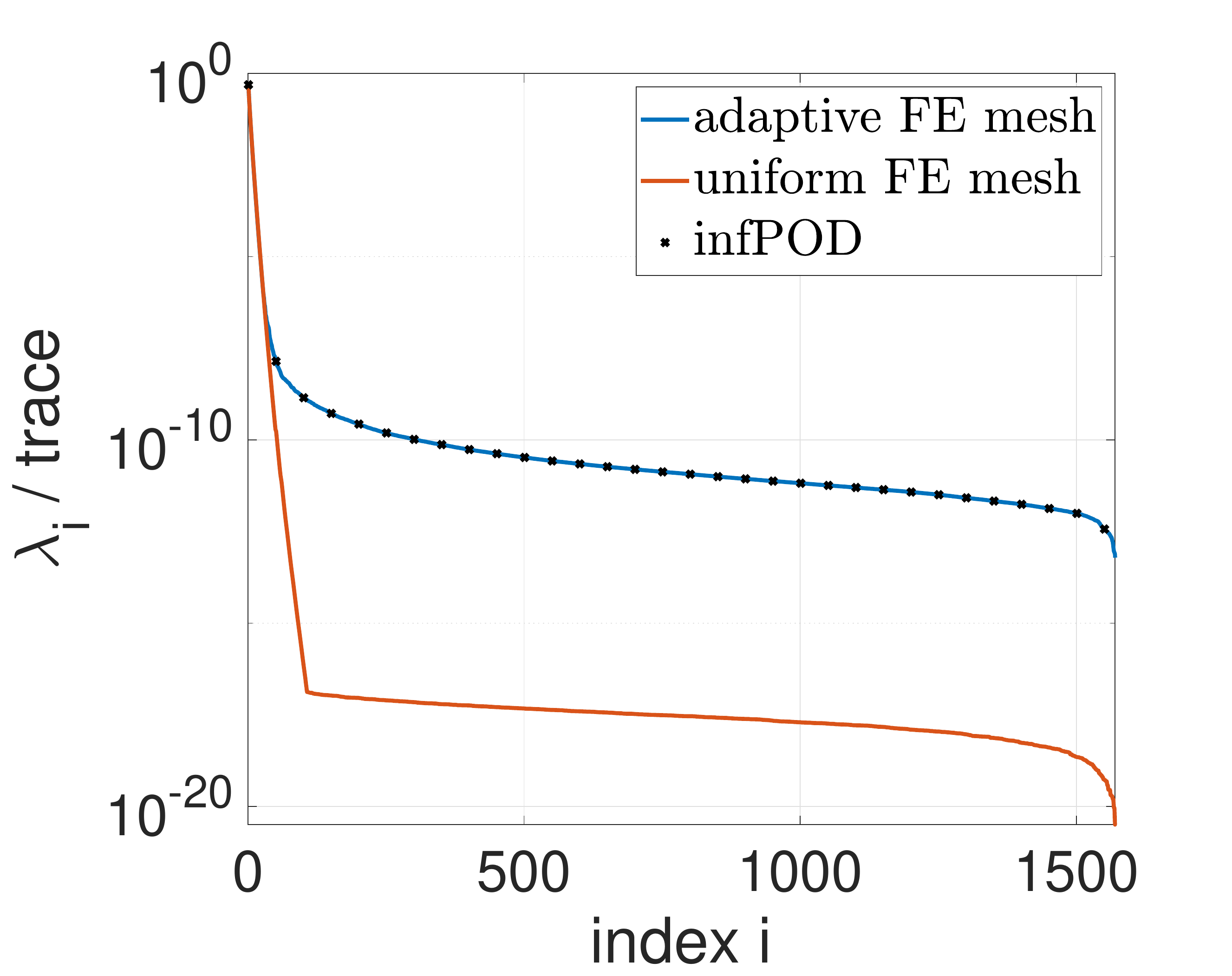}\includegraphics[width=0.32 \textwidth]{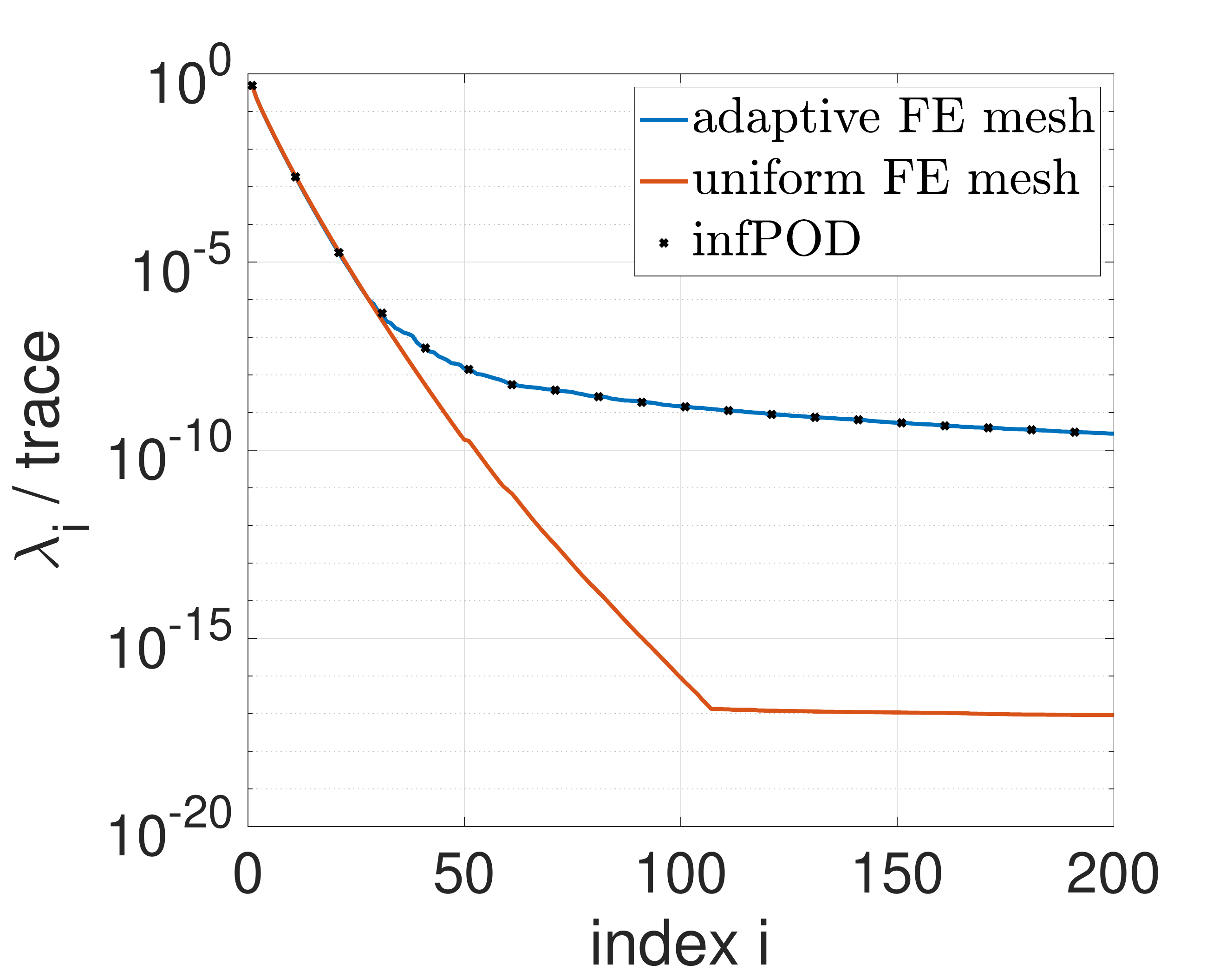} \includegraphics[width=0.351 \textwidth]{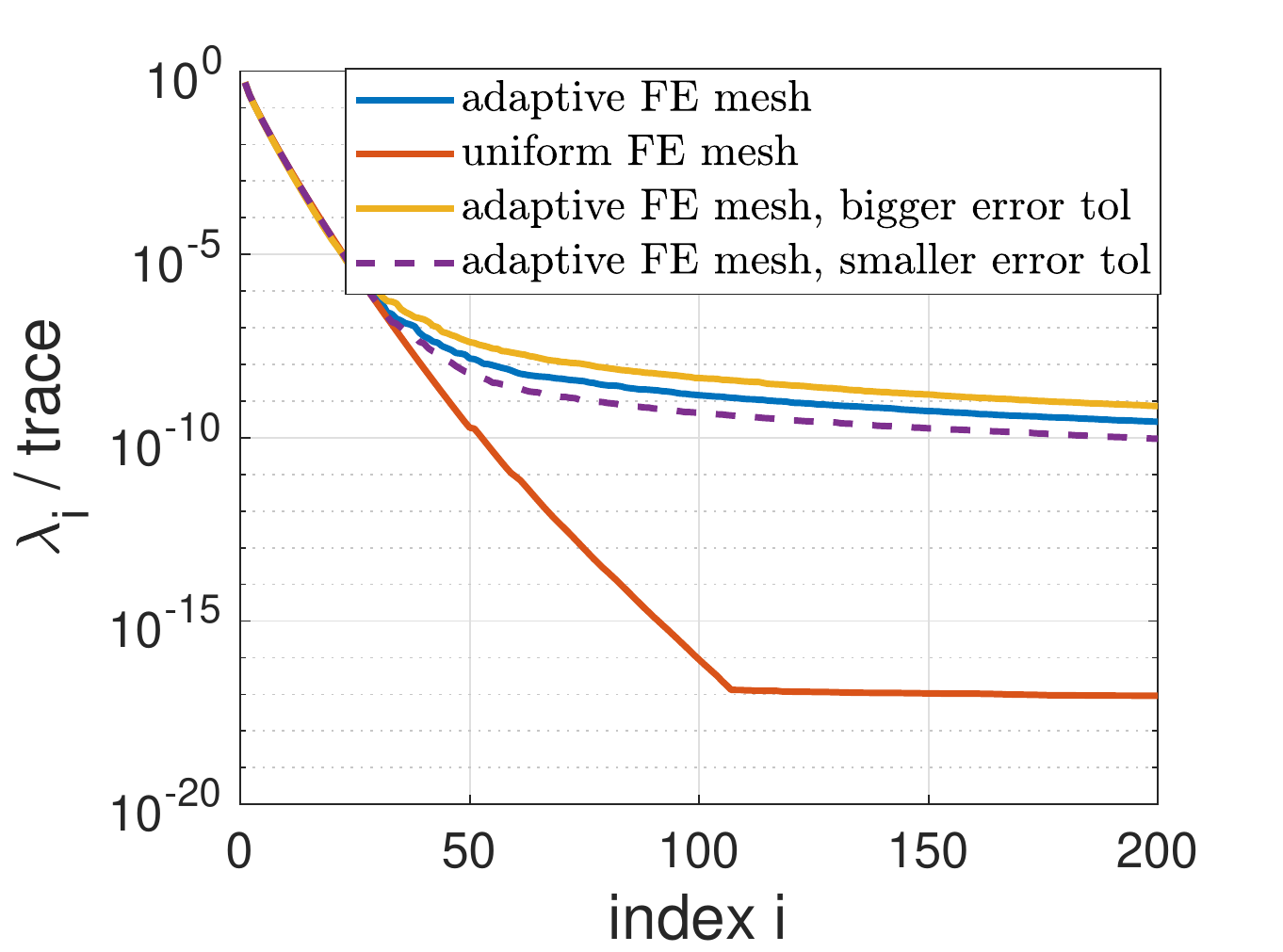}
\caption{\em Run~{\rm\ref{lhe-n}}. Comparison of the normalized eigenvalues using an adaptive and a uniform spatial mesh, respectively. Left: all eigenvalues, middle: first 200 largest eigenvalues, right: first 200 largest eigenvalues with different error tolerances for the adaptivity (1.5 times bigger and smaller error tolerances, respectively).}
\label{fig:Heat_ev}
\end{figure}
\noindent We observe that the eigenvalues for both adaptive approaches coincide. This numerically validates what we expect from theory: the information content which is contained in the matrix $\mathcal{K}$ when we explicitly compute the entries without interpolation is the same as the information content contained within the eigenvalue problem which is formulated when using the finest mesh. No information is added or lost. Moreover, we recognize that about the first 28 eigenvalues computed corresponding to the adaptive simulation coincide with the simulation on a uniform mesh. From index 29 on, the methods deliver different results: for the uniform discretizations, the normalized eigenvalues fall below machine precision at around index 100 and stagnate. In contrary, the normalized eigenvalues for both adaptive approaches flatten in the order around $10^{-10}$. If the error tolerance for the spatial discretization error is chosen larger (or smaller), the stagnation of the eigenvalues in the adaptive method takes place at a higher (or lower) order (see Figure \ref{fig:Heat_ev}, right). Concerning dynamical systems, the magnitude of the eigenvalue corresponds to the characteristic properties of the underlying dynamical system: the larger the eigenvalue, the more information is contained in the corresponding eigenfunction. Since all adaptive meshes are contained in the uniform mesh, the difference in the amplitude of the eigenvalues is due to the interpolation errors during refinement and coarsening. This is the price we have to pay for faster snapshot generation using adaptive methods. A further aspect gained from the decay behavior of the eigenvalues in the adaptive case is the following; the adaptive approach filters out the noise in the system which is related to the modes corresponding to the singular values that are not matched by the eigenvalues of the adaptive approach. This in the language of frequencies means that the overtones in the systems which get lost in the adaptive computations live in the space which is neglected by the POD method based on adaptive finite element snapshots. From this point of view, adaptivity can be interpreted as a smoother. 

The first, second and fifth POD modes of Run \ref{lhe-n} obtained by the adaptive approach are depicted in Figure \ref{fig:Heat_PODmodes}. We observe the classical appearance of the basis functions. The initial condition is reflected by the first POD basis function. The next basis functions admit a number of minima and maxima corresponding to the index in the basis: $\Psi_2$ has two minima and two maxima etc. This behavior is similar to the increasing oscillations in higher frequencies in trigonometric approximations. The POD basis functions corresponding to the uniform spatial discretization have a similar appearance.\hfill$\Diamond$\\
\begin{figure}[H]
\centering
\includegraphics[width=0.3 \textwidth]{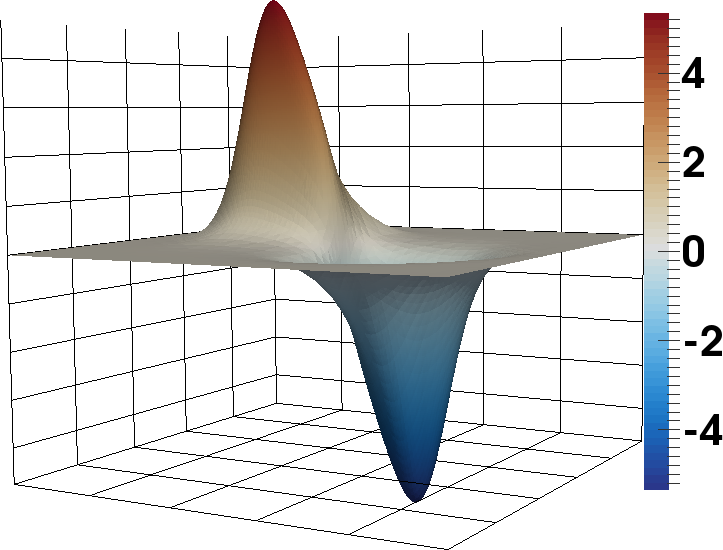} \hfill  \includegraphics[width=0.3 \textwidth]{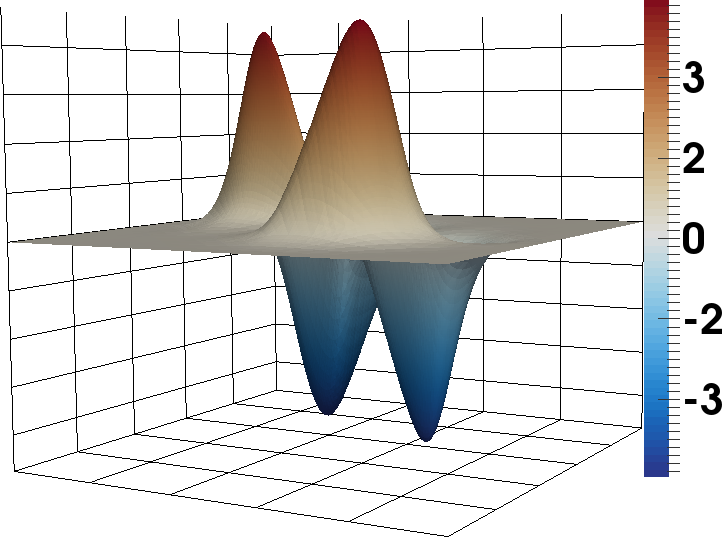} \hfill \includegraphics[width=0.3 \textwidth]{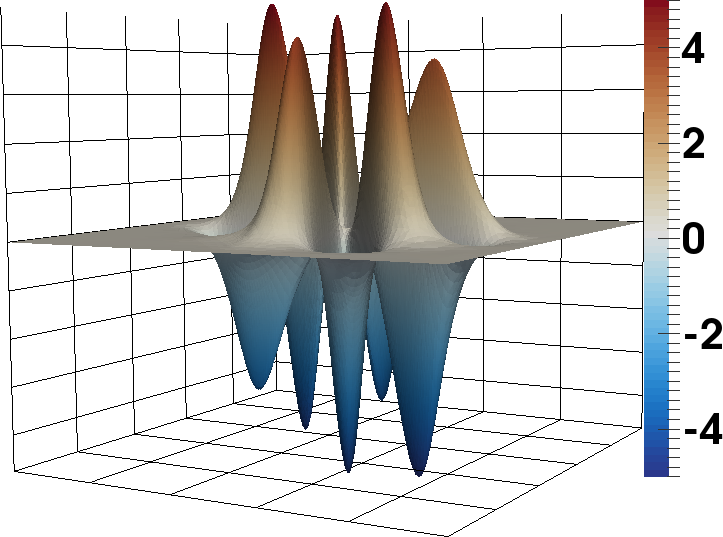} \\ 
\includegraphics[width=0.32 \textwidth]{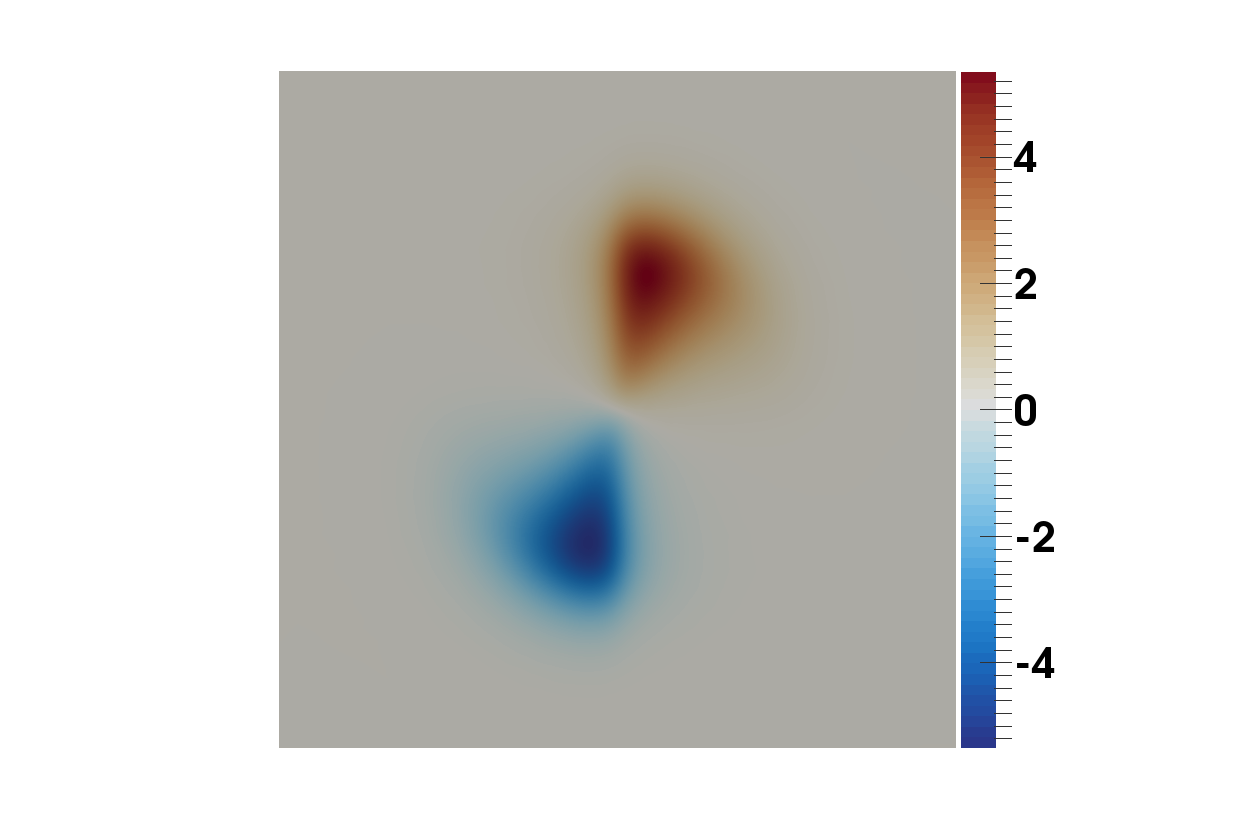} \includegraphics[width=0.32 \textwidth]{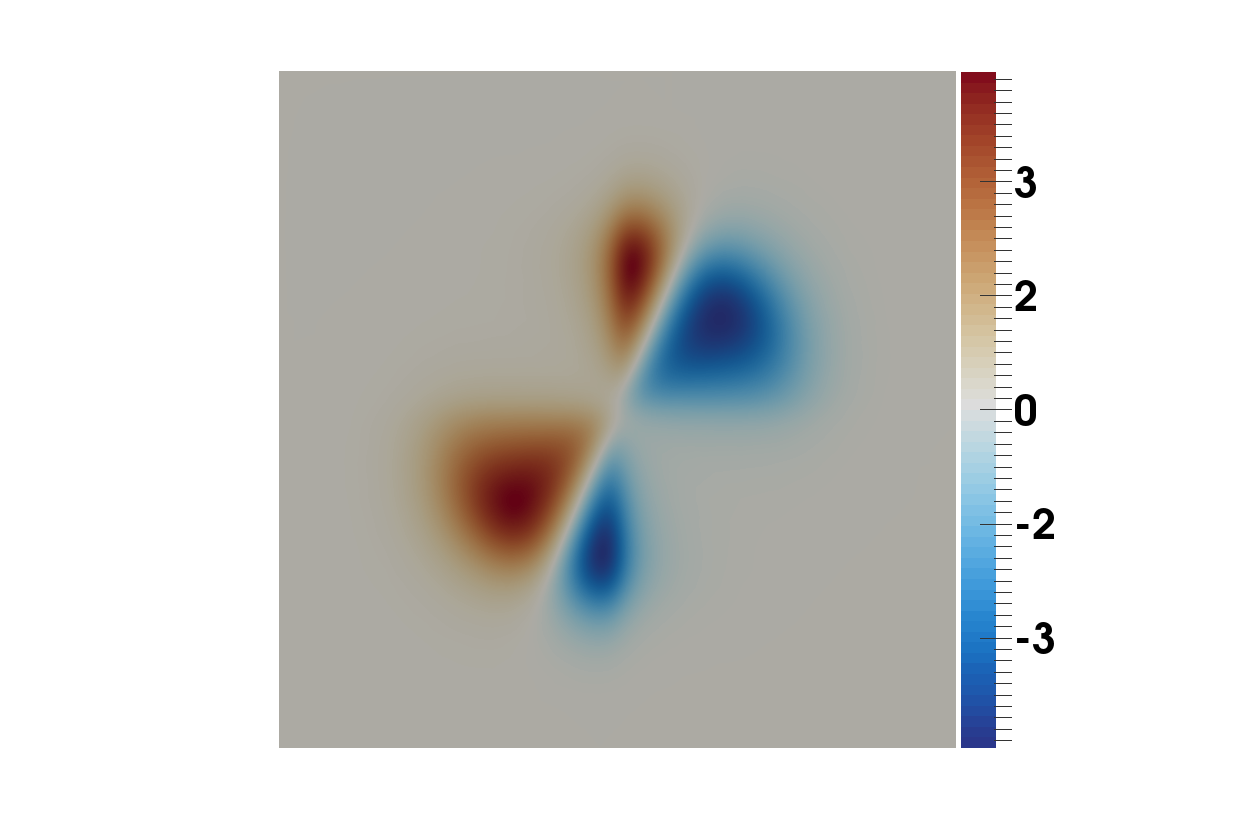} \includegraphics[width=0.32 \textwidth]{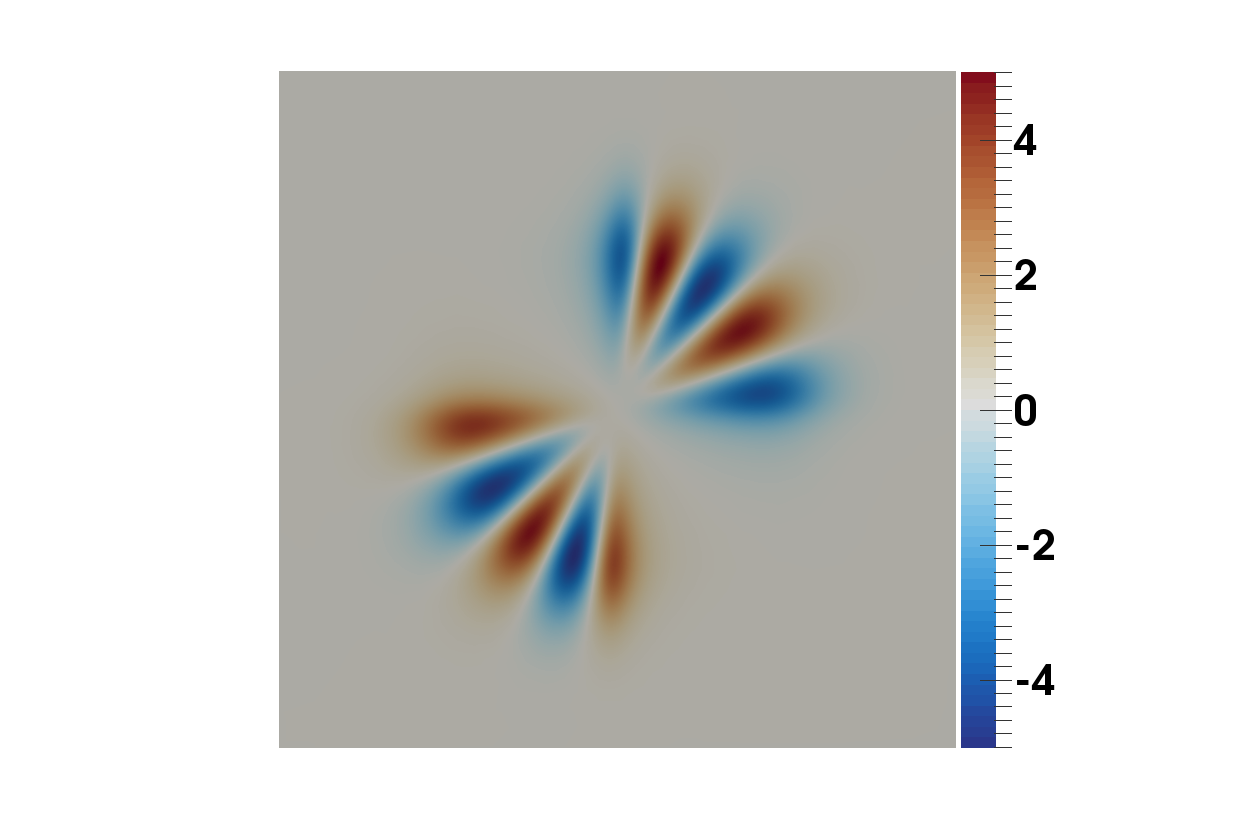} 
\caption{\em Run~{\rm\ref{lhe-n}}. Surface plot (top) and view from above (bottom) of the POD basis functions $\Psi_1$ (left), $\Psi_2$ (middle) and $\Psi_5$ (right).}
\label{fig:Heat_PODmodes}
\end{figure}

\end{run}

\begin{run}[{\cite[Example 6.2]{GH18}}]\label{ch-equation} (Cahn-Hilliard system.) 
\em
We consider Example \ref{che} in the form \eqref{CHcoupled} with $\Omega=(0,1.5) \times (0, 0.75)$, $T=0.025$, constant mobility $\mathsf m\equiv 0.00002$, and  constant surface tension $\sigma\equiv 24.5$. The interface parameter $\varepsilon$ is set to $\varepsilon=0.02$, with resulting interface thickness $\pi\cdot\varepsilon\approx 0.0628$. We use the relaxed double obstacle free energy $W_s^{\text{rel}}$ from \eqref{rdofe} with $s=10^4$. As initial condition, we choose a circle with radius $r = 0.25$ and center $(0.375,0.375)$. The initial condition is transported horizontally with constant velocity $v = (30,0)^T$. We set 
\[
t_j = (j-1) \Delta t\quad\text{for }j=1,\ldots,n_t=1001,
\]
so that $\Delta t=2.5\cdot 10^{-5}$. The numerical computations are performed with the semi-implicit Euler scheme. For this purpose let $c^{j-1} \in V$ and $c^j \in V$ denote the time-discrete solution at $t_{j-1}$ and $t_j$, respectively. Based on the variational formulation \eqref{CH-weak-abstract}  we tackle the time discrete version of \eqref{CHcoupled} in the form: given $c^{j-1}$, find $c^j$, $w^j$ solving
\begin{equation}\label{CH_num}
\left.
\begin{aligned}
\frac{1}{\Delta t}\,{\langle c^j-c^{j-1},\varphi_1\rangle}_{L^2} +{\langle v\cdot\nabla c^{j-1},\varphi_1\rangle}_{L^2} +\mathsf m\,{\langle \nabla w^j, \nabla \varphi_1\rangle}_{L^2}& = 0,\\
-{\langle w^j,\varphi_2\rangle}_{L^2}+\sigma\varepsilon\,{\langle\nabla c^j,\nabla\varphi_2\rangle}_{L^2}+\frac{\sigma}{\varepsilon}\,{\langle W'_+(c^j) + W'_-(c^{j-1}),\varphi_2\rangle}_{L^2}& =0
\end{aligned}
\right\}
\end{equation}
for all $\varphi_1,\varphi_2 \in V$ and $j=2,\ldots,n_t$ with $c^1=c_\circ$. According to \eqref{CH-weak-abstract}, here it is $V=\{v\in H^1(\Omega), \frac{1}{|\Omega|} \int_\Omega v dx = 0\}$. Note that the free energy function $W$ is split into a convex part $W_+$ and a concave part $W_-$, such that $W = W_+ + W_-$ and $W'_+$ is treated implicitly, whereas $W'_-$ is treated explicitly with respect to time. This leads to an unconditionally energy stable time marching scheme, compare \cite{Eyr98}. The system \eqref{CH_num} is discretized in space using piecewise linear and continuous finite elements. The resulting nonlinear equation systems are solved using a semi-smooth Newton method.\\

\noindent Figure \ref{fig:CH_adapt_sim} shows the phase field (left) and the chemical potential (right) for the finite element simulation using adaptive meshes. The initial condition $c_\circ$ is transported horizontally with constant velocity.\\
\begin{figure}[H]
\centering
\includegraphics[width=.45\textwidth]{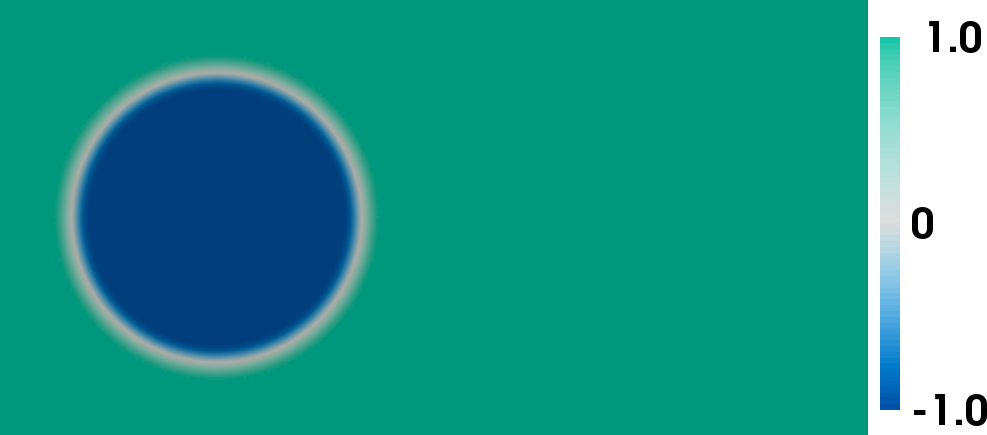} \hfill\includegraphics[width=.45\textwidth]{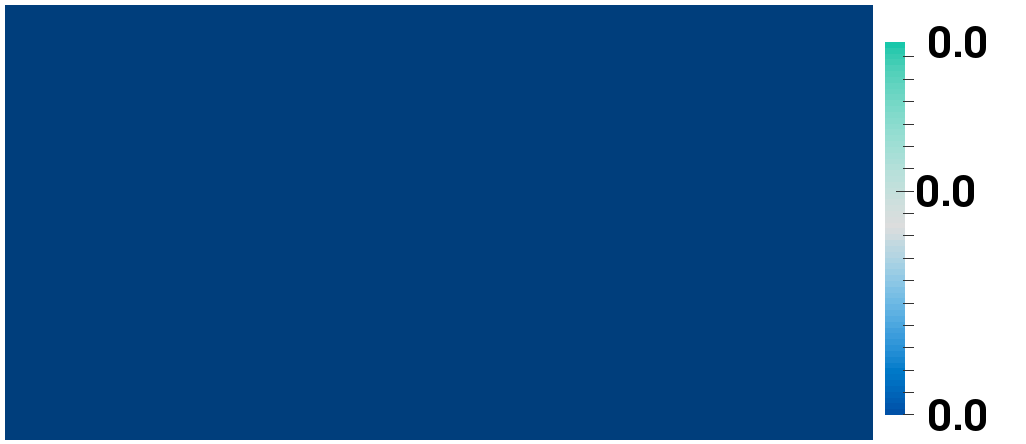}\\
\includegraphics[width=.45\textwidth]{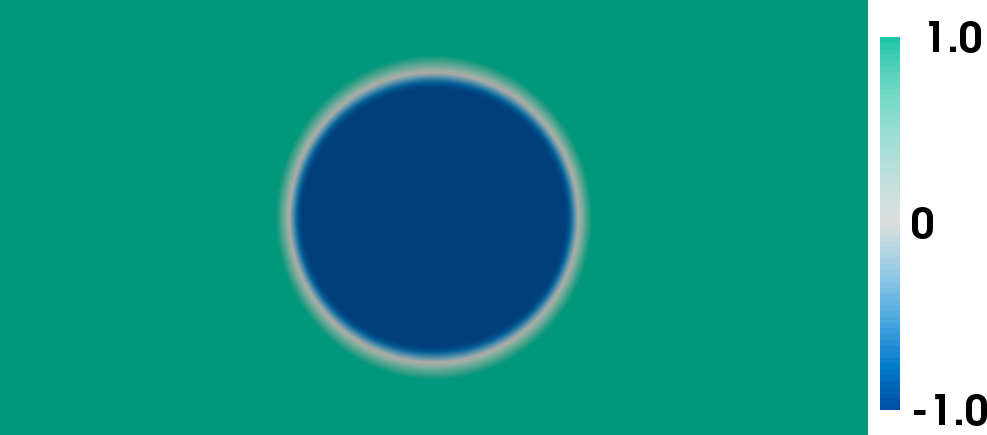} \hfill \includegraphics[width=.45\textwidth]{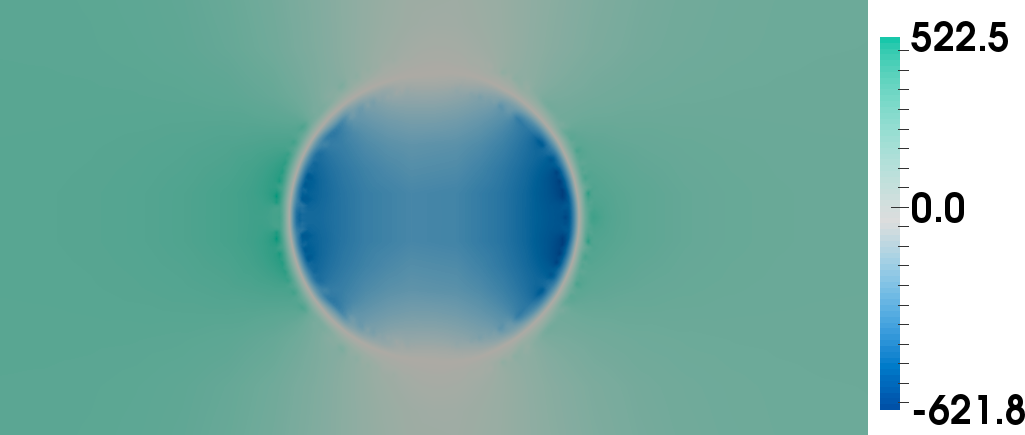}\\
\includegraphics[width=.45\textwidth]{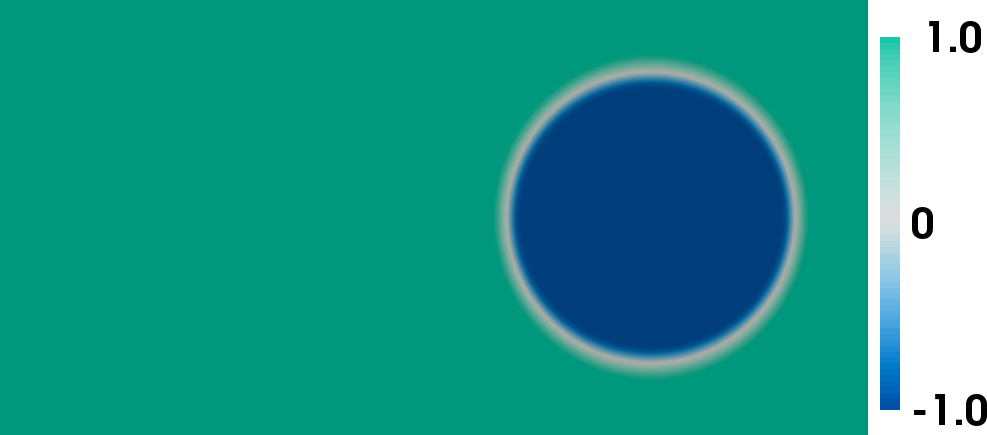} \hfill \includegraphics[width=.45\textwidth]{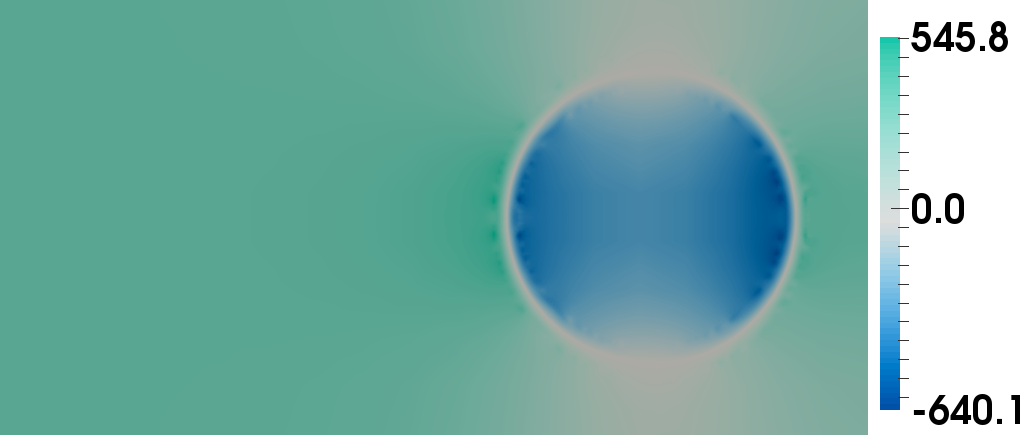}
\caption{\em Run~{\rm\ref{ch-equation}}. Phase field $c$ (left) and chemical potential $w$ (right) computed on adaptive finite element meshes at $t=t_1$ (top), $t=T/2$ (middle) and $t=T$ (bottom).}
\label{fig:CH_adapt_sim}
\end{figure}
 
\noindent The adaptive finite element meshes as well as the finest mesh which is generated during the adaptive finite element simulation are shown in Figure \ref{fig:CH_grids}. The number of degrees of freedom in the adaptive meshes varies between 6113 and 8795. The finest mesh (overlay of all adaptive meshes) has 54108 degrees of freedom, whereas a uniform mesh with discretization fineness as small as the smallest triangle in the adaptive meshes has 88450 degrees of freedom.\\
\begin{figure}[H]
\centering
\includegraphics[width=.45\textwidth]{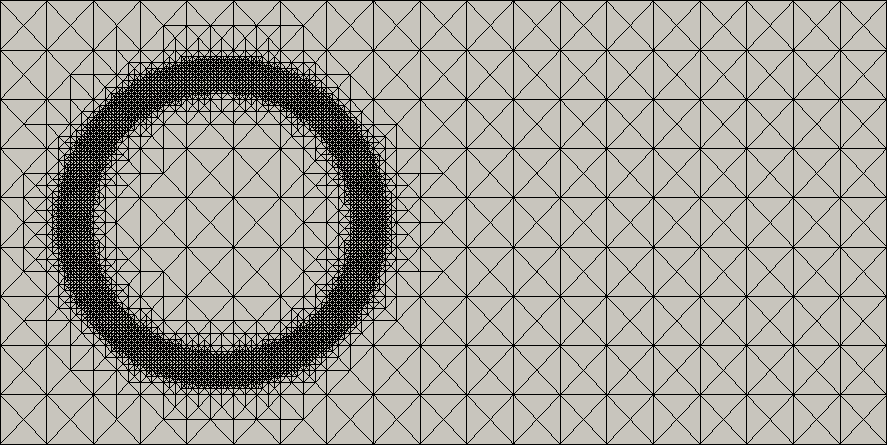} \hfill \includegraphics[width=.45\textwidth]{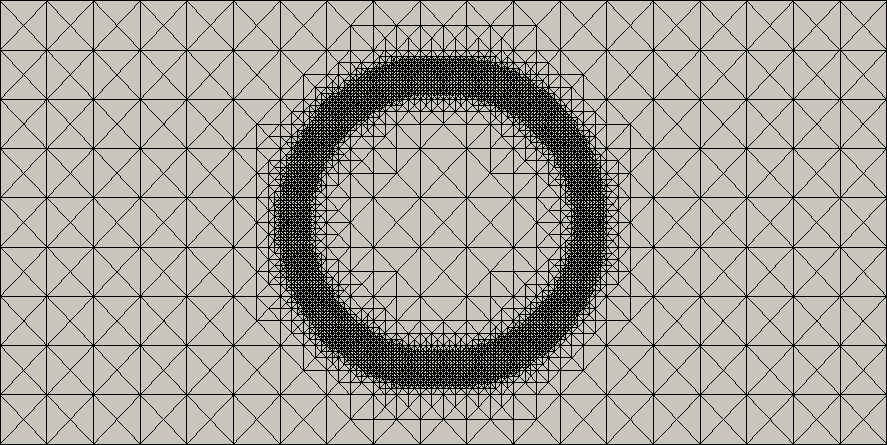}\\
\includegraphics[width=.45\textwidth]{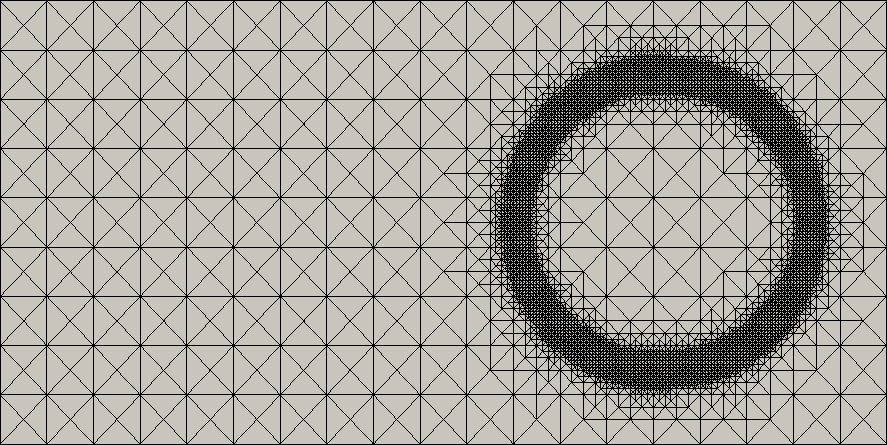} \hfill \includegraphics[width=.45\textwidth]{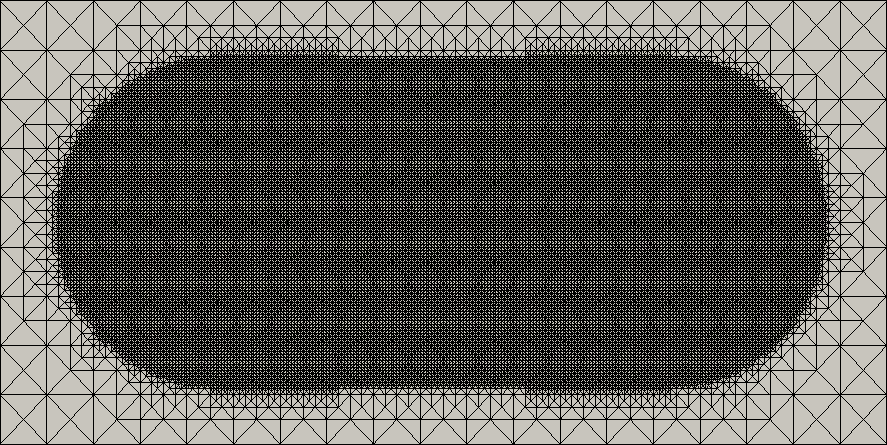}
\caption{\em Run~{\rm\ref{ch-equation}}. Adaptive finite element meshes at $t=t_1$ (top left), $t=T/2$ (top right) and $t=T$ (bottom left) together with the finest mesh (bottom right).}
\label{fig:CH_grids}
\end{figure}

\noindent Figure \ref{fig:CH_POD_bases} shows the first, second and fifth POD mode for the phase field $c$ and the chemical potential $w$. Analogously to Run \ref{lhe-n}, we observe a periodicity in the POD basis functions corresponding to their basis index numbers.\\
\begin{figure}[H]
\centering
\includegraphics[width=.45\textwidth]{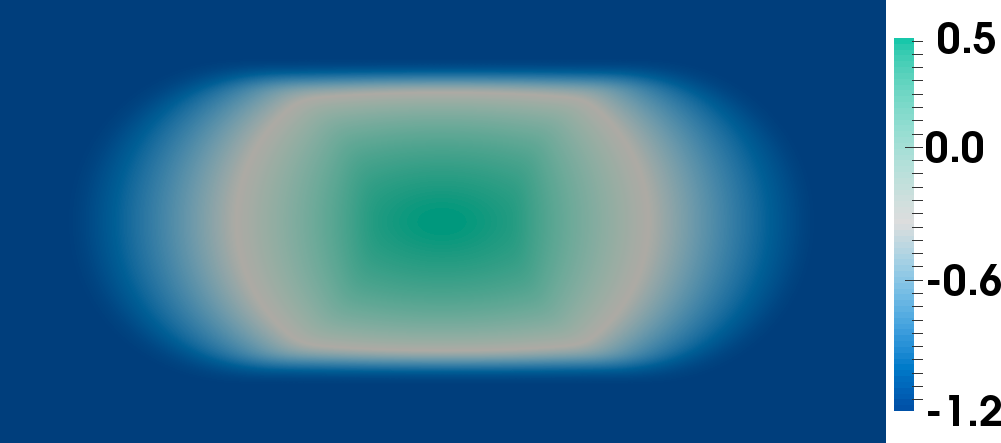} \includegraphics[width=.45\textwidth]{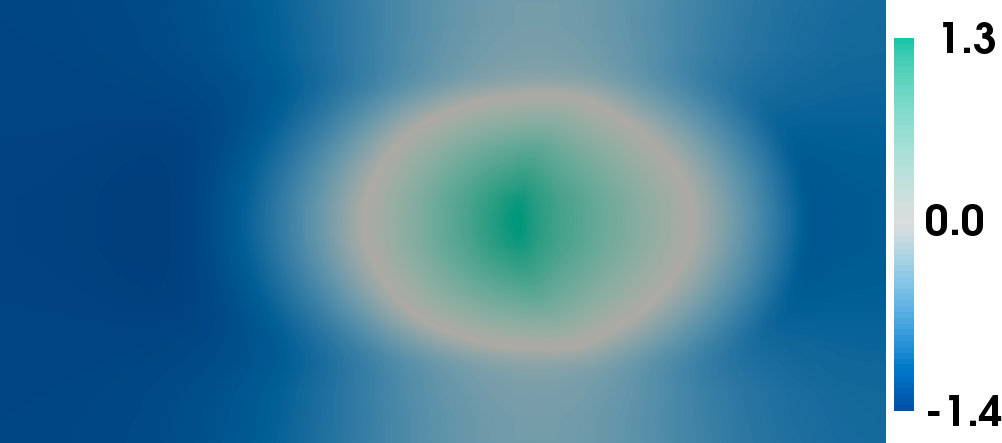}\\
\includegraphics[width=.45\textwidth]{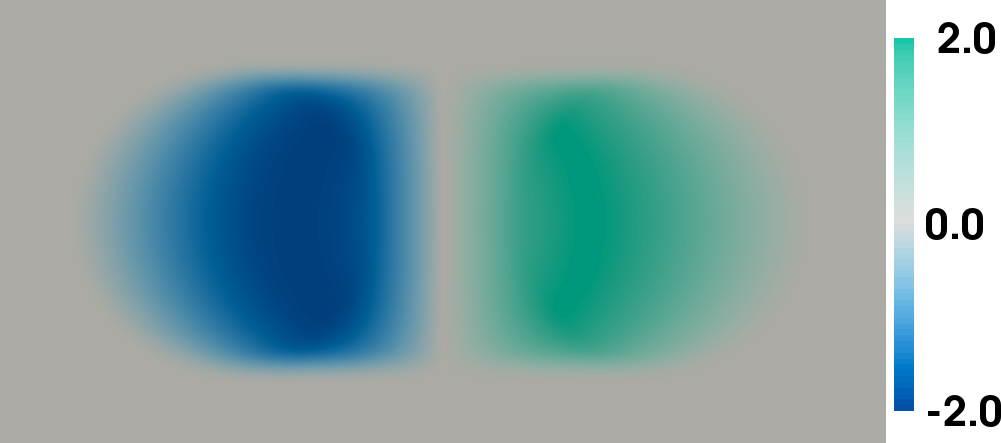} \includegraphics[width=.45\textwidth]{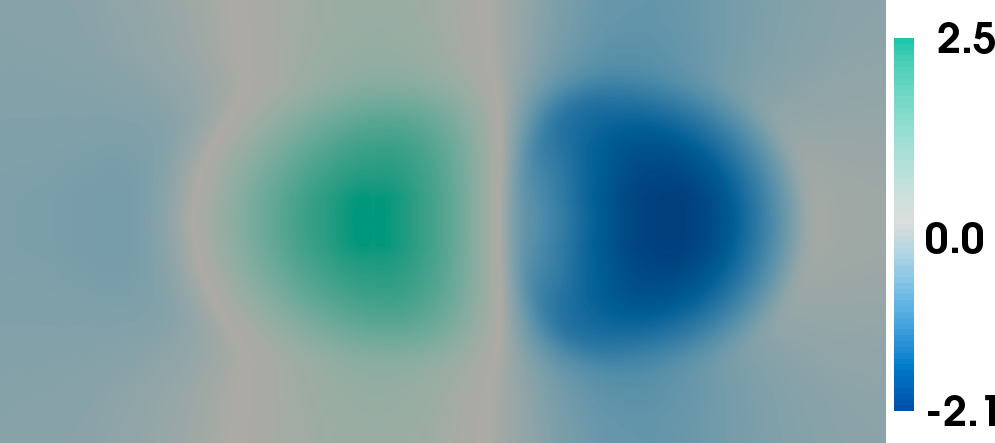}\\
\includegraphics[width=.45\textwidth]{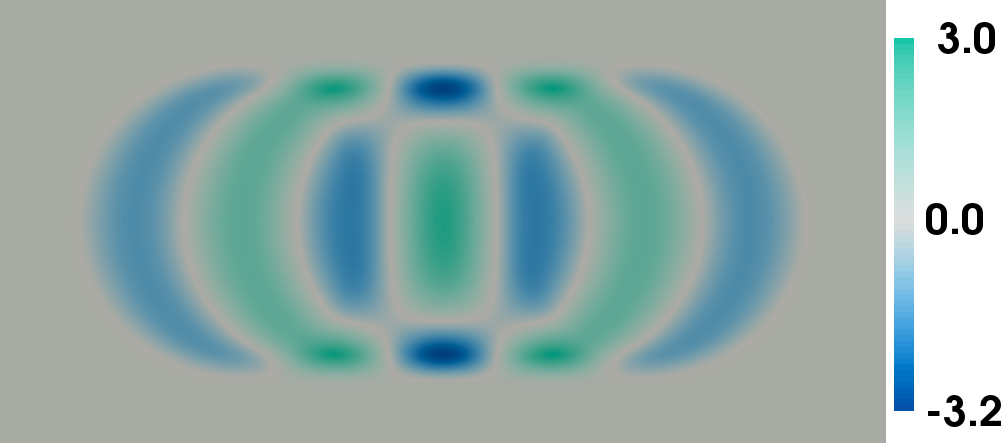} \includegraphics[width=.45\textwidth]{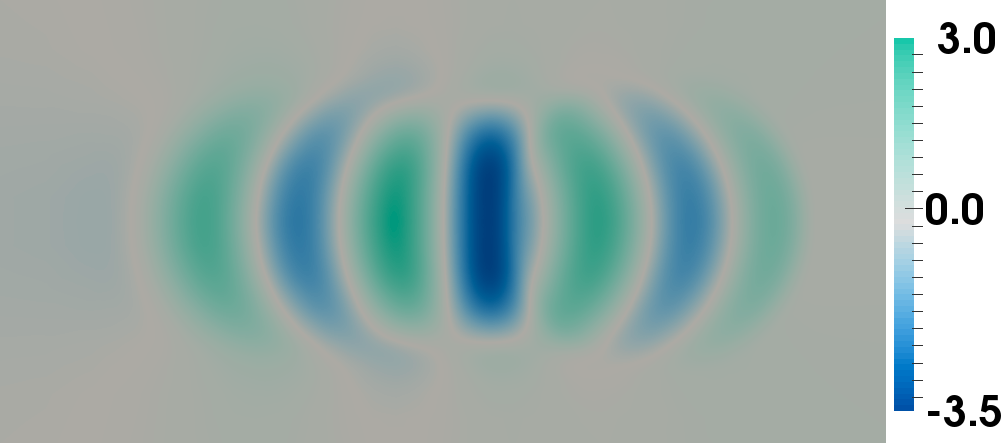}
\caption{\em Run~{\rm\ref{ch-equation}}. First, second and fifth POD modes for $c$ (left) and $w$ (right).}
\label{fig:CH_POD_bases}
\end{figure}

\noindent In the present example we only compare the POD procedure for two kinds of snapshot discretizations, namely the adaptive approach with using a finest mesh, and the uniform mesh approach, where the gridsize is chosen to be of the same size as the smallest triangle in the adaptive meshes. We choose $X=L^2(\Omega)$ and compute a separate POD basis for each of the variables $c$ and $w$.
 
\noindent In Figure \ref{fig:CH_ev}, a comparison is visualized concerning the normalized eigenspectrum for the phase field $c$ and the chemical potential $w$ using uniform and adaptive finite element discretization. We note for the phase field $c$ that about the first 180 eigenvalues computed corresponding to the adaptive simulation coincide with the eigenvalues of the simulation on the finest mesh. Then, the eigenvalues corresponding to the uniform simulation decay faster. Similar observations apply for the chemical potential $w$. 
\begin{figure}[H]
\centering
\includegraphics[width=.46\textwidth]{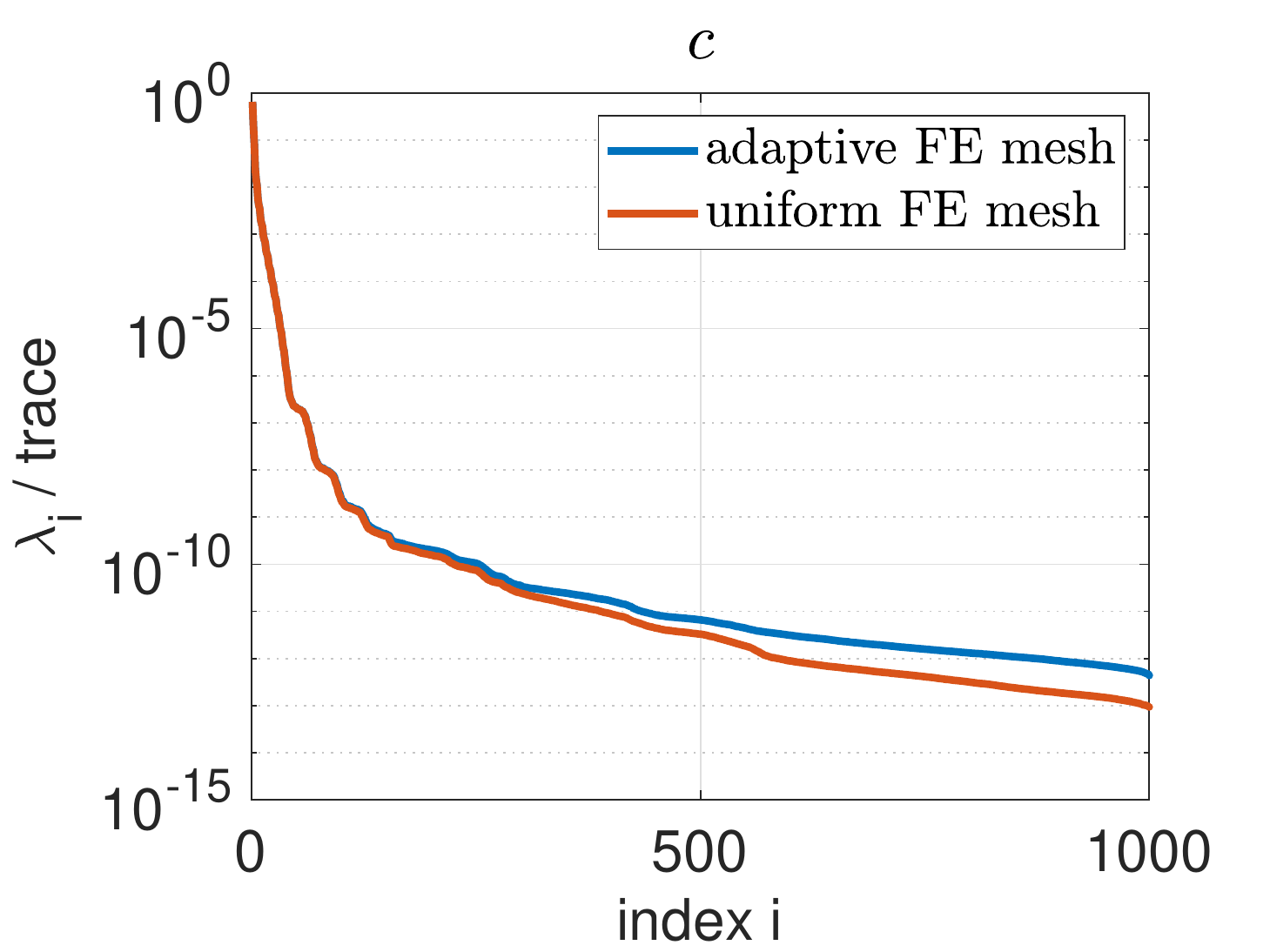} \hfill \includegraphics[width=.46\textwidth]{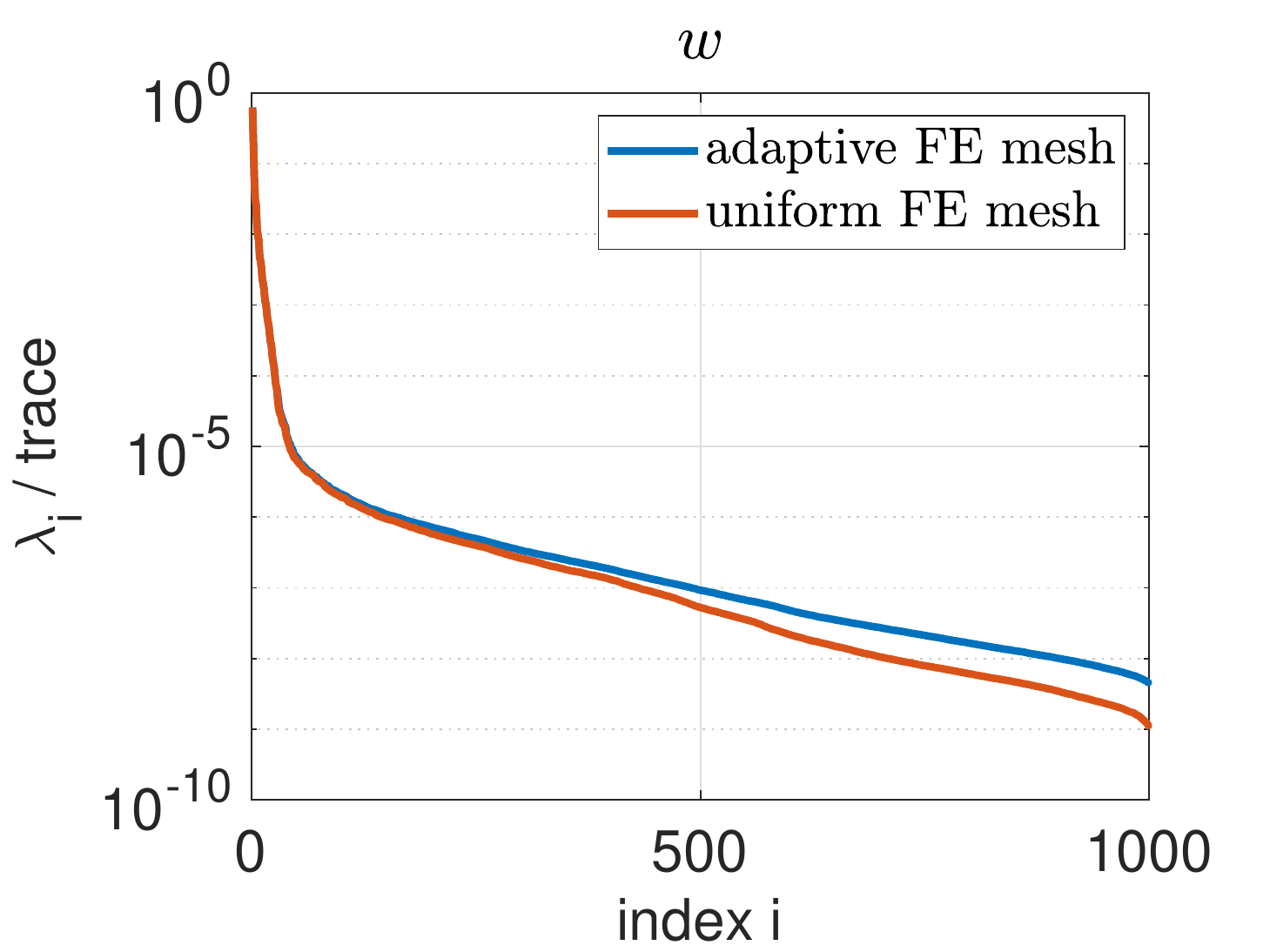}
\caption{\em Run~{\rm\ref{ch-equation}}. Comparison of the normalized eigenvalues for the phase field $c$ (left) and the chemical potential $w$ (right) using an adaptive and a uniform spatial mesh, respectively.}
\label{fig:CH_ev}
\end{figure}

\noindent We use the criterion \eqref{info-content} to determine the basis length $\ell$ which is required to represent a prescribed information content with the respective POD space.
We will choose the POD basis length $\ell_c$ for the phase field $c$ and the number of POD modes $\ell_w$ for the chemical potential, such that
\[
\ell_{\min}=\argmin\,\big\{\mathcal E(\ell):\mathcal E(\ell)> 1 - p \}, \quad  \text{with } \ell = \ell_c \text{ and } \ell_w, \text{ respectively,}
\]
for a given value $p$ representing the loss of information. Alternatively, the POD basis length could be chosen in alignment with the POD projection error \eqref{projection_error} with the expected spatial and/or temporal discretization error, compare e.g. \cite[Theorem 5.1]{GH18}. Let us also refer, e.g., to the recent paper \cite{BMV19}, where different adaptive POD basis extension techniques are discussed. Table~\ref{tab:loss_info} summarizes how to choose $\ell_c$ and $\ell_w$ in order to capture a desired amount of information. Moreover, it tabulates the POD projection error \eqref{projection_error} depending on the POD basis length, where $\lambda_i^c$ and $\lambda_i^w$ denote the eigenvalues for the phase field $c$ and the chemical potential $w$, respectively. The results in Table \ref{tab:loss_info} agree with our expectations: the smaller the loss of information $p$ is, the more POD modes are needed and the smaller is the POD projection error.\hfill$\Diamond$\\
\begin{table}[H]
\centering
\begin{tabular}{  c || c | c || c | c || c | c || c | c }
  $ p  $ & $\ell_c^{\text{ad}} $ &  $\sum_{i>\ell}\lambda_i^c$ & $\ell_w^{\text{ad}} $ &  $\sum_{i>\ell}\lambda_i^w$ &  $\ell_c^{\text{uni}} $ &  $\sum_{i>\ell}\lambda_i^c$ & $\ell_w^{\text{uni}} $ & $\sum_{i>\ell}\lambda_i^w$\\
   \hline
  $10^{-1}$ & 3 &  $2.0 \cdot 10^{-3}$ & 4 &  $156.9 \cdot 10^{0}$ & 3  &  $2.0 \cdot 10^{-3}$ & 4 &  $157.6 \cdot 10^{0}$\\
  $10^{-2}$ & 10 &  $2.1 \cdot 10^{-4}$ & 13 &  $\phantom{1}15.8 \cdot 10^{0}$ & 10 &  $2.1 \cdot 10^{-4}$ & 13 &  $\phantom{1}15.6 \cdot 10^{0}$\\
  $10^{-3}$ & 19 &  $2.5 \cdot 10^{-5}$ & 26 & $\phantom{11}1.8 \cdot 10^{0}$ & 19 &  $2.5 \cdot 10^{-5}$ & 25 &  $\phantom{11}1.8 \cdot 10^{0}$\\  $10^{-4}$ & 29 &  $2.0 \cdot 10^{-6}$ & 211 & $\hspace{2mm}\phantom{11}1.8 \cdot 10^{-1}$ & 28 &  $2.6 \cdot 10^{-6}$ & 160 & $\hspace{2mm}\phantom{11}1.9 \cdot 10^{-1}$\\
  $10^{-5}$ & 37 & $2.5 \cdot 10^{-7}$ & 644 & $\hspace{2mm}\phantom{11}1.1 \cdot 10^{-2}$ & 37 &  $2.4 \cdot 10^{-7}$ & 419 & $\hspace{2mm}\phantom{11}2.5\cdot10^{-2}$\\
 \end{tabular}\hspace{0.5cm}
\vspace{0.4cm} \caption{\em Run~{\rm\ref{ch-equation}}. Number of needed POD bases in order to achieve a loss of information below the tolerance $p$ using adaptive finite element meshes (columns 2-5) and uniform finite element discretization (columns 6-9) and POD projection error.}
\label{tab:loss_info}
\end{table}
\end{run}

\begin{run}[{\cite[Example 6.3]{GH18}}]
\label{r-ad}(Linear heat equation revisited).
\em  
\noindent We again consider Example \ref{slhe} with $c\equiv 0$ . The purpose of this example is to confirm that our POD approach also is applicable in the case of non-nested meshes like it appears in the case of $r$-adaptivity, for example. We set up the matrix $\mathrm K$ for snapshots generated on sequences of non-nested spatial discretizations. This requires the integration over cut elements, see \cite{GH18}. We choose $\Omega=(0,1)\times(0,1)\subset\mathbb R^2$, $[0,T]=[0,1]$, and we apply a uniform temporal discretization with time step size $\Delta t = 0.01$. The analytical solution in the present example is given by
\[
y(t,\bx)=\sin(\pi x_0)\cdot \sin(\pi x_1) \cdot \cos(2 \pi t x_0),
\]
with $\bx = (x_0,x_1)$, source term $f:= y_t-\Delta y$ and the initial condition $g:=y(0,\cdot)$. The initial condition is discretized using piecewise linear and continuous finite elements on a uniform spatial mesh which is shown in Figure~\ref{fig:ex3_meshes} (left). Then, at each time step, the mesh is disturbed by relocating each mesh node according to the assignment
\begin{equation*}
\begin{array}{r c l}
x_0 & \leftarrow & x_0 + \theta \cdot x_0 \cdot (x_0-1) \cdot (\Delta t / 10),\\
x_1 & \leftarrow & x_1 + \theta \cdot 0.5 \cdot x_1 \cdot (x_1-1) \cdot (\Delta t / 10),
\end{array}
\end{equation*}
where $\theta \in \mathbb R_+$ is sufficiently small such that all coordinates of the interior nodes fulfill $0<x_0<1$ and $0<x_1<1$. After relocating the mesh nodes, the heat equation is solved on this mesh for the next time instance. We use Lagrange interpolation to transfer the finite element solution of the previous time step onto the new mesh. The disturbed meshes at $t=0.5$ and $t=1.0$ as well as an overlap of two meshes are shown in Figure \ref{fig:ex3_meshes}. To compute the matrix $\mathrm K$ from \eqref{matrix-k} we have to evaluate the corresponding inner products of the snapshots, where we need to integrate over cut elements.\\
\begin{figure}[H]
\centering
\includegraphics[scale=0.17]{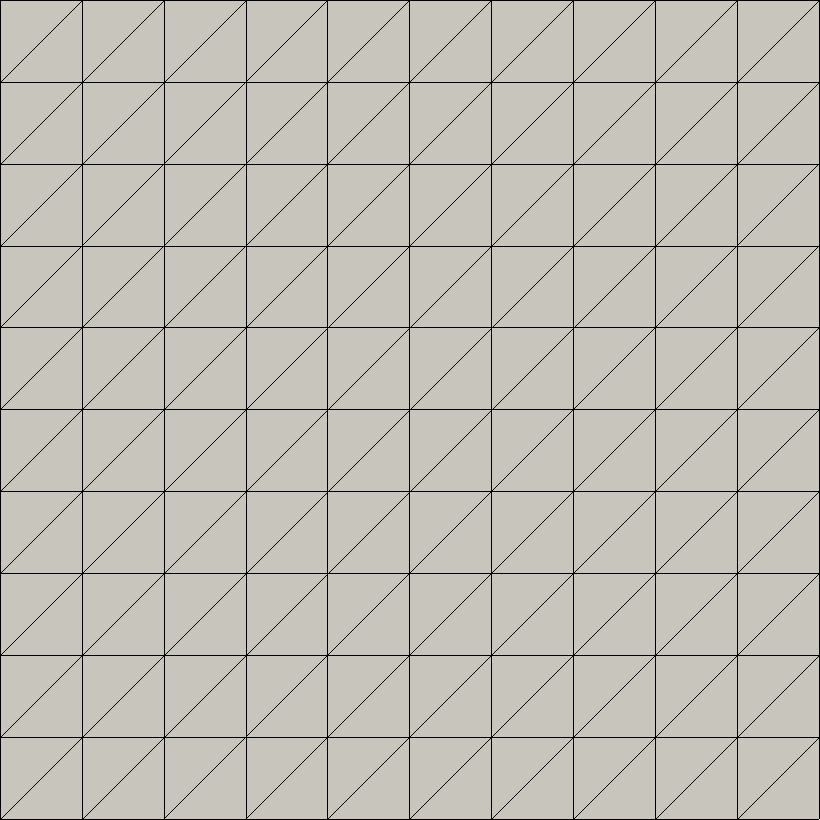} 
\includegraphics[scale=0.17]{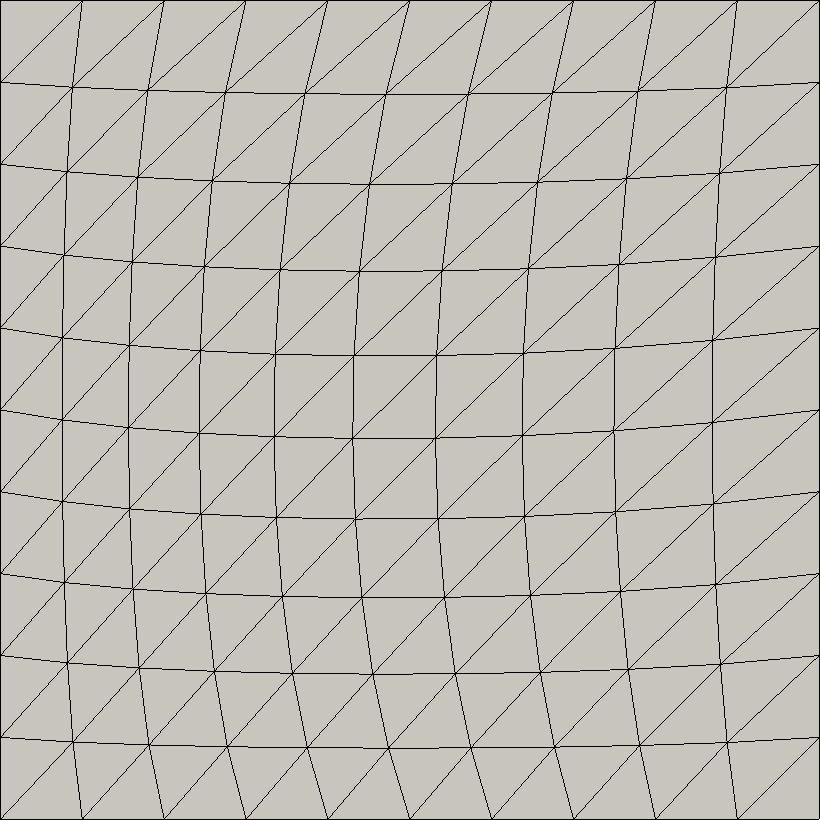} 
\includegraphics[scale=0.17]{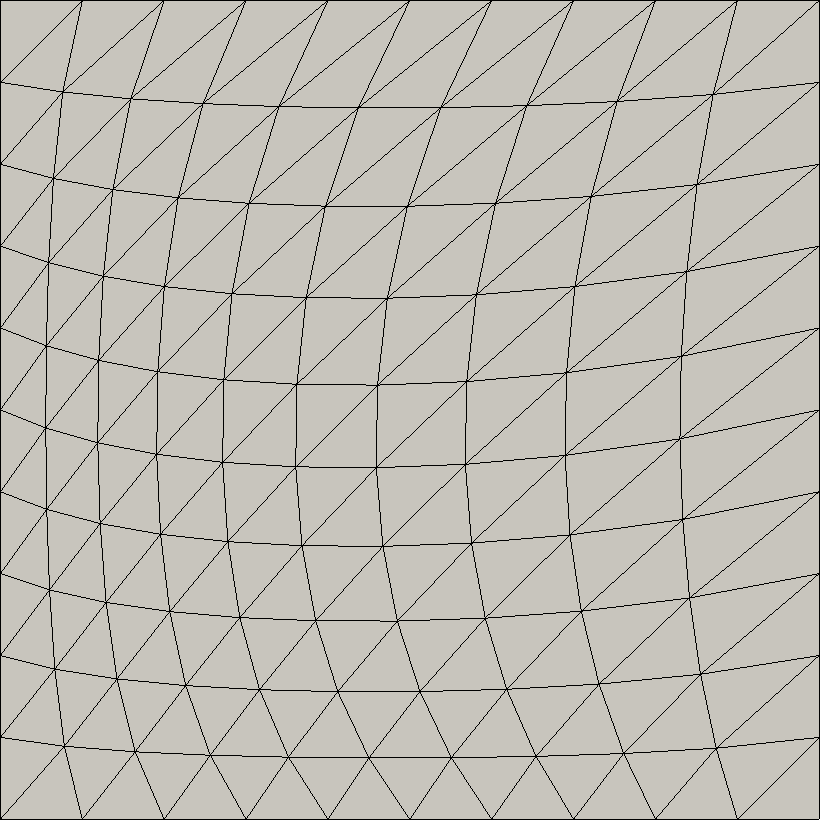} 
\includegraphics[scale=0.17]{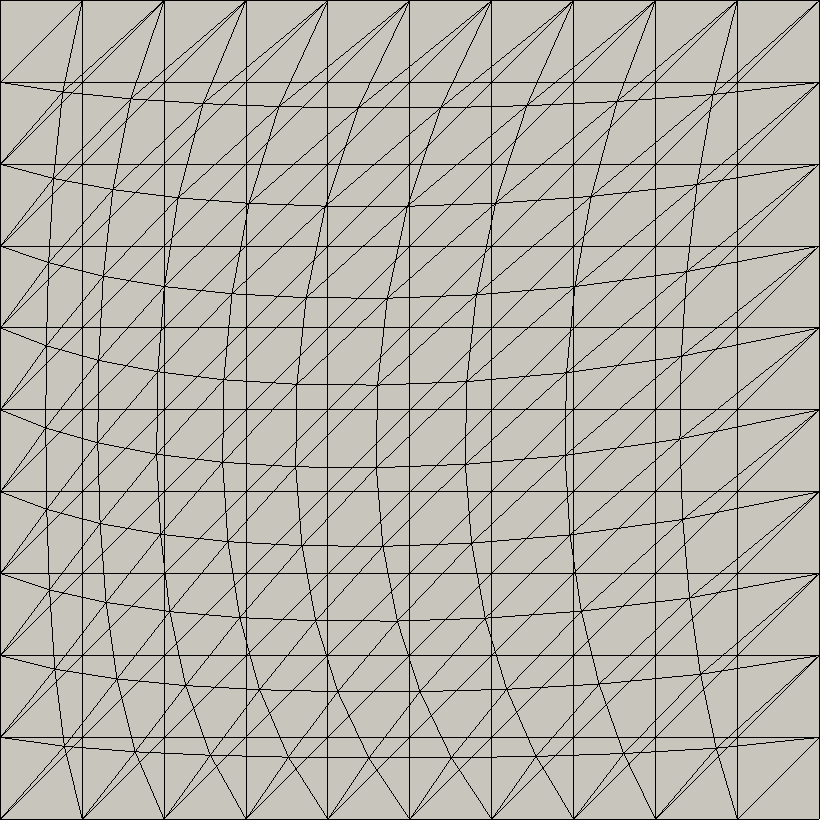} 
\caption{\em Run~{\rm\ref{r-ad}}. Uniform mesh (left), disturbed meshes at $t=0.5$ and $t=1.0$ (middle left, middle right), overlap of the mesh at $t=0$ with the mesh at $t=1.0$ (right). Here, we use $\theta=10$.}
\label{fig:ex3_meshes}
\end{figure}
 
\noindent We compute the eigenvalue decomposition of the matrix representation $\mathrm K$ of the operator $\mathcal K$ (cf. \eqref{matrix-k}) for different values of $\theta$ and compare the results with a uniform mesh (i.e. $\theta = 0$) in Figure \ref{fig:ex3_ev}. We note that the eigenvalues of the disturbed mesh are converging to the eigenvalues of the uniform mesh for $\theta \to 0$. As expected, the eigenvalue spectrum depends only weakly on the underlying mesh given that the mesh size is sufficiently small. 
\begin{figure}[htbp]
\centering
\includegraphics[width=.45\textwidth]{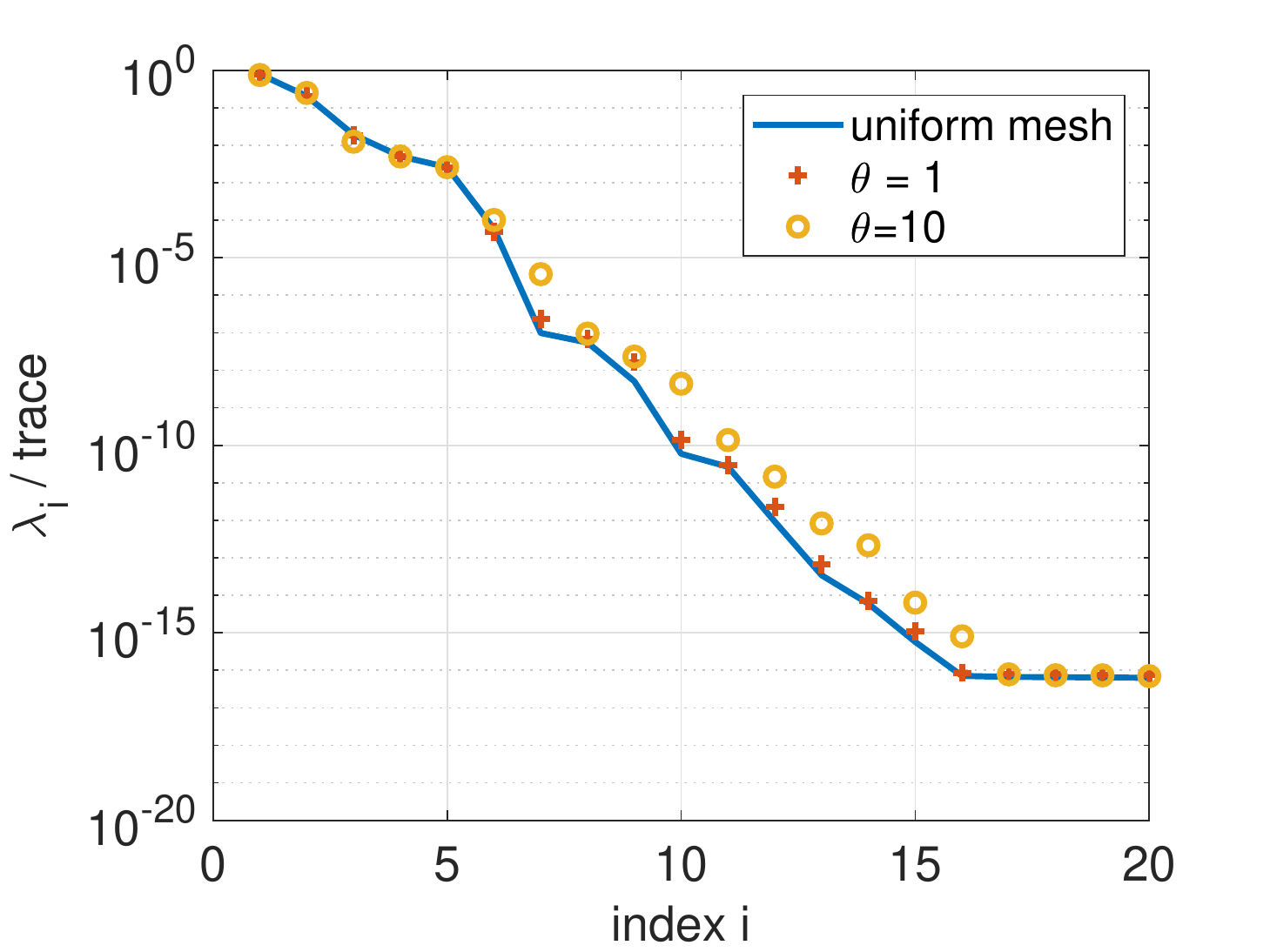} \hfill   \includegraphics[width=.45\textwidth]{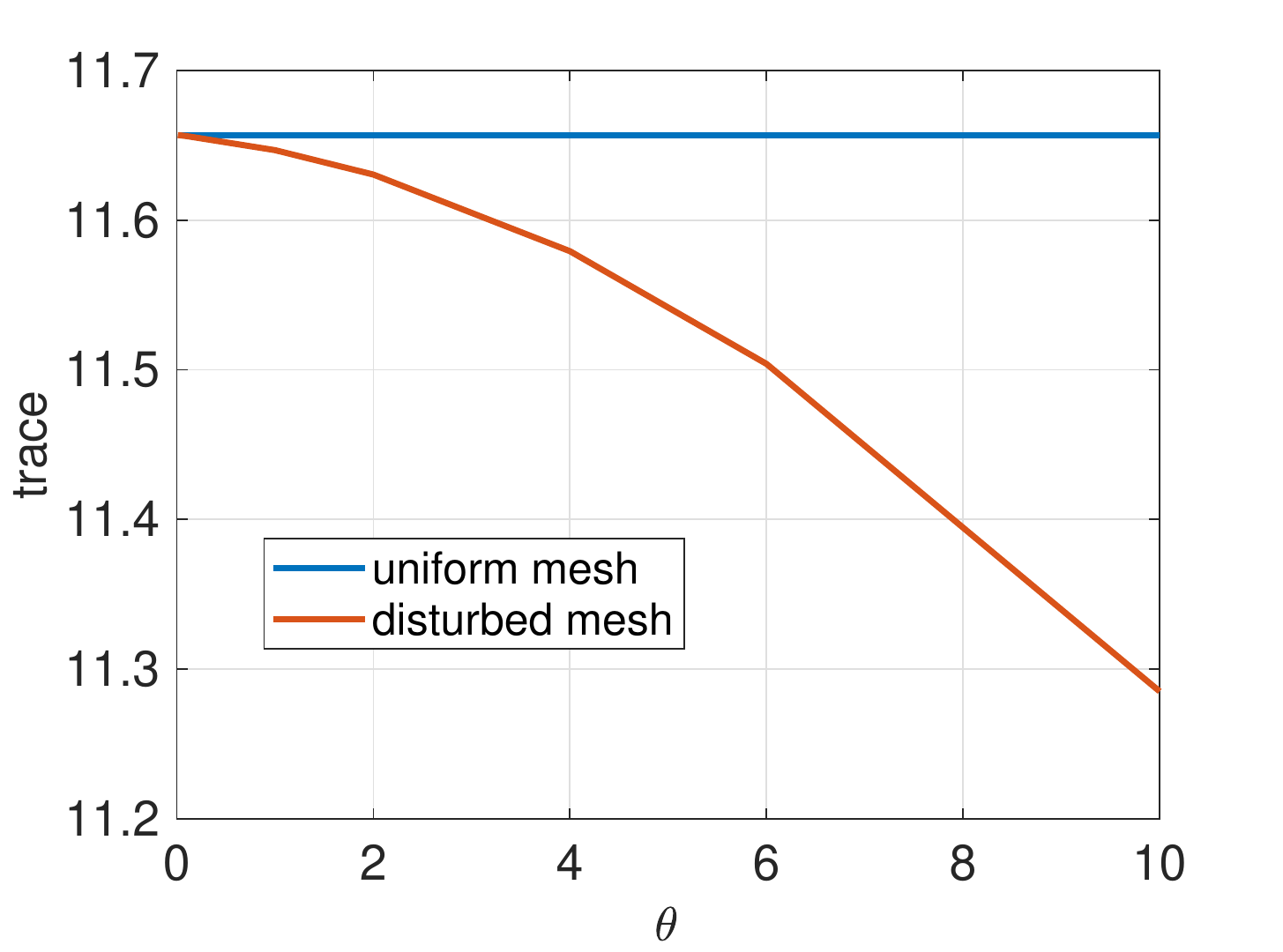} 
\caption{\em Run~{\rm\ref{r-ad}}: Decay of eigenvalues of matrix $K$ with different meshes.}
\label{fig:ex3_ev}
\end{figure}
Concerning the computational complexity of POD with non-nested meshes let us note that solving the heat equation takes 2.1 seconds on the disturbed meshes and 1.8 seconds on the uniform mesh. The computational time needed to compute each entry of the matrix $\mathrm K$ is 0.022 seconds and computing the eigenvalue decomposition for $\mathrm K$ takes 0.0056 seconds. Note that the cut element integration problem for each matrix entry takes a fraction of time required to solve the finite element problem.\hfill$\Diamond$
\end{run}

%%%%%%%%%%%%%%%%%%%%%%%%%%%%%%%%%%%%%%%%%%%%%%%%%%%
\section{The POD Galerkin procedure}
\label{POD-Galerkin}
%%%%%%%%%%%%%%%%%%%%%%%%%%%%%%%%%%%%%%%%%%%%%%%%%%%

Once the POD basis is generated it can be used to set up a POD-Galerkin approximation of the original dynamical system. This is discussed in the present section. In this context we recall that the space spanned by the POD basis is used with a Galerkin method to approximate the original system for e.g. other inputs and/or parameters than those used to generate the snapshots for constructing the POD basis. A typical application is given by PDE-constrained optimization, where the PDE system during the optimization is substituted by POD Galerkin surrogates, see Section~\ref{POD-opt} for more details.

%%%%%%%%%%%%%%%%%%%%%%%%%%%%%%%%%%%%%%%%%%%%%%%%%%%
\subsection{The POD Galerkin procedure}
\label{P2.1}
%%%%%%%%%%%%%%%%%%%%%%%%%%%%%%%%%%%%%%%%%%%%%%%%%%%

Suppose that for given snapshots $y_j^h\in V^{h_j}\subset X$, $1\le j\le n$, we have computed the symmetric matrix
\begin{equation}\label{matrix-k}
\mathrm K=\left(\left(\sqrt{\alpha_i} \sqrt{\alpha_j}\,{\langle y_i^h,y_j^h\rangle}_X\right)\right)_{1\le i,j\le n}\quad\text{with }\mathrm{rank}\, \mathrm K=\mathsf r \le n
\end{equation}
associated to the operator $\mathcal K$ from \eqref{OperatorK} together with its eigensystem. Its $\ell\in\{1,\ldots,\mathsf r\}$ largest eigenvalues are $\{ \lambda_i\}_{i=1}^\ell$ with corresponding eigenvectors $\{ \Phi_i \}_{i=1}^\ell \subset \mathbb R^n$. The POD basis $\{\Psi_i\}_{i=1}^\ell$ is then given by \eqref{R-K}, i.e.,
\[
\Psi_i = \frac{1}{\sqrt{\lambda_i}}\,\mathcal{Y} \Phi_i\quad\text{ for } i=1,\ldots,\ell.
\]
This POD basis is utilized in order to compute a reduced-order model for \eqref{State} along the lines of Section \ref{P1.3.3}, where the space $V^h$ is replaced by the space $V^\ell=\mathrm{span}\,\{\Psi_1,\dots,\Psi_\ell\}\subset V$. More precisely we make the POD Galerkin ansatz
\begin{equation}
\label{PODGalerkin}
y^\ell (t)=\sum\limits_{i=1}^\ell \eta_i (t) \Psi_i=\sum_{i=1}^\ell \eta_i (t)\frac{1}{\sqrt{\lambda_i}}\,\mathcal Y\Phi_i\quad\text{for all } t \in [0,T],
\end{equation}
as an approximation for $y(t)$, with the Fourier coefficients 
\[
\eta_i(t)={\langle y^\ell(t),\Psi_i\rangle}_X=\bigg\langle y^\ell(t),\frac{1}{\sqrt{\lambda_i}}\,\mathcal Y\Phi_i\bigg\rangle_X\quad\text{for }1\le i\le\ell.
\]
Inserting $y^\ell$ into \eqref{State} and choosing $V^\ell \subset V$ as the test space leads to the system
\begin{equation}
\label{ROM_sys_yell}
\left.
\begin{aligned}
\frac{\mathrm d}{\mathrm dt}\,{\langle y^\ell(t),\Psi\rangle}_H+a(y^\ell(t),\Psi)+{\langle \mathcal N(y^\ell(t)),\Psi\rangle}_{V',V}&= {\langle f(t),\Psi \rangle}_{V',V}\\
{\langle y^\ell(0),\Psi\rangle}_H&={\langle y_\circ,\Psi\rangle}_H
\end{aligned}
\right\}
\end{equation}
for all $\Psi\in V^\ell$ and for almost all $t \in (0,T]$.  The system \eqref{ROM_sys_yell} is called POD reduced-order model (POD-ROM). Using the ansatz \eqref{PODGalerkin}, we can write \eqref{ROM_sys_yell} as an $\ell$-dimensional ordinary differential equation system for the POD mode coefficients $\eta(t)=(\eta_i(t))_{1\le i\le\ell}$, $t \in (0,T]$, as follows:
\begin{equation}
\label{ROM_sys}
\left.
\begin{aligned}
\sum_{j=1}^\ell{\dot \eta}_j (t)\,{\langle \Psi_i,\Psi_j\rangle}_H+\sum_{j=1}^\ell\eta_j (t)\,a(\Psi_j,\Psi_i)&={\langle f(t)-\mathcal N(y^\ell(t)),\Psi_i \rangle}_{V',V}\\
\sum_{j=1}^\ell \eta_j(0)\,{\langle\Psi_i,\Psi_j\rangle}_H&={\langle y_\circ,\Psi_i\rangle}_H
\end{aligned}
\right\}
\end{equation}
for $i=1,\ldots,\ell$. Note that $\langle\Psi_i,\Psi_j\rangle_H = \delta_{ij}$ if we choose $X=H$ in the context of Section~\ref{P1.3}. In a next step we rewrite this system using  relation between $\Psi_i$ and $\Phi_i$ given in \eqref{R-K}. This leads to
\begin{equation}
\label{ROM_sys_phi}
\left.
\begin{aligned}
\sum_{j=1}^\ell {\dot \eta}_j (t)\,{\textstyle\frac{{\langle \mathcal Y\Phi_i,\mathcal Y\Phi_j\rangle}_H}{\sqrt{\lambda_i\lambda_j}}}+\sum_{j=1}^\ell\eta_j (t)\,{\textstyle\frac{a(\mathcal Y\Phi_j,\mathcal Y\Phi_i)}{\sqrt{\lambda_i\lambda_j}}}&={\textstyle\frac{{\langle f(t)-\mathcal N(y^\ell(t)), \mathcal Y\Phi_i\rangle}_{V',V}}{\sqrt{\lambda_i}}}\\
&\quad\text{ for } t \in (0,T], \\   
\sum_{j=1}^\ell\eta_j(0)\,{\textstyle\frac{{\langle \mathcal Y\Phi_i, \mathcal Y\Phi_j\rangle}_H}{\sqrt{\lambda_i\lambda_j}}}&={\textstyle\frac{{\langle y_\circ,\mathcal Y\Phi_i\rangle}_H}{\sqrt{\lambda_i}}}
\end{aligned}
\right\}
\end{equation}
for $i=1,\ldots,\ell$. In order to write \eqref{ROM_sys_phi} in a compact matrix-vector form, let us introduce the diagonal matrix
\[
\Uplambda:=\mathrm{diag}\,\bigg(\frac{1}{\sqrt{\lambda_1}},\ldots, \frac{1}{\sqrt{\lambda_\ell}}\bigg)\in \mathbb R^{\ell \times \ell}.
\]
From the first $\ell$ eigenvectors $\{\Phi_i\}_{i=1}^\ell$ of $\mathrm K$ we build the matrix
\[
\Upphi=\big[\Phi_1\,\big|\ldots\big|\,\Phi_\ell]\in \mathbb R^{n\times\ell}.
\]
Then, the system \eqref{ROM_sys_phi} can be written as the system
\begin{equation}
\label{ROM}
\left.
\begin{aligned}
\big(\Uplambda\Upphi^\top \mathrm K \Upphi\Uplambda\big)\dot\eta(t)+\big(\Uplambda\mathrm A^\ell\Uplambda\big)\eta(t)+\Uplambda \mathrm N^\ell(\eta(t))&=\Uplambda\mathrm F^\ell(t)\quad\text{for } t \in (0,T],\\
\Uplambda\Upphi^\top \mathrm K \Upphi\Uplambda\eta(0)&=\Uplambda\eta_\circ
\end{aligned}
\right\}
\end{equation}
for the vector-valued mapping $\eta=(\eta_1,\ldots,\eta_\ell)^\top: [0,T]\to\mathbb R^\ell$, for the nonlinearity $\mathrm N^\ell=(\mathrm N_i^\ell(\cdot))_{1\le i\le\ell}:\mathbb R^\ell\to\mathbb R^\ell$ with
\[
\mathrm N_i^\ell(\mathrm v)=\Big\langle\mathcal N\big({\textstyle\sum_{j=1}^\ell}\mathrm v_j\Psi_j\big),\varphi_i\Big\rangle_{V',V}=\Big\langle\mathcal N\big({\textstyle\sum_{j=1}^\ell}\mathrm v_j\mathcal Y\Phi_j/\sqrt{\lambda_j}\big),\varphi_i\Big\rangle_{V',V}
\]
and for the stiffness matrix $\mathrm A^\ell=((\mathrm A^\ell_{ij}))\in\mathbb R^{\ell\times\ell}$ given as
\[
\mathrm A^\ell_{ij}=a(\mathcal Y\Phi_j,\mathcal Y\Phi_i)\quad\text{for }1\le i,j\le\ell.
\]
Note that the right hand side $\mathrm F^\ell(t)=(\mathrm F_i^\ell(t))_{1\le i\le\ell}$ and the initial condition $\eta_\circ=(\eta_{\circ i})_{1\le i\le\ell}$ are given by
\[
\mathrm F_i^\ell(t)={\langle f(t),\mathcal Y\Phi_i\rangle}_{V',V}={\langle\mathcal Y^*f(t),\Phi_i\rangle}_{\mathbb R^n},\quad t\in[0,T]\text{ a.e.},
\]
and 
\[
\eta_{\circ i}={\langle y_\circ,\mathcal Y\Phi_i\rangle}_H={\langle\mathcal Y^*y_\circ,\Phi_i\rangle}_{\mathbb R^n},
\]
for $i=1,\ldots,\ell$, respectively. Their calculation can be done explicitly for any arbitrary finite element discretization. For a given function $w \in V$ (for example $w=f(t)$ or $w=y_\circ$) with finite element discretization $w=\sum_{i=1}^{m_w} \mathrm w_i \chi_i$, nodal basis $\{\chi_i\}_{i=1}^{m_w} \subset V$ and appropriate mode coefficients $\{\mathrm w_i\}_{i=1}^{m_w}$ we can compute 
\[
\big(\mathcal Y^*w\big)_j={\langle w,y_j\rangle}_X=\bigg\langle\sum_{i=1}^{m_w} \mathrm w_i\chi_i,\sum_{k=1}^{m_j}\mathrm y_k^j\varphi_k^j\bigg\rangle_X=\sum_{i=1}^{m_w}\sum_{k=1}^{m_j}\mathrm w_i \mathrm y_k^j\,{\langle\chi_i,\varphi_k^j\rangle}_X
\]
for $j=1,\ldots,n$ where $ y_j^h=\sum_{k=1}^{m_j}\mathrm y_k^j\varphi_k^j\in V^{h_j}$ denotes the $j$-th snapshot. Again, for any $i =1,\ldots,m_w$ and $k=1,\ldots,m_j$, the computation of the inner product $\langle\chi_i,\varphi_k^j\rangle_X$ can be done explicitly.

Obviously, for linear evolution equations the POD reduced-order model \eqref{ROM} can be set up and solved using snapshots with arbitrary finite element discretizations. The computation of the nonlinear component $\mathrm N^\ell(\eta(t))$ needs particular attention. In Section~\ref{P2.3} we discuss the options to treat the nonlinearity.
 
%%%%%%%%%%%%%%%%%%%%%%%%%%%%%%%%%%%%%%%%%%%%%%%%%%%
\subsection{Time-discrete reduced-order model}
\label{P2.2}
%%%%%%%%%%%%%%%%%%%%%%%%%%%%%%%%%%%%%%%%%%%%%%%%%%%

In order to solve the reduced-order system \eqref{ROM_sys_yell} numerically, we apply the implicit Euler method for time discretization and use for simplicity the same temporal grid $\{t_j\}_{j=1}^n$ as for the snapshots. It is also possible to use a different time grid, cf. \cite{KV02}. The time-discrete reduced-order model reads 
\begin{equation}
\label{ROM_timediscrete}
\left.
\begin{aligned}
\hspace{-3mm}\frac{{\langle y_j^\ell-y_{j-1}^\ell,\Psi\rangle}_H}{\Delta t_j}+a(y_j^\ell,\Psi)+{\langle \mathcal{N}(y_j^\ell),\Psi\rangle}_{V',V}&=\int_{t_{j-1}}^{t_j}\frac{{\langle f(\tau),\Psi\rangle}_{V',V}}{\Delta t_j}\,\mathrm d\tau,\\
{\langle y_1^\ell,\Psi\rangle}_H&={\langle y_\circ,\Psi\rangle}_H
\end{aligned}
\right\}
\end{equation}
for all $\Psi\in V^\ell$ and $j=2,\ldots,n$. Equivalently the following system holds for the coefficient vector $\eta(t)\in\mathbb R^\ell$ (cf. \eqref{ROM}):
\begin{equation}
\label{ROM_infPOD_timediscrete}
\left.
\begin{aligned}
\hspace{-4mm}\big(\Uplambda\Upphi^\top \mathrm K \Upphi\Uplambda\big)\bigg(\frac{\eta^j-\eta^{j-1}}{\Delta t_j}\bigg)+\big(\Uplambda\mathrm A^\ell\Uplambda\big)\eta^j+\Uplambda\mathrm N^\ell(\eta^j)&=\Uplambda\mathrm F^\ell_j,~j=2,\ldots,n\\
\Uplambda\Upphi^\top \mathrm K \Upphi\Uplambda\eta^1&=\Uplambda\eta_\circ
\end{aligned}
\right\}
\end{equation}
with the inhomogeneity $\mathrm F^\ell_j=(\mathrm F_{ji}^\ell)_{1\le i\le\ell}$, $j=2,\ldots,n$, given as
\[
\mathrm F_{ji}^\ell=\int_{t_{j-1}}^{t_j}\frac{{\langle f(\tau),\mathcal Y\Phi_i\rangle}_{V',V}}{\Delta t_j}\,\mathrm d\tau=\int_{t_{j-1}}^{t_j}\frac{{\langle\mathcal Y^*f(\tau),\Phi_i\rangle}_{\mathbb R^n}}{\Delta t_j}\,\mathrm d\tau.
\]
 
%%%%%%%%%%%%%%%%%%%%%%%%%%%%%%%%%%%%%%%%%%%%%%%%%%%
\subsection{Discussion of the computation of the nonlinear term}
\label{P2.3}
%%%%%%%%%%%%%%%%%%%%%%%%%%%%%%%%%%%%%%%%%%%%%%%%%%%

Let us now consider the computation of the nonlinear term $\Uplambda\mathrm N^\ell(\eta^j)\in\mathbb R^\ell$ of the POD-ROM \eqref{ROM}. It holds true
\[
\big(\Uplambda\mathrm N^\ell(\eta^j)\big)_k={\langle \mathcal N(y^\ell(t)),\Psi_k\rangle}_{V',V}=\Big\langle \mathcal N\big({\textstyle\sum_{i=1}^\ell \eta_i(t) \Psi_i}\big),\Psi_k\Big\rangle_{V',V}
\]
for $k=1,\ldots,\ell$. It is well known that the evaluation of nonlinearities in the reduced-order modeling context is computationally expensive. To make this clear, let us assume, we are given a uniform finite element discretization with $m$ degrees of freedom. Then, in the fully discrete setting, the nonlinear term has the form 
\[
\Uppsi^\top W\mathrm N^h(\Uppsi\eta(t))\in\mathbb R^\ell,\quad t\in[0,T]\text{ a.e.},
\]
where $\Uppsi=[\Psi_1\,|\ldots|\,\Psi_\ell ] \in \mathbb R^{m\times\ell}$ is the matrix in which the POD modes are stored columnwise and $W\in\mathbb R^{m\times m}$ is a weighting matrix related to the utilized inner product (cf. \eqref{Wmatrix}). Hence, the treatment of the nonlinearity requires the expansion of $\Uppsi\eta(t)\in\mathbb R^m$ in the full space for $t\in[0,T]$ a.e. Then the 
nonlinearity can be evaluated and finally the result is projected back to the POD space. Obviously, this means that the reduced-order model is not fully independent of the high-order dimension $m$ and efficient simulation cannot be guaranteed. Therefore, it is convenient to seek for hyper reduction, i.e., for a treatment of the nonlinearity, where the model evaluation cost is related to the low dimension $\ell$. Common choices are empirical interpolation methods like, e.g., EIM (\cite{BMNP04}), DEIM (\cite{CS10}), and QDEIM (\cite{DG15}). Another option is dynamic mode decomposition for nonlinear model order reduction, see e.g. \cite{AK16}.  Furthermore, in \cite{Wan15} nonlinear model reduction is realized by replacing the nonlinear term by its interpolation in the finite element space. An alternative approach for the treatment of the nonlinearity is missing point estimation \cite{AWWB08}, or best points interpolation \cite{NPP08}.\\
Most of these methods need a common reference mesh for the computations. To overcome this restriction we propose different paths which allow for more general discrete settings like $r$-adaptivity discussed in Run \ref{r-ad}.

One option is to use EIM \cite{BMNP04}. Alternatively, we can linearize and project the nonlinearity onto the POD space. For this approach, let us consider the linear reduced-order system for a fixed given state $\bar y$, which takes the form
\begin{equation}
\label{P_linearized}
\left.
\begin{aligned}
\frac{\mathrm d}{\mathrm dt}\,{\langle y^\ell(t),\Psi\rangle}_H+a(y^\ell(t),\Psi)+{\langle\mathcal N(\bar y(t)),\Psi\rangle}_{V',V}&={\langle f(t),\Psi\rangle}_{V',V},\\   
{\langle y^\ell(0),\Psi\rangle}_H&={\langle y_\circ,\Psi\rangle}_H   
\end{aligned}
\right\}
\end{equation}
for all $\Psi\in V^\ell$ and for almost all $t \in (0,T]$. The linear evolution problem \eqref{P_linearized} can be set up and solved explicitly without spatial interpolation. In the numerical examples in Section 6, we take the finite element solution as given state in each time step, i.e., $\bar y(t_j)=y_j$ for $j=2,\ldots,n$.

Furthermore, the linearization of the reduced-order model \eqref{ROM_sys_yell} can be considered:
\begin{equation}
\label{P_linearized2}
\left.
\begin{aligned}
&\frac{\mathrm d}{\mathrm dt}\,{\langle y^\ell(t),\Psi\rangle}_H+a(y^\ell(t),\Psi)+{\langle \mathcal N'(\bar y(t))y^\ell(t),\Psi \rangle}_{V',V},\\
&\hspace{30mm}={\langle f(t)-\mathcal N(\bar y(t))+\mathcal N'(\bar y(t)))\bar y(t),\Psi\rangle}_{V',V}\\
&\hspace{13mm}{\langle y^\ell(0),\Psi\rangle}_H={\langle y_\circ,\Psi\rangle}_H,   
\end{aligned}
\right\}
\end{equation}
for all $\Psi\in V^\ell$ and for almost all $t \in (0,T]$, where $\mathcal N'$ denotes the Fr\'echet derivative of the nonlinear operator $\mathcal N$. This linearized problem is of interest e.g. in the context of optimal control, where it occurs in each iteration level within sequential quadratic programming (SQP) methods; see \cite{HPUU09}, for example. Choosing the finite element solution as given state in each time instance and using \eqref{R-K} leads to
\begin{align*}
&{\langle \mathcal N(y_j),\Psi_i\rangle}_{V',V}=\frac{1}{\sqrt{\lambda_i}}\sum_{k=1}^n\sqrt{\alpha_k}(\Phi_i)_k\,{\langle \mathcal N(y_j),y_k\rangle}_{V',V},\\
&{\langle \mathcal N'(y_j)y^\ell(t_j),\Psi_i\rangle}_{V',V}=\bigg\langle\mathcal N'(y_j)\Big(\sum_{k=1}^\ell \eta_k(t_j)\Psi_k\big),\Psi_i \bigg\rangle_{V',V},\\
&\qquad=\sum_{k=1}^\ell \eta_k(t_j) \frac{1}{\sqrt{\lambda_k\lambda_i}} \sum_{\nu=1}^n \sum_{\mu=1}^n \sqrt{\alpha_\nu\alpha_\mu}  (\Phi_k)_\nu (\Phi_i)_\mu\,{\langle \mathcal N'(y_j)y_\nu,y_\mu \rangle}_{V',V}, \\
&{\langle\mathcal N'(y_j)y_j,\Psi_i\rangle}_{V',V}=\frac{1}{\sqrt{\lambda_i}} \sum_{k=0}^n \sqrt{\alpha_k}(\Phi_i)_k\,{\langle \mathcal N'(y_j)y_j,y_k\rangle}_{V',V}
\end{align*}
for $j=2,\ldots,n$ and $i=1,\ldots,\ell$. Finally, we approximate the nonlinearity $\Uplambda\mathrm N^\ell(\eta^j)\in\mathbb R^\ell$ in \eqref{ROM_infPOD_timediscrete} by
\[
\big(\Uplambda\mathrm N^\ell(\eta^j)\big)_i\approx{\langle \mathcal N(y_j)+\mathcal N'(y_j)(y^\ell(t_j)- y_j),\Psi_i\rangle}_{V',V}
\]
for $j=2,\ldots,n$ and $i=1,\ldots,\ell$, which can be written as
\[
\Uplambda\mathrm N^\ell(\eta^j)\approx\Uplambda\Upphi^\top\mathsf{N}^j + \Uplambda\Upphi^\top \mathtt{N}_y^j\Upphi \Uplambda\eta^j-\Uplambda\Upphi^\top \mathsf{N}_y^j,
\]
where 
\[
\mathsf N^j=
\begin{pmatrix}
{\langle\mathcal N(y_j),\sqrt{\alpha_1} y_1\rangle}_{V',V}\\
\vdots\\
{\langle\mathcal N(y_j),\sqrt{\alpha_n} y_n\rangle}_{V',V}
\end{pmatrix} \in \mathbb R^n,~\mathsf N_y^j=  \begin{pmatrix}
{\langle\mathcal N'(y_j) y_j,\sqrt{\alpha_1}y_1\rangle}_{V',V}\\
\vdots\\
{\langle \mathcal N'(y_j) y_j,\sqrt{\alpha_n}y_n\rangle}_{V',V}\\
\end{pmatrix} \in \mathbb R^n,
\]
and with $\tilde y_j= \sqrt{\alpha_j} y_j$, $j=1,\ldots,n$,
\[
\mathtt N_y^j=
\begin{pmatrix}
{\langle \mathcal N'(y_j)\tilde y_1, \tilde y_1\rangle}_{V',V}&\ldots&{\langle\mathcal N'(y_j)\tilde y_n,\tilde y_1\rangle}_{V',V}\\
\vdots&&\vdots\\
{\langle\mathcal N'(y_j)\tilde y_1,\tilde y_1\rangle}_{V',V}&\ldots&{\langle\mathcal N'(y_j)\tilde y_n,\tilde y_n\rangle}_{V',V}
\end{pmatrix}
\in\mathbb R^{n\times n}.
\]
For weakly nonlinear systems this approximation may be sufficient, depending on the problem and its goal. A great advantage of linearizing the semilinear partial differential equation is that only linear equations need to be solved which leads to a further speedup, see Table~\ref{tab:CH_speedups}. However, if a more precise approximation is desired or necessary, we can think of approximations including higher order terms, like quadratic approximation, see, e.g., \cite{Che99} and \cite{RW03}, or Taylor expansions, see, e.g., \cite{Phi00,Phi03} and \cite{FZCZF04}. Nevertheless, the efficiency of higher order approximations is limited due to 
 growing memory and computational costs.
 
\medskip

%%%%%%%%%%%%%%%%%%%%%%%%%%%%%%%%%%%%%%%%%%%%%%%%%%%
\subsection{Expressing the POD solution in the full spatial domain}
\label{P2.4}
%%%%%%%%%%%%%%%%%%%%%%%%%%%%%%%%%%%%%%%%%%%%%%%%%%%

Having determined the solution $\eta(t)$ to \eqref{ROM}, we can set up the reduced solution $y^\ell(t)$ in a continuous framework:
\begin{equation}
\label{solROM}
y^\ell(t)=\sum_{i=1}^\ell\eta_i(t)\bigg(\frac{1}{\sqrt{\lambda_i}} \sum_{j=1}^n\sqrt{\alpha_j}(\Phi_i)_j  y_j\bigg).
\end{equation}
Now, let us turn to the fully discrete formulation of \eqref{solROM}. For a time-discrete setting, we introduce for simplicity the same temporal grid $\{t_j\}_{j=1}^n$ as for the snapshots. The snapshots \eqref{yFE} admit the expansion
\[
y_j = \sum_{i=1}^{m_j} \mathrm y_j^i\varphi_i^j\quad\text{for }j=1,\ldots,n.
\]
Let $\{Q_r^j\}_{r=1}^{l_j}$ denote an arbitrary set of grid points for the reduced system at time level $t_j$. The fully discrete POD solution can be computed by evaluation:
\begin{equation}
\label{yPOD}
y^\ell(t_j,Q_r^j)=\sum_{i=1}^\ell\eta_i(t_j)\bigg(\frac{1}{\sqrt{\lambda_i}}\sum_{\nu=1}^n\sqrt{\alpha_\nu}(\Phi_i)_\nu\Big(\sum_{k=1}^{m_\nu}\mathrm y_k^\nu\varphi_k^\nu(Q_r^j)\Big)\bigg) 
\end{equation}
for $r = 1,\ldots,l_j$ and $j=1,\ldots, n$. This allows us to use any grid for expressing the POD solution in the full spatial domain. For example, we can use the same nodes at time level $j$ for the POD simulation as we have used for the snapshots, i.e., for $j=1,\ldots,n$ it holds $l_j=m_j$ and $Q_r^j=P_k^j$ for all $r,k=1,\ldots,m_j$. Another option can be to choose 
\[
\big\{Q_r^j\big\}_{r=1}^{l_j}=\bigcup_{j=1}^n\bigcup_{k=1}^{m_j}\big\{P_k^j\big\}\quad\text{for }j=1,\ldots,n,
\]
i.e., the common finest grid. Obviously, a special and probably the easiest case concerning the implementation is to choose snapshots which are expressed with respect to the same finite element basis functions and utilize the common finest grid for the simulation of the reduced-order system, which is proposed by \cite{URL16}. After expressing the adaptively sampled snapshots with respect to a common finite element space, the subsequent steps coincide with the common approach of taking snapshots which are generated without adaptivity. Then, expression \eqref{yPOD} simplifies to
\begin{equation}
\label{POD_full}
y^\ell(t_j, P_r)=\sum_{i=1}^\ell \eta_i(t_j)\bigg( \frac{1}{\sqrt{\lambda_i}} \sum_{\nu=1}^n\sqrt{\alpha_\nu}(\Phi_i)_\nu\mathrm y_\nu\bigg)\quad\text{for }j=1,\ldots,n,
\end{equation}
where $\{P_r\}_{r=1}^m$ are the nodes of the common finite element space.

\medskip

\begin{run}[{\cite[Example 6.1]{GH18}}]
\label{lhe-n-II}
\em
Let us revisit Run~\ref{lhe-n} and consider its POD Galerkin solutions. The POD solutions for $\ell = 10$ and $\ell = 50$ POD basis functions using spatial adaptive snapshots which are interpolated onto the finest mesh are shown in Figure~\ref{fig:Heat_POD}. As expected, the more POD basis functions we use (until stagnation of the corresponding eigenvalues), the less oscillations appear in the POD solution and the better the approximation is. 
\begin{figure}[H]
\centering
\includegraphics[width=.3\textwidth]{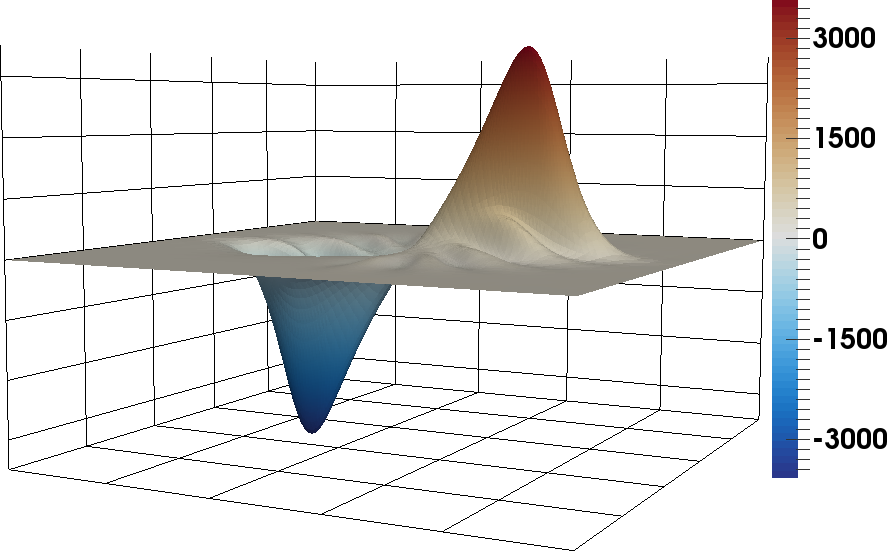} \hfill \includegraphics[width=.3\textwidth]{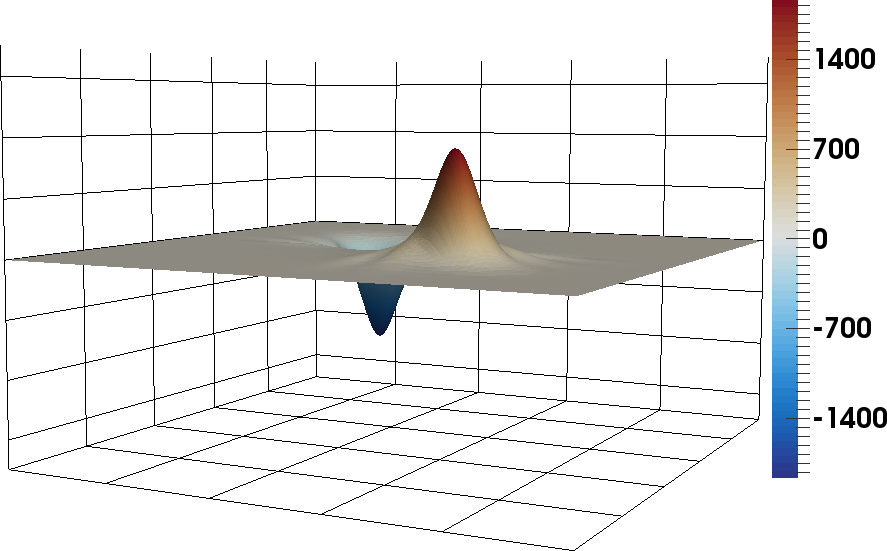} \hfill \includegraphics[width=.3\textwidth]{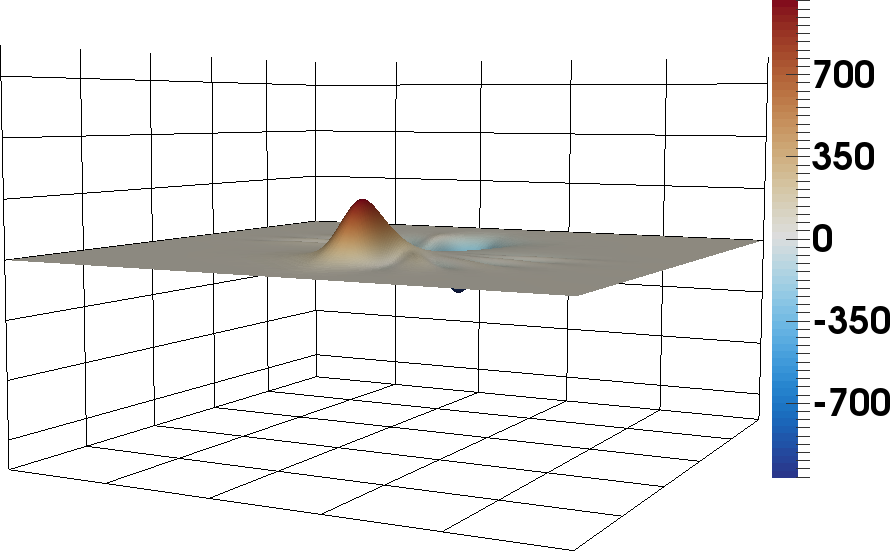}\\
\includegraphics[width=.3\textwidth]{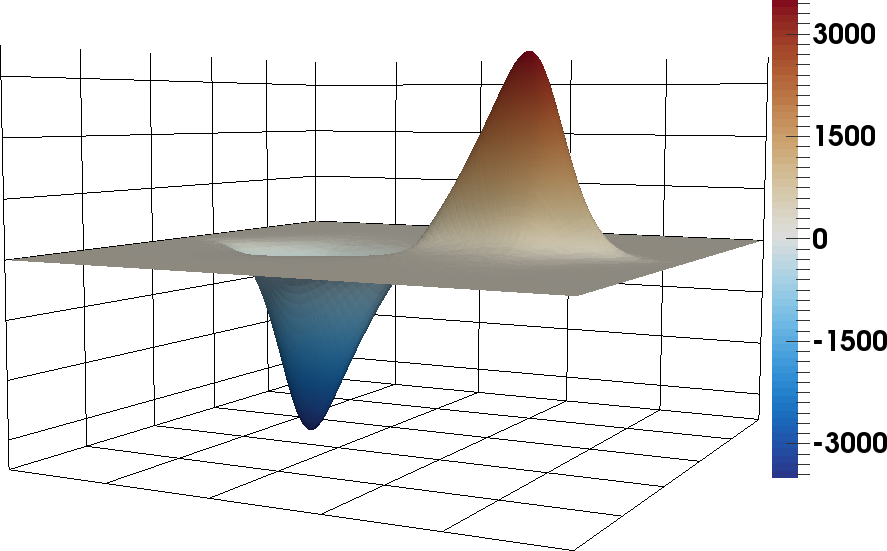} \hfill \includegraphics[width=.3\textwidth]{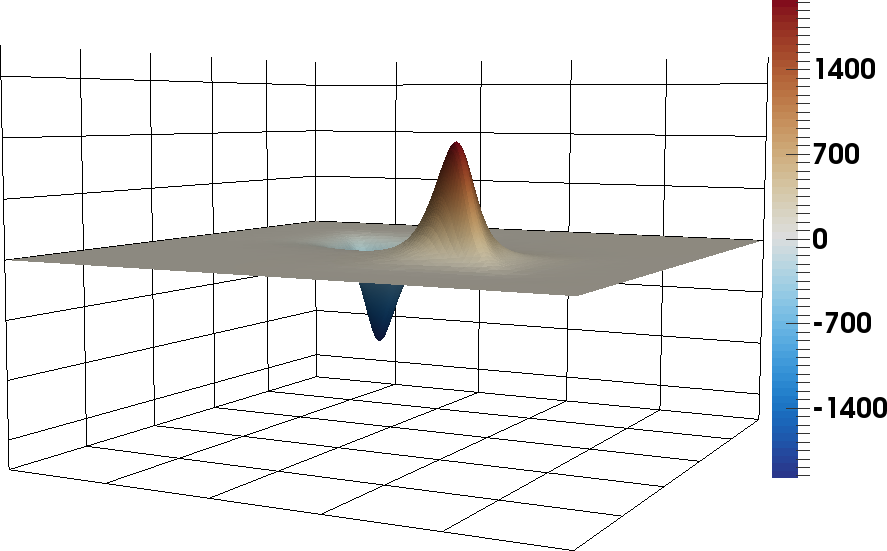}  \hfill \includegraphics[width=.3\textwidth]{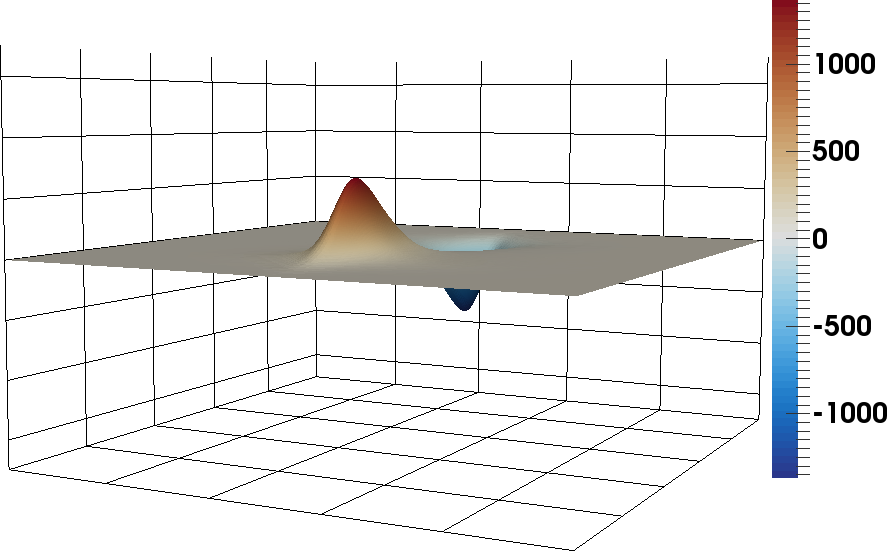}    
\caption{\em Run~{\rm\ref{lhe-n-II}}. Surface plot of the POD solution using $\ell = 10$ (top) and $\ell = 50$ (bottom) POD basis functions at $t=t_1$ (left), $t=T/2$ (middle) and $t=T$ (right).}
\label{fig:Heat_POD}
\end{figure}
Table~\ref{tab:errors_heat} compares the approximation quality in the relative $L^2(0,T;L^2(\Omega))$-norm of the POD solution using adaptively generated snapshots which are interpolated onto the finest mesh with snapshots of uniform spatial discretization depending on different POD basis lengths. Then, for $\ell = 20$ we obtain a relative $L^2(0,T;L^2(\Omega))$-error between the POD solution and the finite element solution of size $\varepsilon_{\text{FE}}^\text{ad} = 3.08 \cdot 10^{-2}$, and a relative $L^2(0,T;L^2(\Omega))$-error between the POD solution and the true solution of size $\varepsilon_{\text{true}}^\text{ad} = 2.17 \cdot 10^{-2}$.\\
\begin{table}[H]
\centering
\begin{tabular}{c|c|c|c|c}
$\ell$&$\varepsilon_{\text{FE}}^{\text{ad}}$&$\varepsilon_{\text{FE}}^{\text{uni}}$&$\varepsilon_{\text{true}}^{\text{ad}}$&$\varepsilon_{\text{true}}^{\text{uni}}$\\\hline
 1  & \hspace{-0.2cm}$1.30 \cdot 10^{0}$ & \hspace{-0.3cm} $1.30 \cdot 10^{0}$ & \hspace{-0.3cm}  $1.28 \cdot 10^{0}$  & \hspace{-0.3cm}  $1.30 \cdot 10^{0}$ \\
 3  & $7.49 \cdot 10^{-1}$  & $7.58 \cdot 10^{-1}$ &  $7.46 \cdot 10^{-1}$  &  $7.60 \cdot 10^{-1}$\\
 5   & $4.39 \cdot 10^{-1}$  & $4.45 \cdot 10^{-1}$ &  $4.39 \cdot 10^{-1}$  & $4.46 \cdot 10^{-1}$ \\
 10   & $1.37 \cdot 10^{-1}$   & $1.37 \cdot 10^{-1}$  & $1.36 \cdot 10^{-1}$  &  $1.38 \cdot 10^{-1}$ \\
 20  & $3.08 \cdot 10^{-2}$  & $1.56 \cdot 10^{-2}$  & $2.17 \cdot 10^{-2}$  & $1.60 \cdot 10^{-2}$   \\
 30  & $2.59 \cdot 10^{-2}$  &  $2.04 \cdot 10^{-3}$ & $1.49   \cdot 10^{-2}$  &  $3.00 \cdot 10^{-3}$ \\
 50  & $2.63 \cdot 10^{-2}$   &  $5.67 \cdot 10^{-5}$ & $1.41 \cdot 10^{-2}$  &    $ 2.07 \cdot 10^{-3}$  \\
 100  & $2.61 \cdot 10^{-2}$  &  $6.48 \cdot 10^{-8}$  & $1.40 \cdot 10^{-2}$  & $2.06 \cdot 10^{-3}$\\
 150  & $2.61 \cdot 10^{-2}$  &  $8.13 \cdot 10^{-7}$  & $1.39 \cdot 10^{-2}$   &   $2.07 \cdot 10^{-3}$\\
\end{tabular}
\vspace{0.4cm}
\caption{\em Run~{\rm\ref{lhe-n-II}}. Relative $L^2(0,T;L^2(\Omega))$-error between the POD solution and the finite element solution (columns 2-3) and the true solution (columns 4-5), respectively, using adaptive finite element snapshots which are interpolated onto the finest mesh and using a uniform mesh.}
\label{tab:errors_heat}
\end{table}
\begin{table}[htbp]
\centering
\begin{tabular}{c|c|c|c|c}
$\ell$&$\varepsilon_{\text{FE}}^{\text{ad}}$&$\varepsilon_{\text{FE}}^{\text{uni}}$&$\varepsilon_{\text{true}}^{\text{ad}}$&$\varepsilon_{\text{true}}^{\text{uni}}$\\\hline
 1  & \hspace{-0.2cm}$1.46 \cdot 10^{0}$ & \hspace{-0.3cm} $1.46 \cdot 10^{0}$ & \hspace{-0.3cm}  $1.46 \cdot 10^{0}$  & \hspace{-0.3cm}  $1.47 \cdot 10^{0}$ \\
 3  & \hspace{-0.3cm} $1.21 \cdot 10^{0}$  & \hspace{-0.3cm} $1.22 \cdot 10^{0}$ &  \hspace{-0.3cm} $1.22 \cdot 10^{0}$  &  \hspace{-0.3cm} $1.22 \cdot 10^{0}$\\
 5   & $9.39 \cdot 10^{-1}$  & $9.45 \cdot 10^{-1}$ &  $9.47 \cdot 10^{-1}$  & $9.51 \cdot 10^{-1}$ \\
 10   & $4.22 \cdot 10^{-1}$   & $4.25 \cdot 10^{-1}$  & $4.33 \cdot 10^{-1}$  &  $4.31 \cdot 10^{-1}$ \\
 20  & $7.76 \cdot 10^{-2}$  & $7.27 \cdot 10^{-2}$  & $1.02 \cdot 10^{-1}$  & $8.19 \cdot 10^{-2}$   \\
 30  & $2.92 \cdot 10^{-2}$  &  $1.22 \cdot 10^{-2}$ & $7.26   \cdot 10^{-2}$  &  $3.52 \cdot 10^{-2}$ \\
 50  & $2.61 \cdot 10^{-2}$   &  $4.74 \cdot 10^{-4}$ & $7.05 \cdot 10^{-2}$  &    $ 3.27 \cdot 10^{-2}$  \\
 100  & $2.79 \cdot 10^{-2}$  &  $4.78 \cdot 10^{-7}$  & $6.94 \cdot 10^{-2}$  & $3.27 \cdot 10^{-2}$\\
 150  & $2.93 \cdot 10^{-2}$  &  $2.84 \cdot 10^{-7}$  & $6.87 \cdot 10^{-2}$   &   $3.27 \cdot 10^{-2}$\\
\end{tabular}
\vspace{0.4cm}
\caption{\em Run~{\rm\ref{lhe-n-II}}. Relative $L^2(0,T;H^1(\Omega))$-error between the POD solution and the finite element solution (columns 2-3) and the true solution (columns 4-5), respectively, using adaptive finite element snapshots which are interpolated onto the finest mesh and using a uniform mesh.}
\label{tab:errors_heat_l2h1}
\end{table}

We note that $\varepsilon_{\text{FE}}^\text{uni}$ decays down to $10^{-8}$ ($\ell = 100$) and then stagnates if using a uniform mesh. This behavior is clear, since the more POD basis elements we include (up to stagnation of the corresponding eigenvalues), the better is the POD solution an approximation for the finite element solution. On the other hand, both $\varepsilon_{\text{true}}^\text{uni}$ and $\varepsilon_{\text{true}}^\text{ad}$ start to stagnate after $\ell = 30$ in Table \ref{tab:errors_heat}, columns 4 and 5. This is due to the fact that at this point the spatial (and temporal) discretization error dominates the modal error. This is in accordance with the decay of the eigenvalues shown in Figure \ref{fig:Heat_ev} and is accounted for e.g. in the error estimation presented in \cite[Theorem 5.1]{GH18}. Similar observations hold true for the relative $L^2(0,T;H^1(\Omega))$-error listed in Table~\ref{tab:errors_heat_l2h1} with the difference that the $L^2(0,T;H^1(\Omega))$-error is larger than the respective $L^2(0,T;L^2(\Omega))$-error.

The computational times for the full and the low order simulation using uniform finite element discretizations and adaptive finite element snapshots, which are interpolated onto the finest mesh, respectively, are listed in Table \ref{tab:CPU_times_heat}.\\
\begin{table}[H]
\centering
\begin{tabular}{ l | c | c | c}
  & adaptive FE mesh  & uniform FE mesh & speedup factor \\
 \hline
 FE simulation & 944 sec & 8808 sec & 9.3 \\
 POD offline computations & 264 sec & 1300 sec & 4.9 \\
 POD simulation & \multicolumn{2}{c|}{ 0.07 sec}  & --\\ 
 \hline
 speedup factor & 13485 &  125828 & --\\
\end{tabular}
\vspace{0.4cm}
\caption{\em Run~{\rm\ref{lhe-n-II}}. CPU times for FE and POD simulation using uniform finite element meshes and adaptive finite element snapshots which are interpolated onto the finest mesh, respectively, and using $\ell = 50$ POD modes.}
\label{tab:CPU_times_heat}
\end{table}
\noindent Once the POD basis is computed in the offline phase, the POD simulation corresponding to adaptive snapshots is 13485 times faster than the FE simulation using adaptive finite element meshes. This speedup factor is important when one  considers e.g. optimal control problems with time-dependent PDEs, where the POD-ROM can be used as surrogate model in repeated solution of the underlying PDE model. In the POD offline phase, the most expensive task is to express the snapshots with respect to the common finite element space, which takes 226 seconds. Since $\mathrm K$ \eqref{matrix-k} is symmetric, it suffices to calculate the entries on and above the diagonal, which are $\sum_{k=1}^nk=(n^2+n)/2$ entries. Thus, the computation of each entry in the correlation matrix $\mathrm K$ using a common finite element space takes around 0.00018 seconds. We note that in the approach explained in Sections \ref{POD-adapt} and \ref{POD-Galerkin}, the computation of the matrix $\mathrm K$ is expensive. For each entry the calculation time is around 0.03 seconds, which leads to a computation time of around 36997 seconds for the matrix $\mathrm K$. The same effort is needed to build $\mathrm A^\ell=a(\mathcal{Y}\Phi_j,\mathcal{Y}\Phi_i)$. In this case, the offline phase takes therefore around 88271 seconds. For this reason, the approach to interpolate the adaptively generated snapshots onto the finest mesh is computationally more favorable. But since the computation of $\mathrm K$ can be parallelized, the offline computation time can be reduced provided that the appropriate hardware is available.\hfill$\Diamond$
\end{run}

\begin{run}[Cahn-Hilliard equations]
\label{che-II}
\em
Now let us revisit Run~\ref{che}, where we in the following run the numerical simulations for different combinations of numbers for $\ell_c$ and $\ell_w$ of Table \ref{tab:loss_info}. The approximation quality of the POD solution using adaptive meshes is compared to the use of a uniform mesh in Table \ref{tab:CH_adaptiveFE}. As expected, Table \ref{tab:CH_adaptiveFE} shows that the error between the POD surrogate solution and the high-fidelity solution gets smaller for an increasing number of utilized POD basis functions. Moreover, a larger number of POD modes is needed for the chemical potential $w$ than for the phase field $c$ in order to get an error in the same order which is in accordance to the fact that the decay of the eigenvalues for $w$ is slower than for $c$ as seen in Figure \ref{fig:CH_ev}.\\
\begin{table}[H]
\centering
\begin{tabular}{  c | c || c |  c  | c | c }
  $\ell^c$ & $\ell^w$ & $c: \varepsilon_{\text{FE}}^{\text{ad}}$ &    $w: \varepsilon_{\text{FE}}^{\text{ad}}$ & $c: \varepsilon_{\text{FE}}^{\text{uni}}$ &   $w: \varepsilon_{\text{FE}}^{\text{uni}}$  \\
 \hline
 3 & 4 & $8.44 \cdot 10^{-3}$ & \hspace{-0.3cm} $3.00 \cdot 10^{0}$ &  $8.44 \cdot 10^{-3}$  & \hspace{-0.3cm} $3.75 \cdot 10^{0}$\\
 10 & 13 & $3.30 \cdot 10^{-3}$  & $3.77 \cdot 10^{-1}$ & $ 3.30 \cdot 10^{-3} $ & $4.32 \cdot 10^{-1}$ \\
 19 & 26 & $1.57 \cdot 10^{-3}$  & $2.12 \cdot 10^{-1}$ & $1.57 \cdot 10^{-3}$ &  $2.39 \cdot 10^{-1}$ \\
 29 & 26 & $7.34 \cdot 10^{-4}$  & $1.09 \cdot 10^{-1}$ & $7.32 \cdot 10^{-4}$ &  $1.16 \cdot 10^{-1}$\\
 37 & 26 & $3.57 \cdot 10^{-4}$  & $4.82 \cdot 10^{-2}$ & $ 3.55 \cdot 10^{-4}$ &  $5.04 \cdot 10^{-2}$  \\

 50 & 50 & $1.88 \cdot 10^{-4}$  & $2.17 \cdot 10^{-2}$ & $1.86 \cdot 10^{-4}$ &  $2.33 \cdot 10^{-2} $ \\
 65 & 26 & $9.74 \cdot 10^{-5}$  & $1.11 \cdot 10^{-2}$ & $9.56 \cdot 10^{-5}$ &  $1.15 \cdot 10^{-2}$\\
 100 & 100 & $3.37 \cdot 10^{-5}$ & $3.56 \cdot 10^{-3}$ & $3.22 \cdot 10^{-5}$  & $3.42 \cdot 10^{-3}$ \\
\end{tabular}
\vspace{0.4cm}
\caption{\em Run~{\rm\ref{che-II}}. Relative $L^2(0,T;L^2(\Omega))$-error between the POD solution and the finite element solution using adaptive meshes (columns 3-4) and using a uniform mesh (columns 5-6), respectively.}
\label{tab:CH_adaptiveFE}
\end{table}
\noindent
We now discuss the treatment of the nonlinearity and also investigate the influence of non-smoothness of the model equations to the POD procedure. Using the convex-concave splitting for $W$, we obtain for the Moreau-Yosida relaxed double obstacle free energy the concave part $W_-^{\text{rel}} (c) = \frac{1}{2}(1-c^2)$ and the convex part $W_+^{\text{rel}}(c) = \frac{s}{2}(\max (c-1,0)^2 + \min (c+1,0)^2)$. This means that the first derivative of the concave part is linear with respect to the phase field variable $c$. The challenging part is the convex term with non-smooth first derivative. For a comparison, we consider the smooth polynomial free energy with concave part $W_-^p(c) = \frac{1}{4}(1-2c^2)$ and convex part $W_+^p(c) = \frac{1}{4}c^4$.

\noindent
Figure~\ref{fig:ev_psiprime} shows the decay of the normalized eigenspectrum for the phase field $c$ (left) and the first derivative of the convex part $W'_+(c)$ (right) for the polynomial and the relaxed double obstacle free energy. Obviously, in the non-smooth case more POD modes are needed for a good approximation than in the smooth case. This behavior is similar to the decay of the Fourier coefficients in the context of trigonometric approximation, where the decay of the Fourier coefficients depends on the smoothness of the approximated object. 
\begin{figure}[H]
\centering
\includegraphics[width=.45\textwidth]{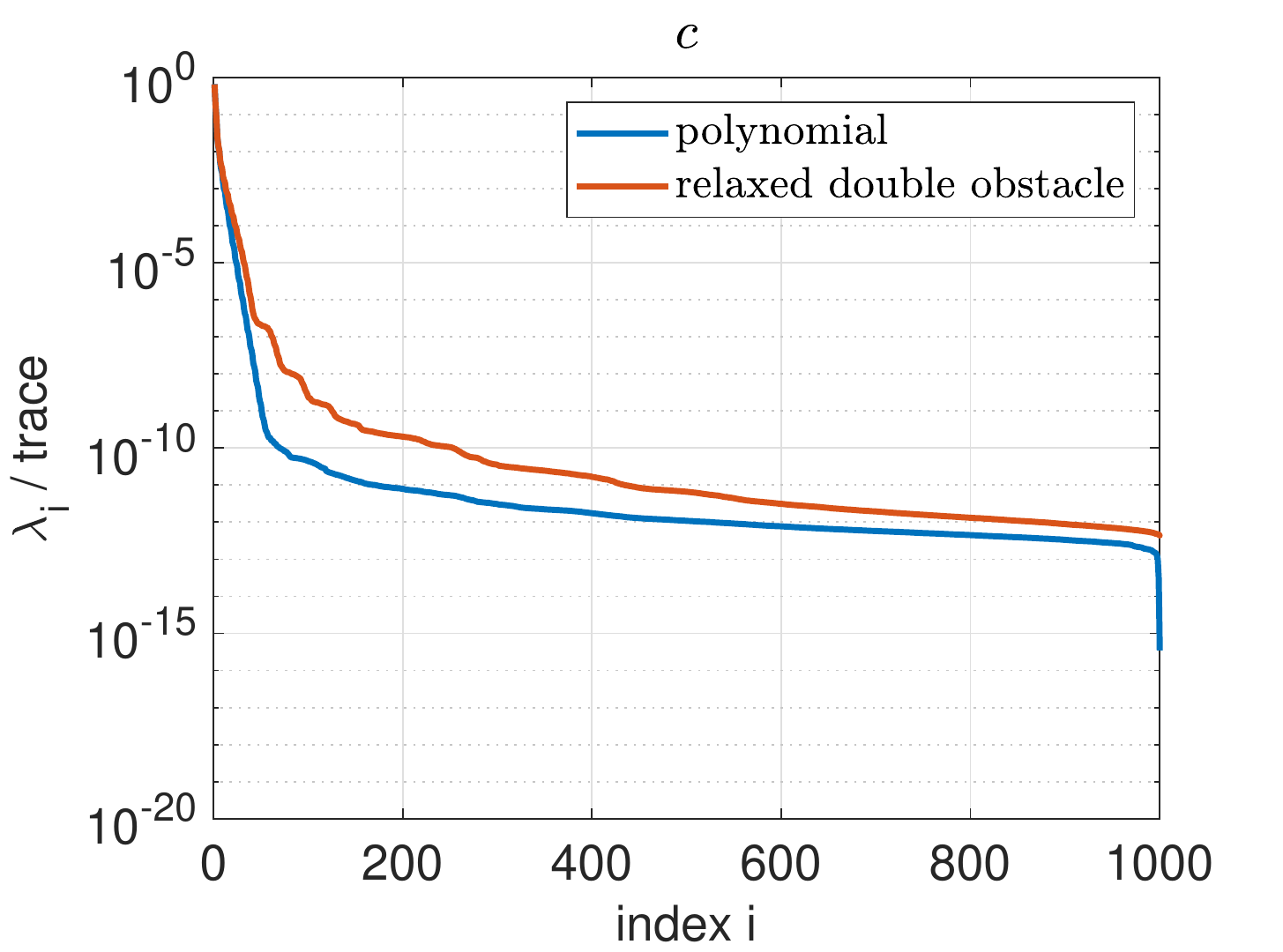} \hfill \includegraphics[width=.45\textwidth]{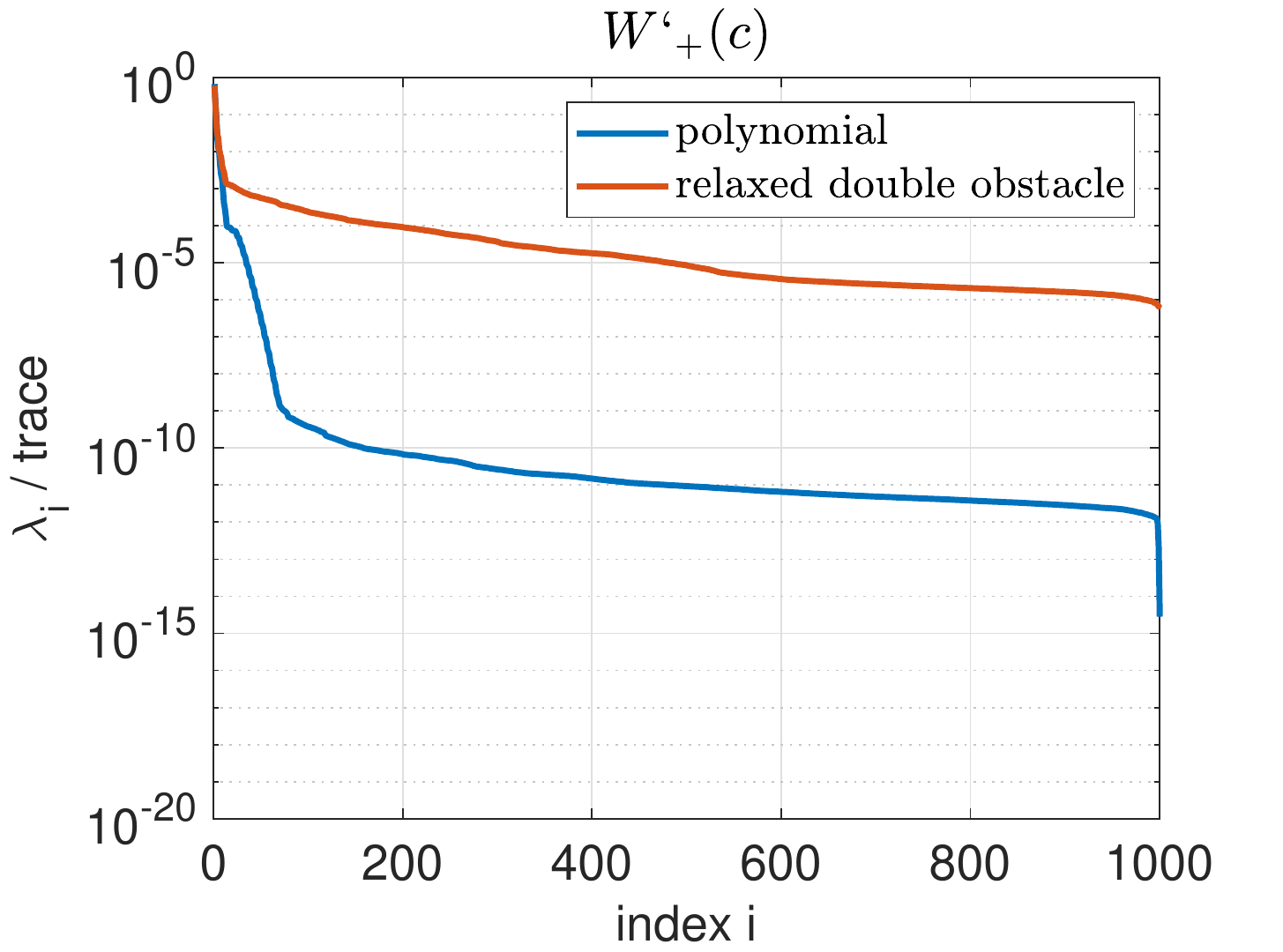} 
\caption{\em Run~{\rm\ref{che-II}}. Comparison of the normalized eigenvalues for $c$ (left) and the first derivative of the convex part $W'_+$ of the free energy (right) using polynomial and relaxed double obstacle energy, respectively.}
\label{fig:ev_psiprime}
\end{figure}

\noindent
Table \ref{tab:CH_speedups} summarizes computational times for different finite element runs as well as reduced-order simulations using the polynomial and the relaxed double obstacle free energy, respectively. In addition, the approximation quality is compared. The computational times are rounded averages from various test runs. It turns out that the finite element simulation (row 1) using the smooth potential is around two times faster than using the non-smooth potential. This is due to the fact that in the smooth case, two to three Newton steps are needed for convergence in each time step, whereas in the non-smooth case six to eight iterations are needed in the semismooth Newton method.\\
\begin{table}[H]
\centering
\begin{tabular}{  l | r | r ||  r  | r  }
  & \multicolumn{2}{c||}{$W^p$} & \multicolumn{2}{c}{$W_s^{\text{rel}}$} \\
 \hline
  FE & \multicolumn{2}{c||}{ 1644 s}  & \multicolumn{2}{c}{ 3129 s} \\[0.2cm]
 \hline
  & $\ell_c=3$ & $\ell_c=19$ &  $\ell_c=3$ & $\ell_c=19$\\
  & $\ell_w=4$ & $\ell_w=26$ &  $\ell_w=4$ & $\ell_w=26$\\
\hline
 POD offline &  355 s & 355 s & 350 s & 349 s\\
 DEIM offline & 8 s & 8 s & 9 s & 10 s \\
 ROM & 183 s & 191 sec & 2616 s & 3388 s\\
 ROM-DEIM & 0.05 s & 0.1 s & 0.04 s & no conv. \\
 ROM-proj & 0.008 s & 0.03 s & 0.01 s & 0.03 s\\
 \hline
 speedup FE-ROM & 8.9 & 8.6 & 1.1 &  none \\
  speedup FE-ROM-DEIM &  32880 &  16440  &  78225 & --\\
 speedup FE-ROM-proj &  205500 &  54800 &  312900 & 104300 \\
 \hline
  rel $L^2(Q)$ error ROM & $5.46 \cdot 10^{-3}$ & $3.23 \cdot 10^{-4}$ &  $8.44 \cdot 10^{-3}$ &  $1.57 \cdot 10^{-3}$  \\
    rel $L^2(Q)$ error ROM-DEIM & $1.46 \cdot 10^{-2}$ & $ 3.83 \cdot 10^{-4}$  & $8.84 \cdot 10^{-3}$  & --\\
   rel $L^2(Q)$ error ROM-proj & $4.70 \cdot 10^{-2}$ & $ 4.18 \cdot 10^{-2}$  &  $8.72 \cdot 10^{-3}$  &  $ 9.80 \cdot 10^{-3}$ \\
\end{tabular}
\vspace{0.4cm}
\caption{\em Run~{\rm\ref{che-II}}. Computational times, speedup factors and approximation quality for different POD basis lengths and using different free energy potentials.}
\label{tab:CH_speedups}
\end{table}

Using the smooth polynomial free energy, the reduced-order simulation is 8-9 times faster than the finite element simulation, whereas using the relaxed double obstacle free energy only delivers a very small speedup. The inclusion of DEIM (we use $\ell_{\text{deim}}=\ell_c$) in the reduced-order model leads to immense speedup factors for both free energy functions (row 8). This is due to the fact that the evaluation of the nonlinearity in the reduced-order model is still dependent on the full spatial dimension and hyper reduction methods are necessary for useful speedup factors. Note that the speedup factors are of particular interest in the context of optimal control problems. At the same time, the relative $L^2(0,T;L^2(\Omega))$-error between the finite element solution and the ROM-DEIM solution is close to the quality of the reduced-order model solution (row 10-11).

\noindent
However, in the case of the non-smooth free energy function using $\ell_c=19$ POD modes for the phase field and $\ell_w=26$ POD modes for the chemical potential, the inclusion of DEIM has the effect that the semismooth Newton method does not converge. For this reason, we treat the nonlinearity by applying the technique explained in Section~\ref{P2.1}, i.e. we project the finite element snapshots for $W'_+(c)$ (which are interpolated onto the finest mesh) onto the POD space. Since this leads to linear systems, the computational times are very small (row 6). The error between the finite element solution and the reduced-order solution using projection of the nonlinearity is of the magnitude $10^{-02}/10^{-03}$. Depending on the motivation, this approximation quality might be sufficient. Nevertheless, we note that for large numbers of POD modes, using the projection of the nonlinearity onto the POD space leads to a large increase of the error.\hfill$\Diamond$
\end{run}

To summarize, a POD reduced-order model construction approach is proposed which can be set up and solved for snapshots originating from arbitrary FE (and also other) spaces. The method is applicable for $h$-, $p$- and $r$-adaptive finite elements. It is motivated from an infinite-dimensional perspective. Using the method of snapshots we are able to set up the correlation matrix $K$ from \eqref{matrix-k} by evaluating the inner products of snapshots which live in different FE spaces. For non-nested meshes, this requires the detection of cell collision and integration over cut finite elements. A numerical strategy how to implement this practically is elaborated and numerically tested. Using the eigenvalues and eigenvectors of this correlation matrix, we are able to set up and solve a POD surrogate model that does not need the expression of the snapshots with respect to the basis of a common FE space or the interpolation onto a common reference mesh. Moreover, an error bound for the error between the true solution and the solution to the POD-ROM using spatially adapted snapshots is available in \cite[Theorem 5.1]{GH18}. The numerical tests show that the POD projection error decreases if the number of utilized POD basis functions is increased. However, the error between the POD solution and the true solution stagnates when the spatial discretization error dominates. Moreover, the numerics show that using the correlation matrix calculated explicitly without interpolation in order to build a POD-ROM gives the same results as the approach where the snapshots are interpolated onto the finest mesh. From a computational point of view, sufficient hardware should be available in order to compute the correlation matrix in parallel and make the offline computational time competitive. For semilinear evolution problems, the nonlinearity is treated by linearization. This is of interest in view of optimal control problems, in which a linearized state equation has to be solved in each SQP iteration level. An appropriate treatment of the nonlinearity in our applications gains significant speedup of the ROM in computational times when compared to the full simulations. This makes POD-MOR with adaptive finite elements an ideal approach for the construction of surrogate models in e.g. optimal control with nonlinear PDE systems as they arise e.g. in the context of multi-phase flow control problems.

%%%%%%%%%%%%%%%%%%%%%%%%%%%%%%%%%%%%%%%%%%%%%%%%%%%
\section{Certification with a priori and a posteriori error estimates}
\label{POD-cert}
%%%%%%%%%%%%%%%%%%%%%%%%%%%%%%%%%%%%%%%%%%%%%%%%%%%

As we have seen in Section~\ref{POD-Galerkin} POD provides a method for deriving low order models of dynamical systems. It can be thought of as a Galerkin approximation in the spatial variable, built from functions corresponding to the solution of the physical system at prespecified time instances. After carrying out a singular value decomposition the leading $\ell$ generalized eigenfunctions are chosen as the POD basis $\{\Psi_j\}_{j=1}^\ell$ of rank $\ell$. As soon as one uses POD, questions concerning the quality of the approximation properties, convergence, and rate of convergence become relevant. Let us refer, e.g., to the literature \cite{CGS12,GV17,KV01,KV02,KV02b,SS13,Sch12,Sin14,Rav11} for a priori error analysis for POD Galerkin approximations. It turns out that the error depends on the decay of the sum $\sum_{i>\ell}\lambda_i$, the error $\Delta t^\beta$ (with an appropriate $\beta\ge1$) due to the used time integration method, the used Galerkin spaces $\{V^{h_j}\}_{j=1}^n$  and the choice $X=H$ or $X=V$. In particular, best approximation properties hold provided the time differences $\dot y^h(t_j)$ (or the finite difference discretizations) are included in the snapshot ensembles; cf. \cite{KV01,KV02,Sin14}.

Let us recall numerical test examples from \cite[Section~1.5]{GV17}. The programs are written in {\sc Matlab} using the {\sc Partial Differential Equation Toolbox} for the computation of the piecewise linear FE discretization. For the temporal integration the implicit Euler method is applied based on the equidistant time grid $t_j=(j-1)\Delta t$, $j=1,\ldots,n$ and $\Delta t=T/(n-1)$.

\begin{run}[POD for the heat equation; cf. {\cite[Run~1]{GV17}}]
\label{Run-GV1}
\em
We choose the final time $T=3$, the spatial domain $\Omega=(0,2)\subset\mathbb R$, the Hilbert spaces $H=L^2(\Omega)$, $V=H^1_0(\Omega)$, the source term $f(t,\bx)=t^3-\bx^2$ for $(t,\bx)\in Q=(0,T)\times\Omega$ and the discontinuous initial value $y_\circ(\bx)=\chi_{(0.5,1.0)}-\chi_{(1,1.5)}$ for $\bx\in\Omega$, where, e.g., $\chi_{(0.5,1)}$ denotes the characteristic function on the subdomain $(0.5,1)\subset\Omega$, $\chi_{(0.5,1)}(\bx)=1$ for $\bx\in (0.5,1)$ and $\chi_{(0.5,1)}(\bx)=0$ otherwise. We consider a discretization of the linear heat equation (compare \eqref{heat} with $c \equiv 0$)
\begin{equation}
\label{HeatEqGV}
\begin{aligned}
y_t(t,\bx)-\Delta y(t,\bx) &=f(t,\bx)&&\text{for  }(t,\bx)\in Q,\\
y(t,\bx)&= 0&& \text{for }(t,\bx)\in\Sigma=(0,T)\times\partial\Omega,\\
y(0,\bx)&=y_\circ(\bx)&& \text{for }\bx\in\Omega.
\end{aligned}
\end{equation}
To obtain an accurate approximation of the exact solution we choose $n=4000$ so that $\Delta t\approx 7.5\cdot 10^{-4}$ holds. For the FE discretization we choose $m=500$ spatial grid points and  the equidistant mesh size $h=2/(m+1)\approx 4\cdot 10^{-3}$. Thus, the FE error -- measured in the $H$-norm -- is of the order $10^{-4}$. In the left graphic of Figure~\ref{Figure5.1}, the FE solution $y^h$ to the state equation \eqref{HeatEqGV} is visualized.
\begin{figure}
\includegraphics[height=32mm,width=52mm]{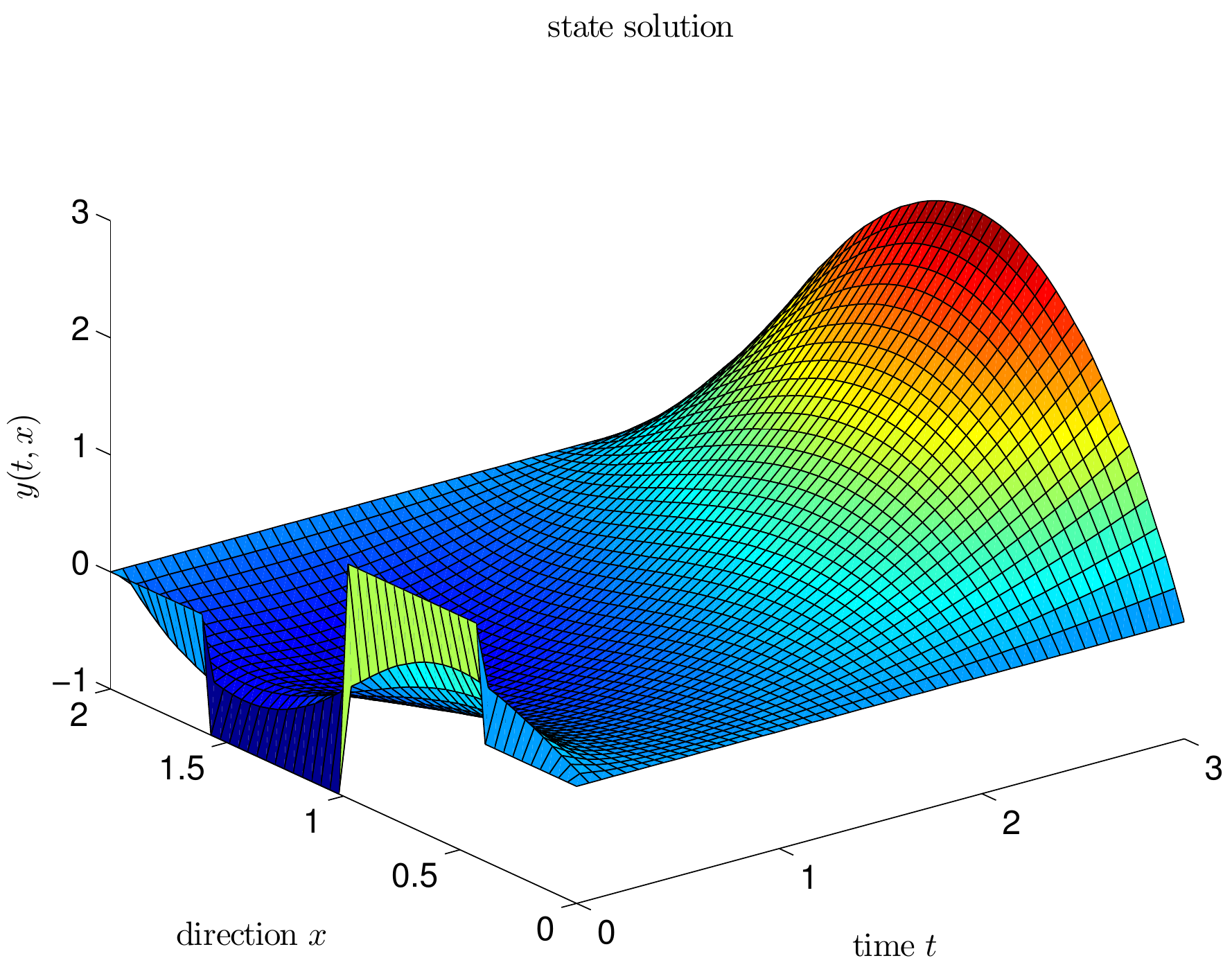}\hspace{5mm}
\includegraphics[height=32mm,width=52mm]{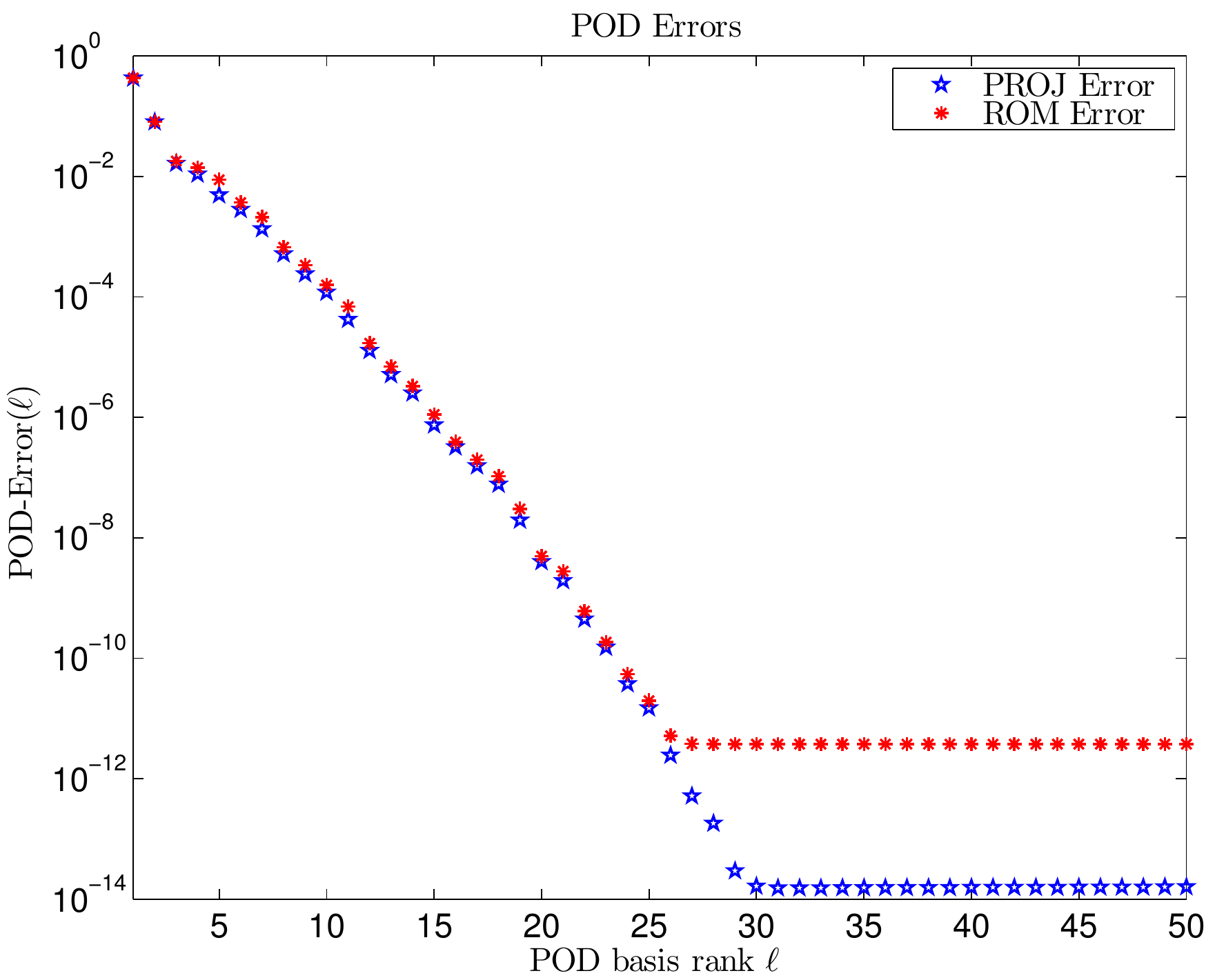}
\caption{Run~\ref{Run-GV1} (cf. \cite[Figure~1.1]{GV17}). The FE solution $y^h$ (left) and the residuals corresponding to the POD basis rank $\ell$ (right).}
\label{Figure5.1}
\end{figure}
To compute a POD basis $\{\Psi_i\}_{i=1}^\ell$ of rank $\ell$ we utilize the multiple discrete snapshots $y^1_j=y^h(t_j)$ for $1\le j\le n_t$ as well $y^2_1=0$ and $y^2_j=(y^h(t_j)-y^h(t_{j-1})/\Delta t$, $j=2,\ldots,n_t$, i.e., we include the temporal difference quotients in the snapshot ensemble and $K=2$, $n_1=n_2=n_t$. We choose $X=H$ and utilize the (stable) SVD to determine the POD basis of rank $\ell$; compare Section~\ref{P1.2}. We address this issue in a more detail in Run~\ref{Run-GV4}. Since the snapshots are FE functions, the POD basis elements are also FE functions. In the right plot of Figure~\ref{Figure5.1}, the projection and reduced-order error given by
\begin{align*}
\mathrm{PROJ~Error}(\ell)&=\bigg(\sum_{j=1}^{n_t}\alpha_j\,\Big\|y^h(t_j)-\sum_{i=1}^\ell {\langle y^h(t_j),\psi_i\rangle}_H\,\psi_i\Big\|_H^2\bigg)^{1/2},\\
\mathrm{ROM~Error}(\ell)&=\bigg(\sum_{j=1}^{n_t}\alpha_j\,\big\|y^h(t_j)-y^\ell(t_j)\big\|_H^2\bigg)^{1/2}
\end{align*}
are plotted for different POD basis ranks $\ell$. The chosen trapezoidal weights $\alpha_j$ have been introduced in \eqref{alpha}. We observe that both errors decay rapidly and coincide until the accuracy $10^{-12}$, which is already significant smaller than the FE discretization error. These numerical results reflect the a priori error estimates presented in \cite[Theorem~1.29]{GV17}.\hfill$\Diamond$
\end{run}

\begin{run}[POD for a convection dominated heat equation; cf. {\cite[Run~2]{GV17}}]
\label{Run-GV2}
\em
Now we consider a more challenging example. We study a convection-reaction-diffusion equation with a source term which is close to being singular: Let $T$, $\Omega$, $y_\circ$, $H$ and $V$ be given as in Run~\ref{Run-GV1}. The parabolic problem reads as follows
\begin{align*}
y_t(t,\bx)-cy_{\bx\bx}(t,\bx)+\beta y_x(t,\bx)+ay(t,\bx) &=f(t,\bx)&&\text{for  }(t,\bx)\in Q,\\
y(t,\bx) &= 0&&\text{for  }(t,\bx)\in \Sigma,\\
y(0,\bx)&=y_\circ(\bx)&&\text{for }\bx\in\Omega.
\end{align*}
We choose the diffusivity $c=0.025$, the velocity $\beta=1.0$ that determines the speed in which the initial profile $y_\circ$ is shifted to the boundary and the reaction rate $a=-0.001$. Finally, $f(t,\bx)=\mathbb P(\frac1{1-t})\cos(\pi \bx)$ for $(t,\bx)\in Q$, where $(\mathbb Pz)(t)=\min(+l,\max(-l,z(t)))$ restricts the image of $z$ on a bounded interval. In this situation, the state solution $y$ develops a jump at $t=1$ for $l\to\infty$; see the left plot of Figure~\ref{Figure5.2}.
\begin{figure}
\includegraphics[height=32mm,width=52mm]{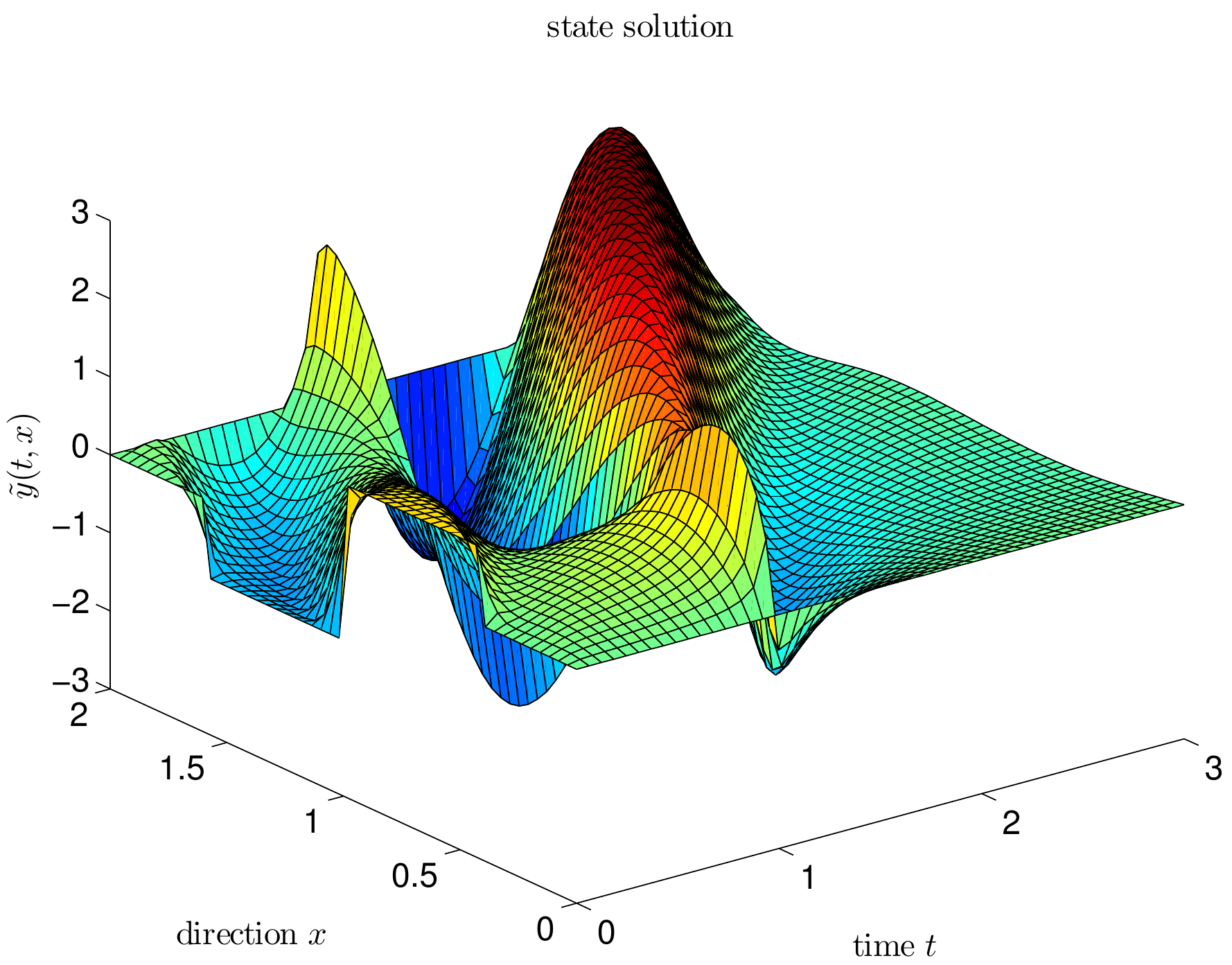}\hspace{5mm}
\includegraphics[height=32mm,width=52mm]{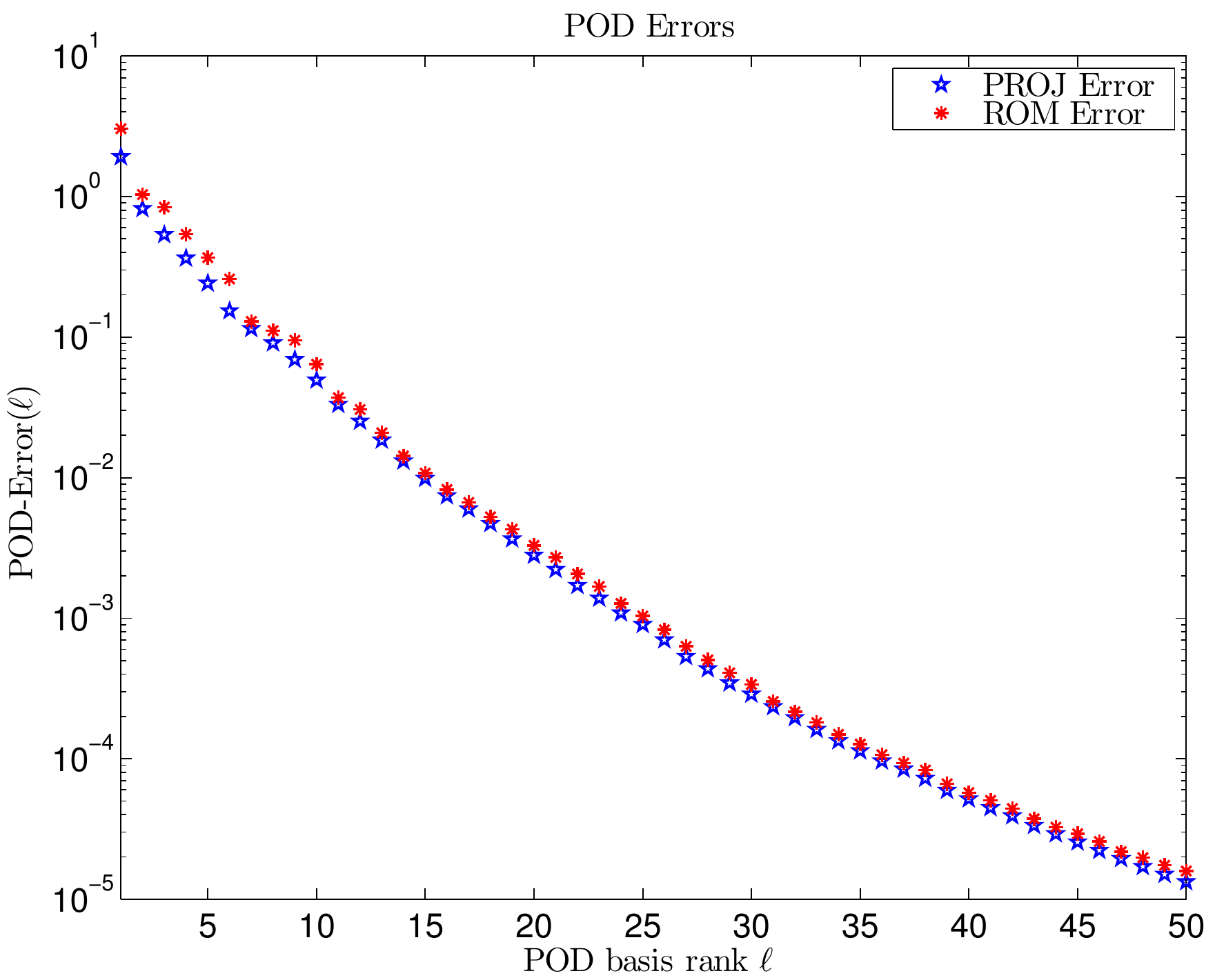}
\caption{Run~\ref{Run-GV2} (cf. \cite[Figure~1.2]{GV17}. The FE solution $y^h$ (left) and the residuals corresponding to the POD basis rank $\ell$ (right).}
\label{Figure5.2}
\end{figure}
The right plot of Figure~\ref{Figure5.2} demonstrates that in this case, the decay of the reconstruction residuals and the decay of the errors are much slower as in the right plot of Figure~\ref{Figure5.1}. The manifold dynamics of the state solution require an inconvenient large number of POD basis elements. Since the supports of these ansatz functions in general cover the whole domain $\Omega$, the corresponding system matrices of the reduced model are not sparse. This is different for the matrices arising in the FE Galerkin framework. Model order reduction is not effective for this example if a good accuracy of the solution function $y^\ell$ is required. Strategies to improve the accuracy and robustness of the POD-ROM in those situations are discussed in e.g. \cite{BBI09,WWXI18}\hfill$\Diamond$
\end{run}

\begin{run}[True and exact approximation error; cf. {\cite[Run~3]{GV17}}]\label{Run-GV3}
\em
We consider the setting introduced in Run~\ref{Run-GV1} again. The exact solution to  \eqref{HeatEqGV} does not possess a representation by elementary functions. Hence, the presented reconstruction and reduction errors actually are the residuals with respect to a high-order FE solution $y^h$. To compute an approximation $y$ of the exact solution $y_\mathsf{ex}$ we apply a Crank-Nicolson method (with Rannacher smoothing \cite{Ran84}) ensuring $\|y-y_\mathsf{ex}\|_{L^2(0,T;H)}=\mathcal O(\Delta t^2+h^2)\approx10^{-5}$. In the context of model reduction, such a state is sometimes called the ``true'' solution. To compute the FE state $y^h$ we apply the Euler method. In the left plot of Figure~\ref{Figure5.3} we compare the true solution $y_\mathsf{ex}$ with the associated POD approximation for different values $n_t\in\{64,128,256,...,8192\}$ of the time integration and for the spatial mesh size $h=4\cdot10^{-3}$.
\begin{figure}
\includegraphics[height=32mm,width=52mm]{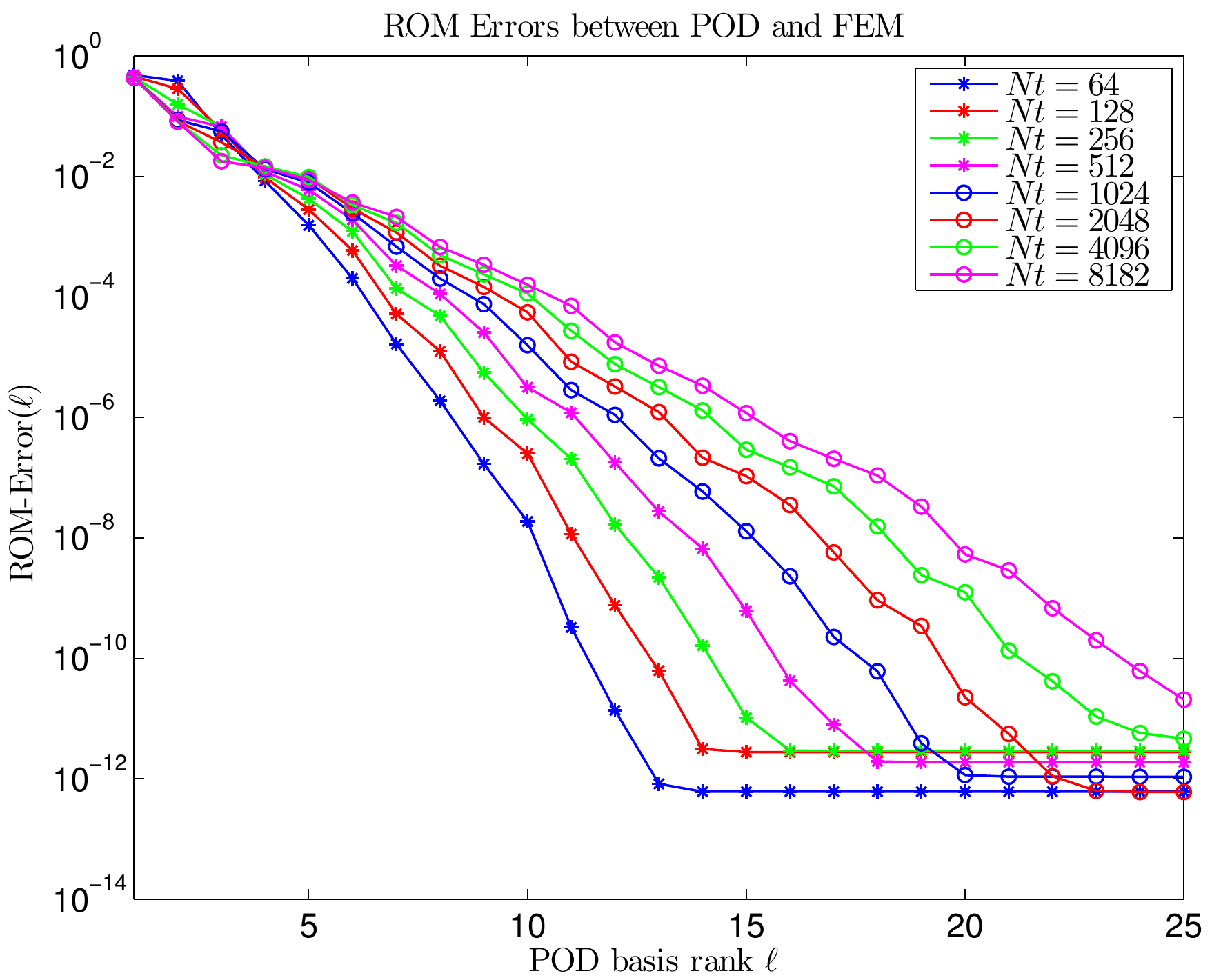}\hspace{5mm}
\includegraphics[height=32mm,width=52mm]{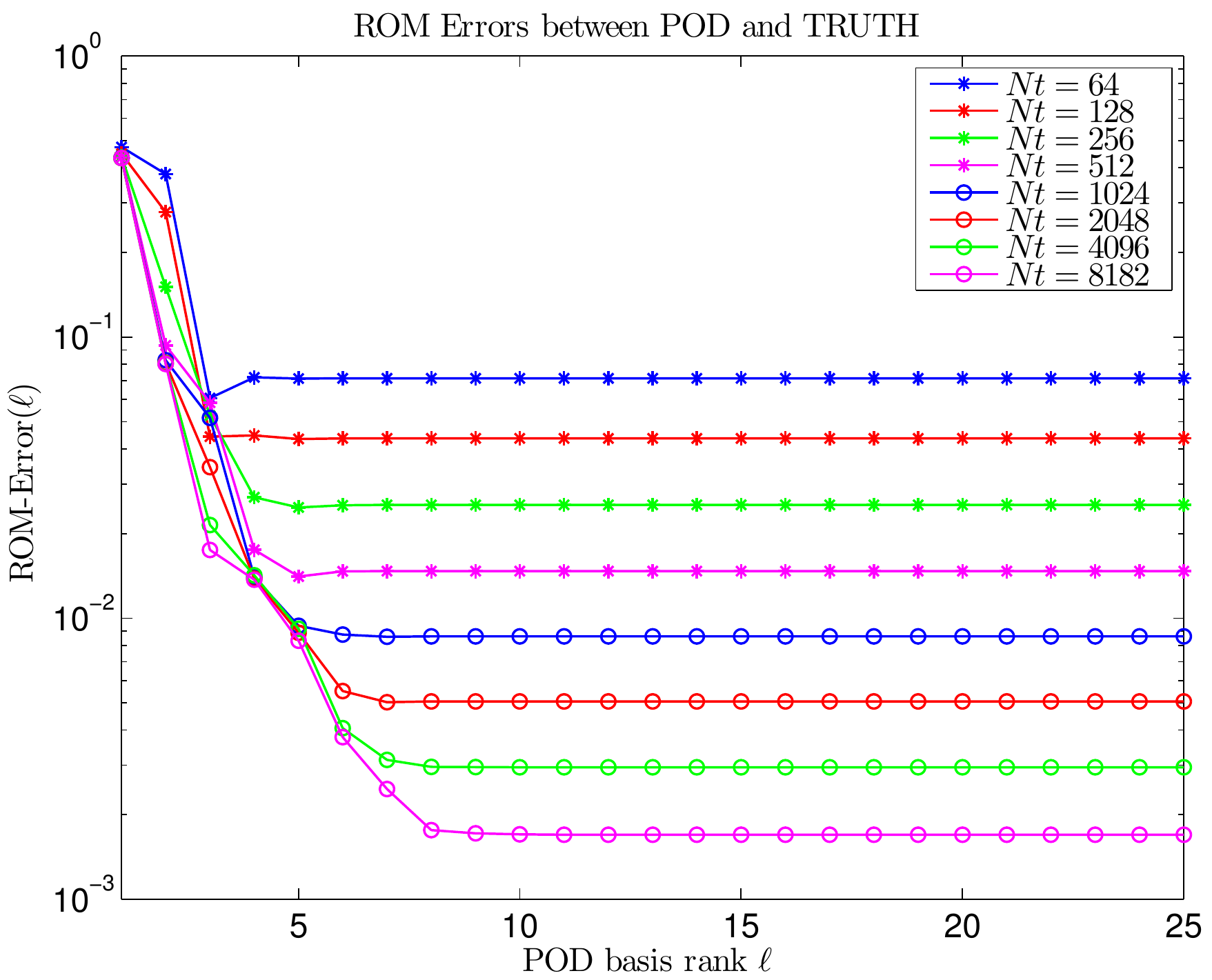}
\caption{Run~\ref{Run-GV3} (cf. \cite[Figure~1.3]{GV17}. The ROM errors with respect to the true solution (left) and the exact one (right).}
\label{Figure5.3}
\end{figure}
For the norm we apply a discrete $L^2(0,T;H)$-norm as in Run~\ref{Run-GV1}. Let us mention that we compute for every $n_t$ a corresponding FE solution $y^h$. We observe that the residuals ignore the errors arising by the application of time and space discretization schemes for the full-order model. The errors decay below the discretization error $10^{-5}$. If these discretization errors are taken into account, the residuals stagnate at the level of the full-order model accuracy instead of decaying to zero; cf. right plot of Figure~\ref{Figure5.3}. Due to the implicit Euler method we have $\|y^h-y_\mathsf{ex}\|_{L^2(0,T;H)}=\mathcal O(\Delta t+h^2)$ with the mesh-size $h=4\cdot10^{-3}$. In particular, from $n_t\in\{64,128,256,...,8192\}$ it follows that $\Delta t>3\cdot10^{-4}>h^2=1.6\cdot10^{-5}$. Therefore, the spatial error is dominated by the time error for all values of $n_t$. We can observe that the exact residuals do not decay below a limit of the order $\Delta t$. One can observe that for fixed POD basis rank $\ell$, the residuals with respect to the true solution increase if the high-order accuracy is improved by enlarging $n_t$, since the reduced-order model has to approximate a more complex system in this case, where the residuals with respect to the exact solution decrease due to the lower limit of stagnation $\Delta t=3/(n_t-1)$.\hfill$\Diamond$
\end{run}

\begin{run}[Different strategies for a POD basis computation; cf. {\cite[Run~4]{GV17}}]
\label{Run-GV4}
\em
As we have explained in Section~\ref{P1.2}, let $Y\in\mathbb R^{m\times n}$ denote the matrix of snapshots with rank $\mathsf r$, $W\in\mathbb R^{m\times m}$ be the (sparse) spatial weighting matrix consisting of the elements $\langle\varphi_j,\varphi_i\rangle_X$ (introduced Section~\ref{P1.3.3}) and $D\in\mathbb R^{n\times n}$ be the diagonal matrix containing the nonnegative weighting parameters $\alpha_j^k$. As we have explained in Section~\ref{P1.2}, the POD basis $\{\Psi_i\}_{i=1}^\ell$ of rank $\ell\le\mathsf r$ can be determined by providing an eigenvalue decomposition of the matrix $\bar Y\bar Y^\top=W^{1/2}YDY^\top W^{1/2}\in\mathbb R^{m\times m}$,  one of $\bar Y^\top\bar Y=D^{1/2}Y^\top WYD^{1/2}\in\mathbb R^{n\times n}$, or a singular value decomposition of $\bar Y=W^{1/2}YD^{1/2}\in\mathbb R^{m\times n}$. Since $n\gg m$ in Runs~\ref{Run-GV1}-\ref{Run-GV3}, the first variant is the cheapest one from a computational point of view. In case of multiple space dimensions or if a second-order time integration scheme such as some Crank-Nicolson technique is applied, the situation is converse. On the other hand, a singular value decomposition is more accurate and stable than an eigenvalue decomposition if the POD elements corresponding to eigenvalues/singular values which are close to zero are taken into account: Since $\lambda_i=\sigma_i^2$ holds for all eigenvalues $\lambda_i$ and singular values $\sigma_i$, the singular values are able to decay to machine precision, where the eigenvalues stagnate significantly above. This is illustrated in the left graphic of Figure~\ref{Figure5.4}.
\begin{figure}
\includegraphics[height=32mm,width=52mm]{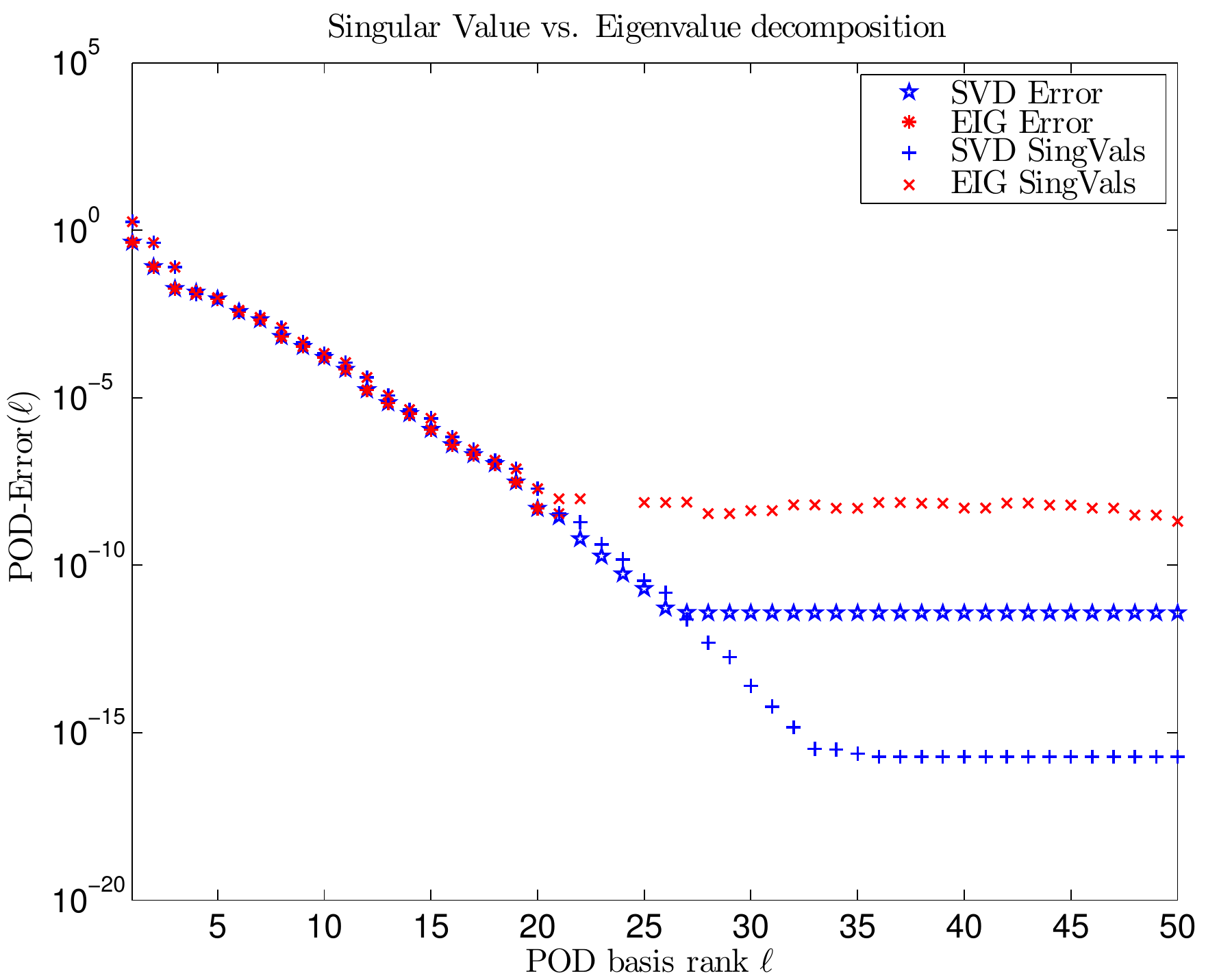}\hspace{5mm}
\includegraphics[height=32mm,width=52mm]{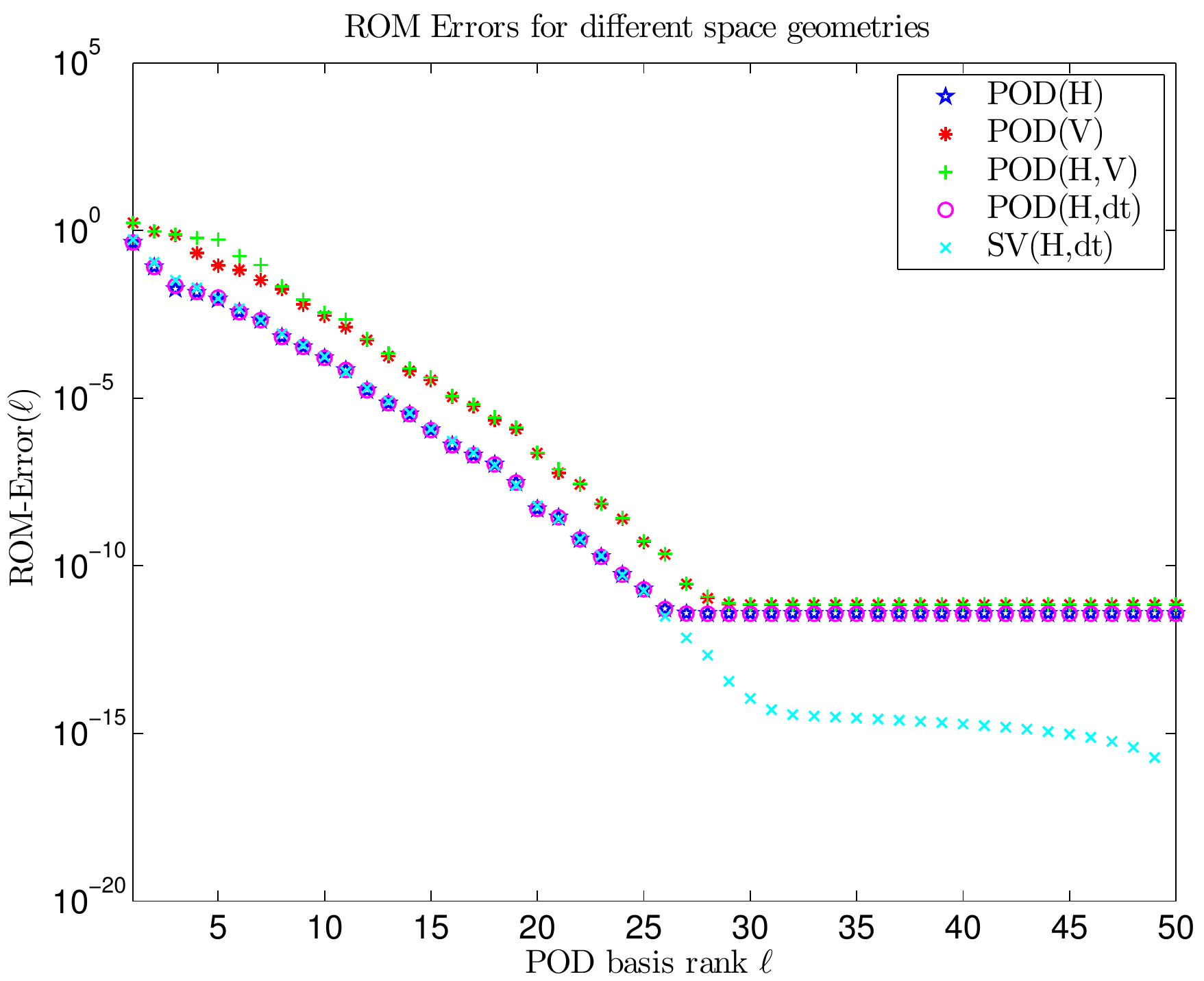}
\caption{Run~\ref{Run-GV4} (cf. \cite[Figure~1.4]{GV17}. Singular values $\sigma_i$ using the SVD (SVD Vals) or the eigenvalue decomposition (EIG Vals) and the associated ROM errors (SVD error and EIG Error, respectively) (left); ROM errors for different the choices for $X$, the error norm and the snapshot ensembles (right).}
\label{Figure5.4}
\end{figure}
Indeed, for $\ell>20$ the EIG-ROM system matrices become singular due to the numerical errors in the eigenfunctions and the reduced-order system is ill-posed in this case while the SVD-ROM model remains stable. In the right plot of Figure~\ref{Figure5.4} POD elements are constructed with respect to different scalar products and the resulting ROM errors are compared: $\|\cdot\|_H$-residuals for $X=H$ (denoted by POD(H)), $\|\cdot\|_V$-residuals for $X=V$ (denoted by POD(V)), $\|\cdot\|_V$-residuals for $X=H$ (denoted by POD(H,V)), which also works quite well, the consideration of time derivatives in the snapshot sample (denoted by POD(H,dt)) which allows to apply the a priori error estimate given in \cite[Theorem~1.29-2)]{GV17} and the corresponding sums of singular values (denoted by SV(H,dt)) corresponding to the unused eigenfunctions in the latter case which indeed nearly coincide with the ROM errors. \hfill$\Diamond$
\end{run}

Notice that in many applications, the quality of the reduced-order model does not vary significantly if the weights matrix $W$ refers to the space $X=H$ or $X=V$ and if time derivatives of the used snapshots are taken into account or not. Especially, the ROM residual decays with the same order as the sum over the remaining singular values, independent of the chosen geometrical framework.

%%%%%%%%%%%%%%%%%%%%%%%%%%%%%%%%%%%%%%%%%%%%%%%%%%%
\section{Optimal snapshot location for computing POD basis functions}
\label{POD-snap}
%%%%%%%%%%%%%%%%%%%%%%%%%%%%%%%%%%%%%%%%%%%%%%%%%%%

The construction of reduced-order models for nonlinear dynamical systems using proper orthogonal decomposition (POD) is based on the information carried of the so-called snapshots. These provide the spatial distribution of the nonlinear system at discrete time instances. Thus, we are interested in optimizing the choice of these time instances in such a manner that the error between the POD-solution and the trajectory of the dynamical system is minimized. This approach was suggested in \cite{KV10} and was extended in \cite{LV14} to parametrized elliptic problems. Let us briefly mention some related issues of interest. In \cite{BDW04,ES95} the situation of missing snapshot data is investigated and gappy POD is introduced for their reconstruction. An important alternative to POD model reduction is given by reduced basis approximations; we refer to \cite{PR06} and references given there. In \cite{GMNP07} a reduced model is constructed for a parameter dependent family of large scale problems by an iterative procedure that adds new basis variables on the basis of a greedy algorithm. In the Ph.D thesis \cite{Bui07} a model reduction is sought of a class for a family of models corresponding to different operating stages.

Suppose that we are given the $n_t$ snapshots $\{y(t_j)\}_{j=1}^{n_t}\subset V\subset X$. The goal is to determine additional $\mathsf k$ snapshots at time instances $\tau=(\tau_1,\ldots,\tau_\mathsf k)$ with $0\le\tau_j\le T$, $j=1,\ldots,\mathsf k$. In \cite{KV10} we propose to determine $\tau=(\tau_1,\ldots,\tau_\mathsf k)$ by solving the optimization problem
\begin{equation}
\label{Poptsnap}
\min_{0\le\tau_1,\ldots,\tau_\mathsf k\le T}\int_0^T{\|y(t)-y^\ell(t)\|}_V^2\,\mathrm dt,
\end{equation}
where $y$ and $y^\ell$ are the solutions to \eqref{State} and its POD Galerkin approximation, respectively. Clearly, the definition of the operator $\mathcal R$ given in \eqref{OperatorR} has to be modified as follows:
\[
\mathcal R^\tau\Psi=\sum_{j=1}^{n_t}\alpha_j^\tau\,{\langle y(t_j),\Psi\rangle}_X\,y(t_j)+\sum_{j=1}^{\mathsf k}\alpha^\tau_{n_t+j}\,{\langle y(\tau_j),\Psi\rangle}_X\,y(\tau_j)
\]
with appropriately modified (trapezoidal) weights $\alpha_j^\tau$, $j=1,\ldots,\mathsf k+n_t$. Consequently, \eqref{Poptsnap} becomes an optimization problem subject to the equality constraints
\[
\mathcal R^\tau\Psi_i=\lambda_i\Psi_i,\quad i=1,\ldots,\ell.
\]
Note that no precautions are made in \eqref{Poptsnap} to avoid multiple appearance of a snapshot. In fact, this would simply imply that a specific snapshot location should be given a higher weight than others. While the presented approach shows how to choose optimal snapshots in evolution equations, a similar strategy is applicable in the context of parameter dependent systems.

It turns out in our numerical tests carried out in \cite{KV10} that the proposed criterion is sensitive with respect to the choice of the time instances. Moreover, the tests demonstrate the feasibility of the method in determining optimal snapshot locations for concrete diffusion equations.

\begin{run}[cf. {\cite[Run~1]{KV10}}]
\label{Run-KV1}
\em
For $T=1$ let $Q=(0,T) \times \Omega$ and $\Omega= (0,1) \times (0,1) \subset \mathbb R^2$. For the FE triangulation we choose a uniform grid with mesh size $h=1/40$, i.e., we have 900 degrees of freedom for the spatial discretization. Then, we consider
\begin{align*}
y_t(t,\bx)-c\Delta y(t,\bx)+\beta \cdot \nabla y(t,\bx)+ y(t,\bx)&=f(\bx)&&\text{for }(t,\bx)\in Q,\\
c\,\frac{\partial y}{\partial\bn}(t,\bx)+q(\bx)y(t,\bx)&=g(\bx)&&\text{for }(t,\bx)\in\Sigma,\\
y(0,\bx)&=y_\circ(\b x)&&\text{for }\bx\in\Omega,
\end{align*}
where $c=0.1$, $\beta=(0.1,-10)^\top\in\mathbb R^2$,
\[
f(\bx)= \left\{
\begin{array}{ll}
4 & \text{for all } \bx=(x_1,x_2) \text{ with } (x_1-0.25)^2+(x_2-0.65)^2 \le 0.05,\\
0 & \text{otherwise,}
\end{array}
\right.
\]
and $y_\circ(\bx)=\sin(\pi x_1)\cos(\pi x_2)$ for $\bx=(x_1,x_2) \in \Omega$ (see Figure~\ref{Figure-Run-KV1a}, left plot).
\begin{figure}
\begin{center}
\includegraphics[height=32mm,width=33mm]{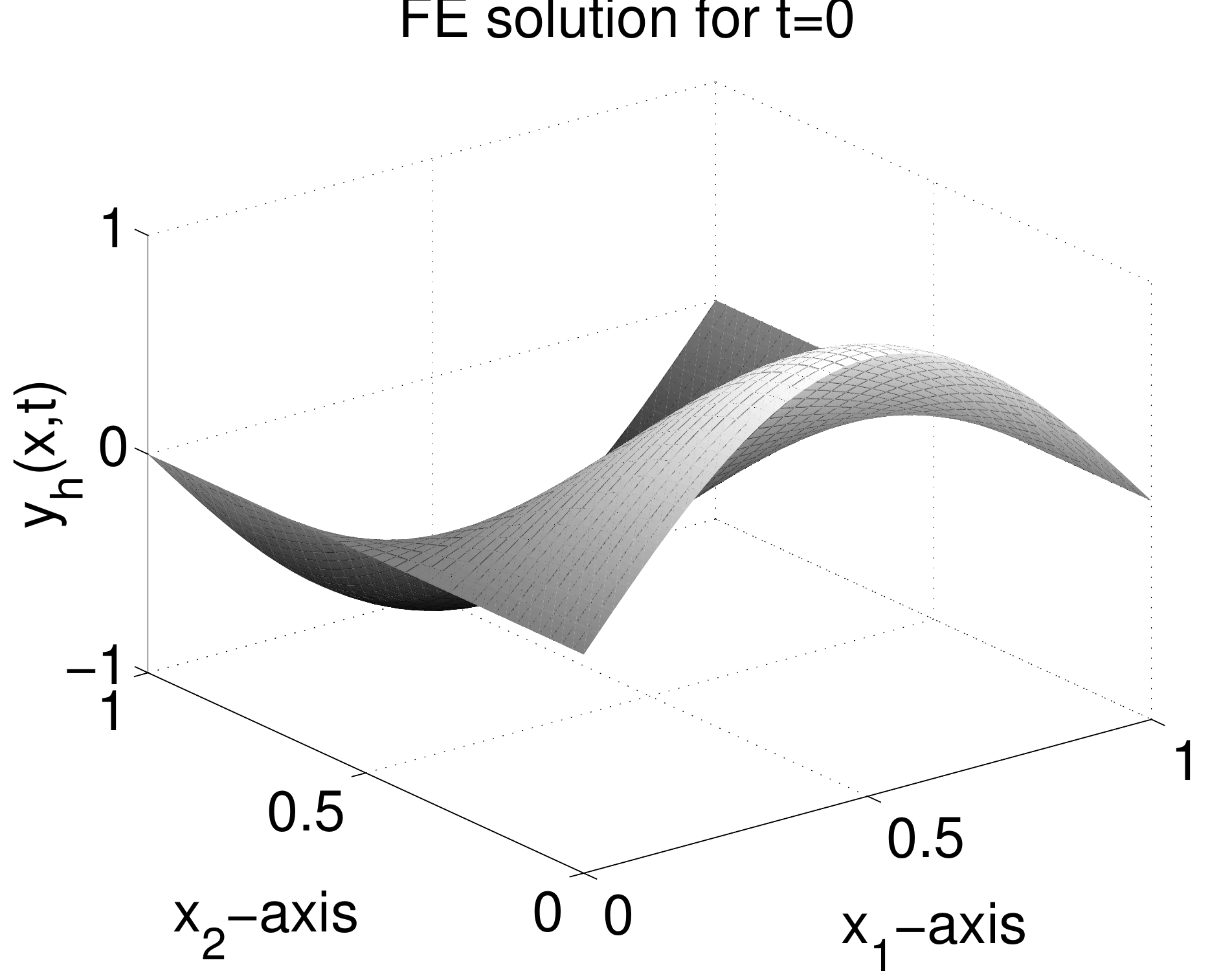}\hspace{2mm}
\includegraphics[height=32mm,width=33mm]{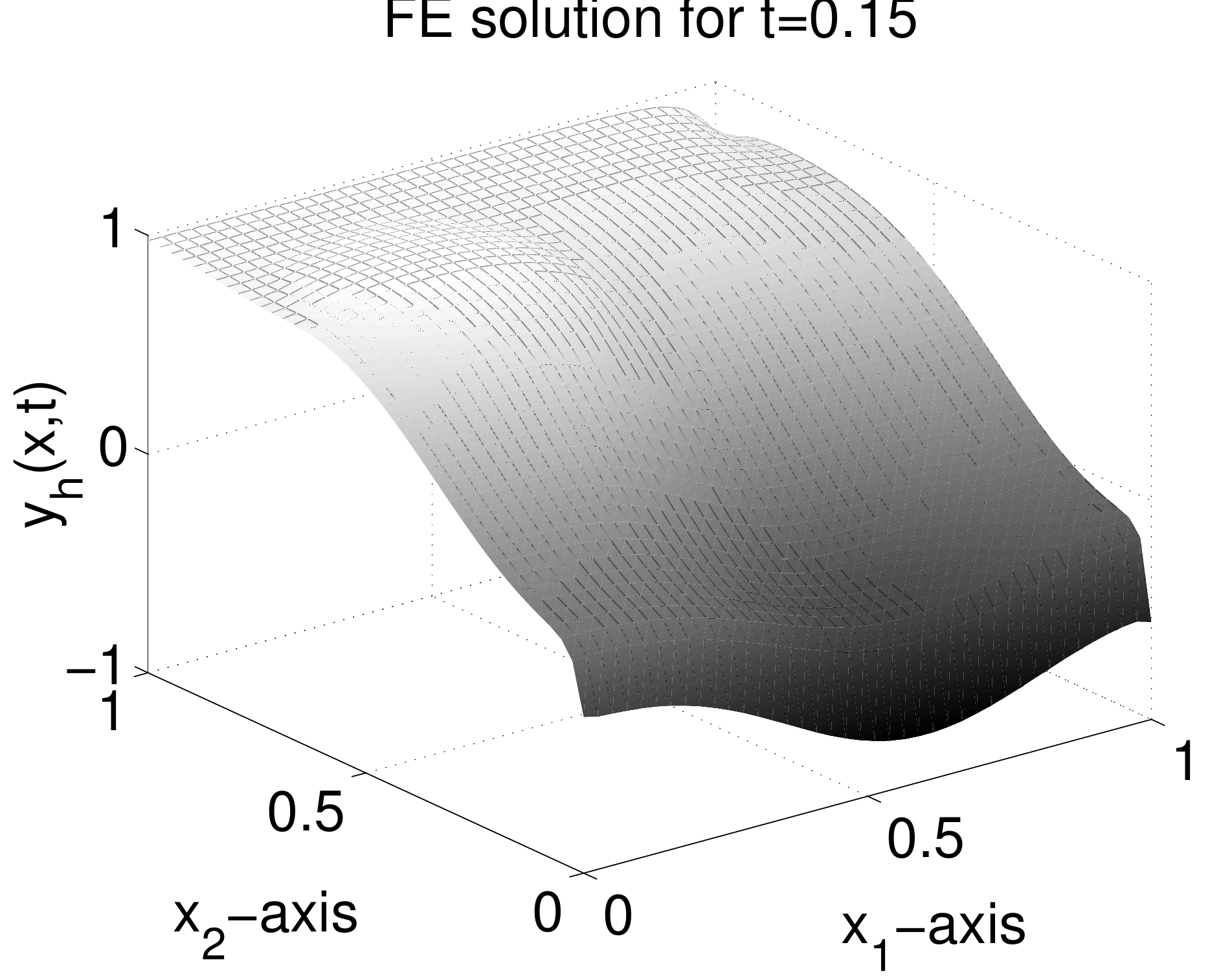}\hspace{2mm}
\includegraphics[height=32mm,width=33mm]{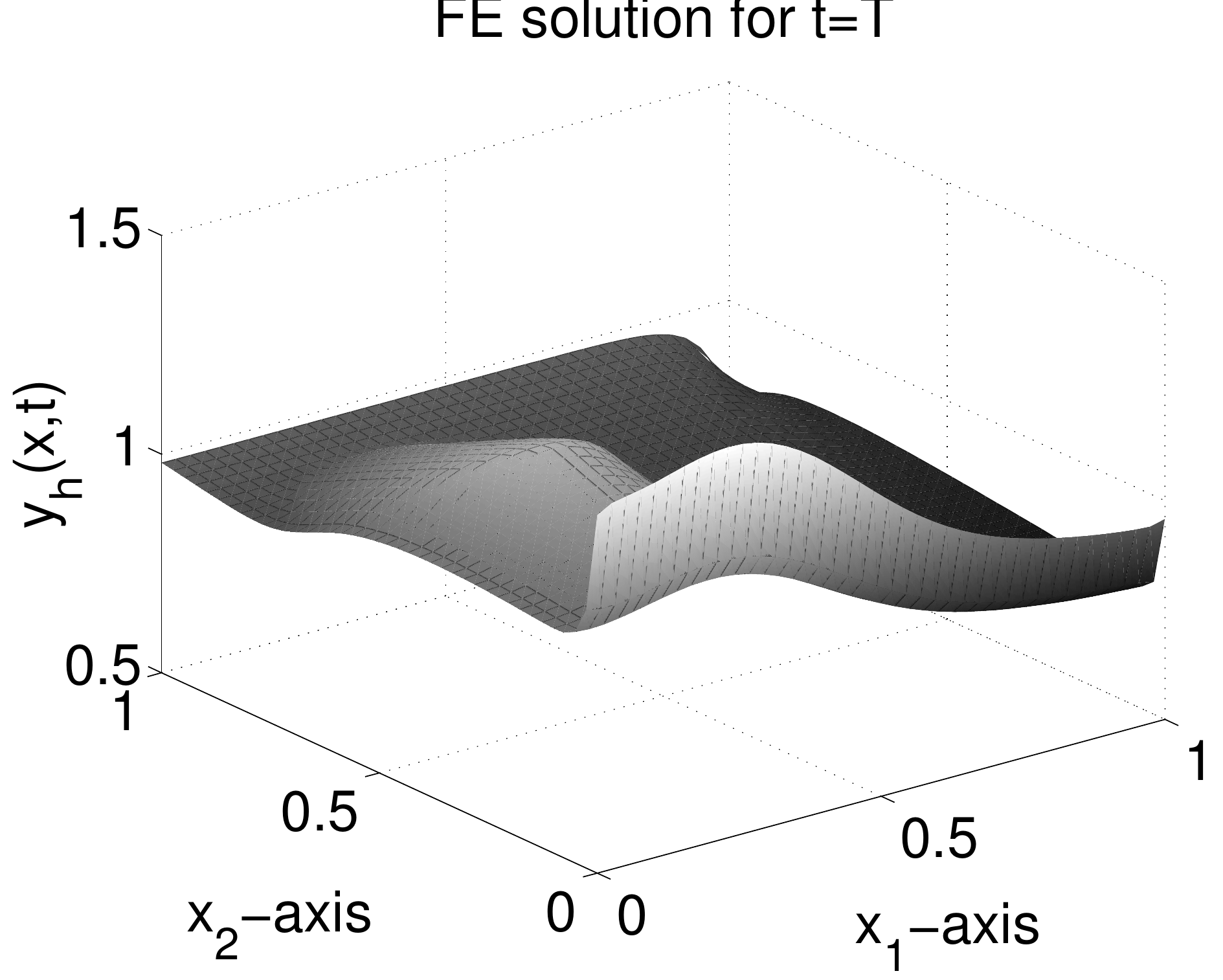}
\end{center}
\caption{Run~\ref{Run-KV1} (cf. \cite[Figures~3 and 4]{KV10}). initial condition $y_\circ$ (left plot), FE solution $y^h$ for $t=0.3$ (middle) and $t=T$ (right plot).}
\label{Figure-Run-KV1a}
\end{figure}
Furthermore, we have
\begin{align*}
q(\bx)&=\left\{
\begin{array}{cl}
1 & \text{for } \bx=(x_1,1) \text{ with } 0 < x_1 < 1,\\
x_2 & \text{for } \bx=(1,x_2) \text{ with } 0 < x_2 < 1,\\
-2 &\text{for } \bx=(x_1,0) \text{ with } 0 < x_1 < 1,\\
0 & \text{for } \bx=(0,x_2) \text{ with } 0 < x_2< 1,
\end{array}
\right.\\
g(\bx)&=\left\{
\begin{array}{cl}
1 & \text{for } \bx=(x_1,1) \text{ with } 0 < x_1 < 1,\\
0 & \text{for } \bx=(1,x_2) \text{ with } 0 < x_2 < 1,~\text{for } \bx=(0,x_2) \text{ with } 0 < x_2< 1,\\
-1 &\text{for } \bx=(x_1,0) \text{ with } 0 < x_1 < 1.
\end{array}
\right.
\end{align*}
We utilize piecewise linear FE functions. The FE solutions $y^h=y^h(t,\bx)$ for $t=0.15$ and $t=T$ are shown in Figure~\ref{Figure-Run-KV1a}. Next we take snapshots on the fixed uniform time grid $t_j=(j-1)\Delta t$, $1 \le j \le n_t$, with $n_t=10$ and $\Delta t=T/n_t=0.1$. The goal is to determine four additional time instances $\bar t=(\bar t_1,\ldots,\bar t_4) \in [0,T]$ based on a FE approximation for \eqref{Poptsnap}. Since the behavior of the solution exhibits more change during the initial time interval $[0,0.3]$ than later on, we initialize our Quasi-Newton method by the starting value $\tau^0=(0.05,0.15,0.25,0.35) \in [0,T]$. The number of POD ansatz functions is fixed to be $\ell=3$. The corresponding value of the ROM error is approximately $0.1093$. The optimal solution is given as $\bar\tau=(0.0092,0.0076,0.1336,0.2882) \in [0,T]$, while the associated ROM error is approximately $0.0165$, which is a reduction of about 85\,\%. In Figure~\ref{Figure-Run-KV1b} we can see that the shapes of the three POD bases changes significantly from the initial time instances $\tau^0\in\mathbb R^4$ to the optimal ones $\bar\tau\in\mathbb R^4$.\hfill$\Diamond$
\begin{figure}
\begin{center}
\includegraphics[height=32mm,width=33mm]{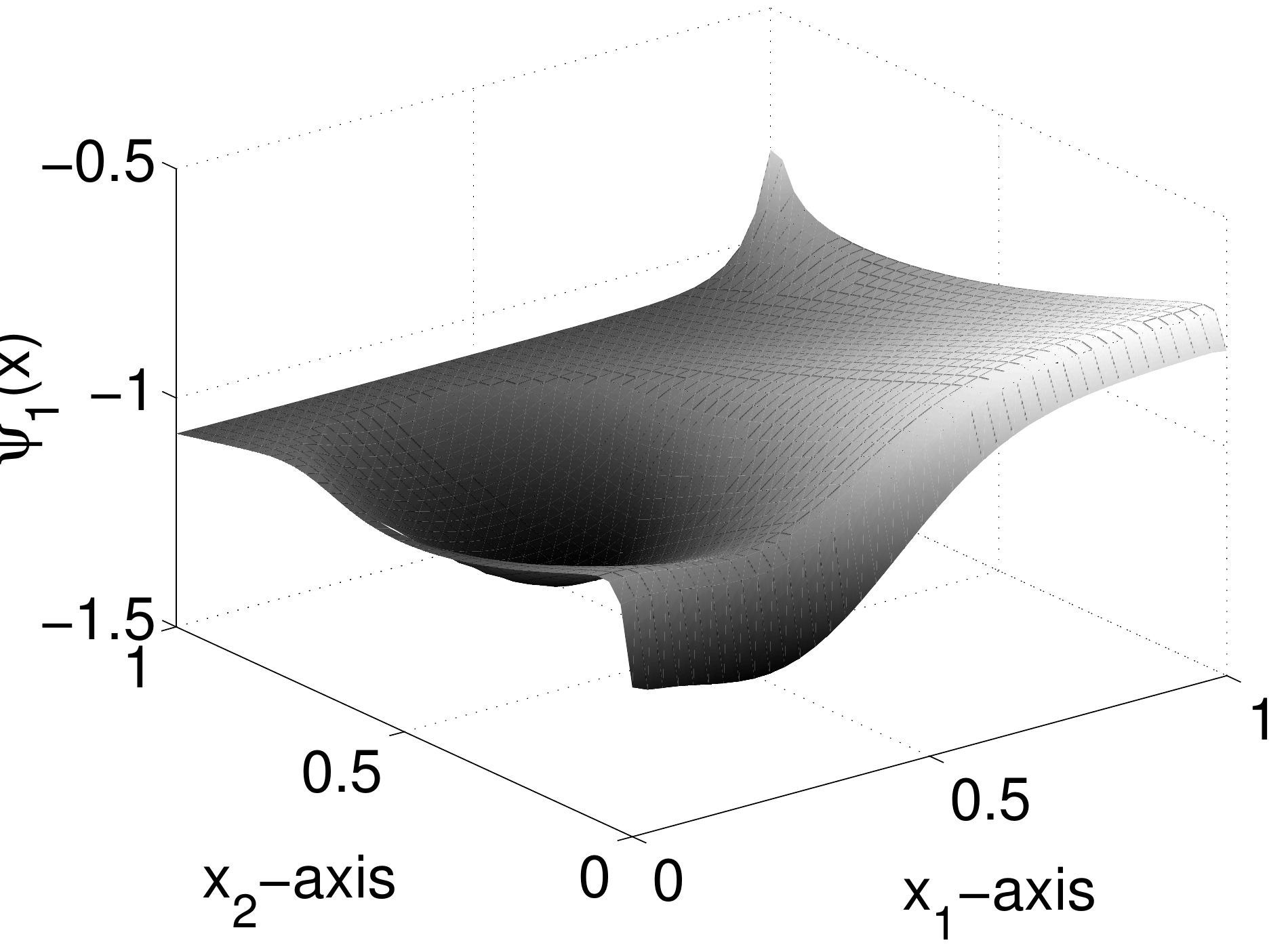}\hspace{2mm}
\includegraphics[height=32mm,width=33mm]{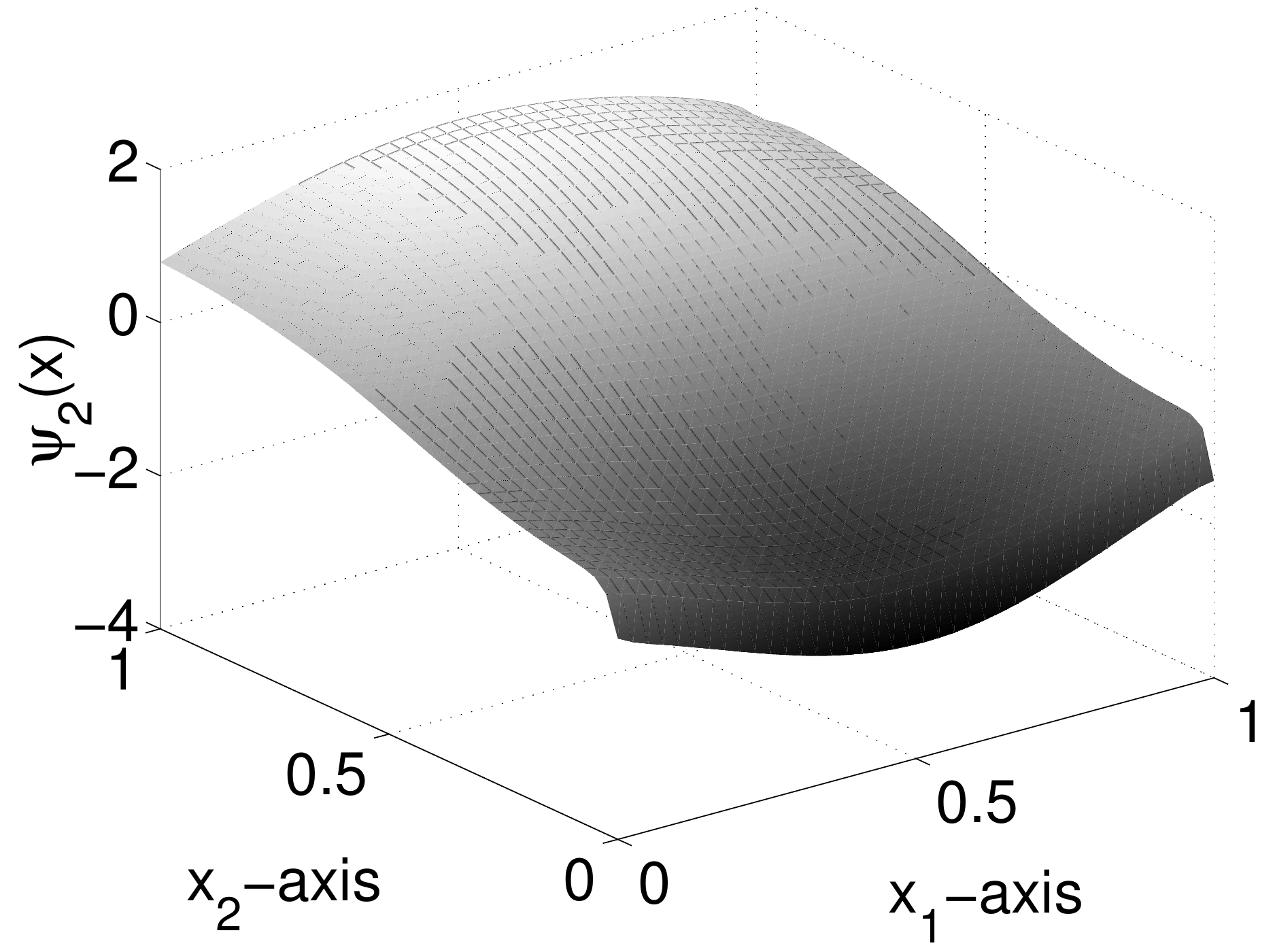}\hspace{2mm}
\includegraphics[height=32mm,width=33mm]{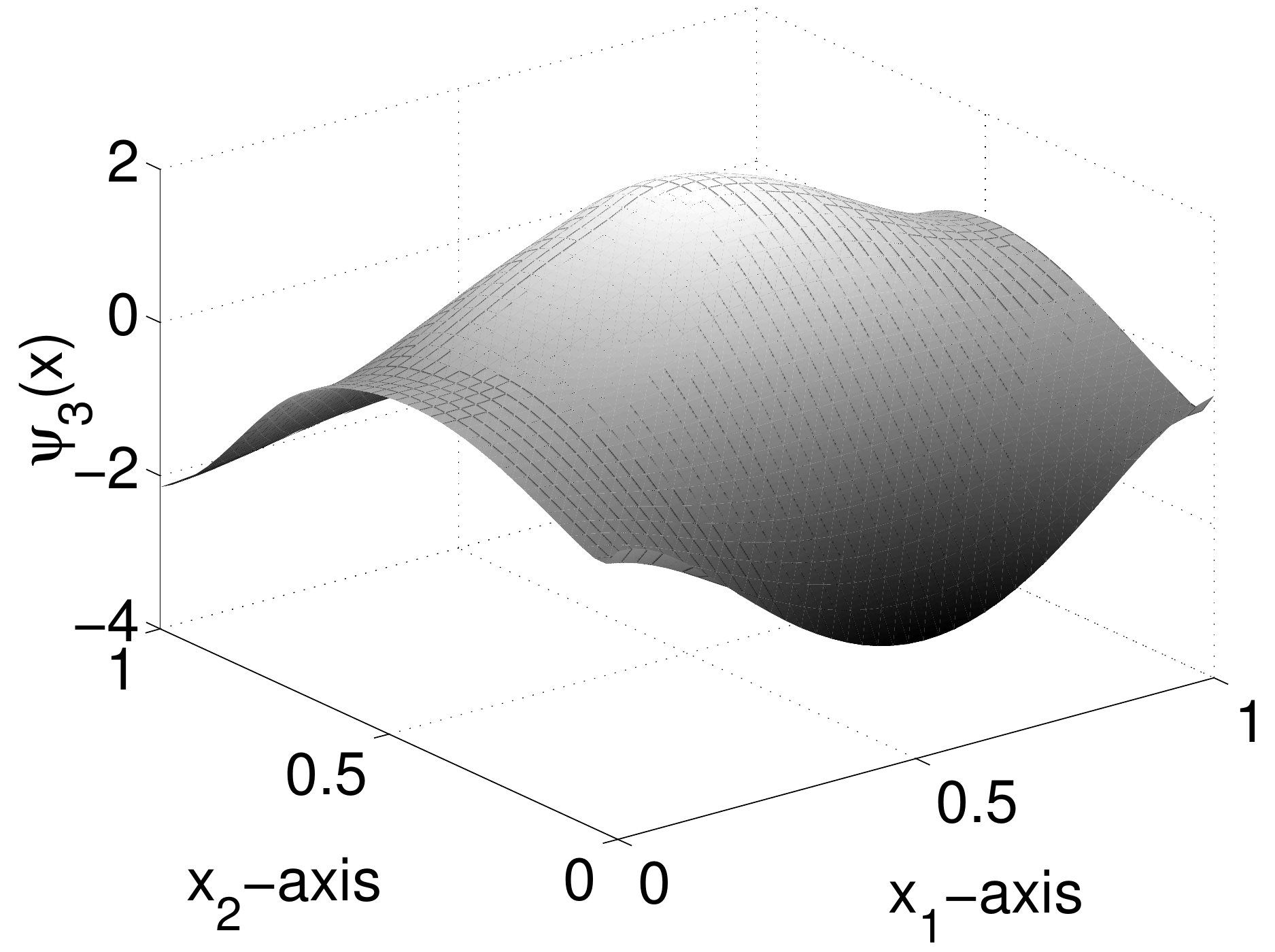}\\
\includegraphics[height=32mm,width=33mm]{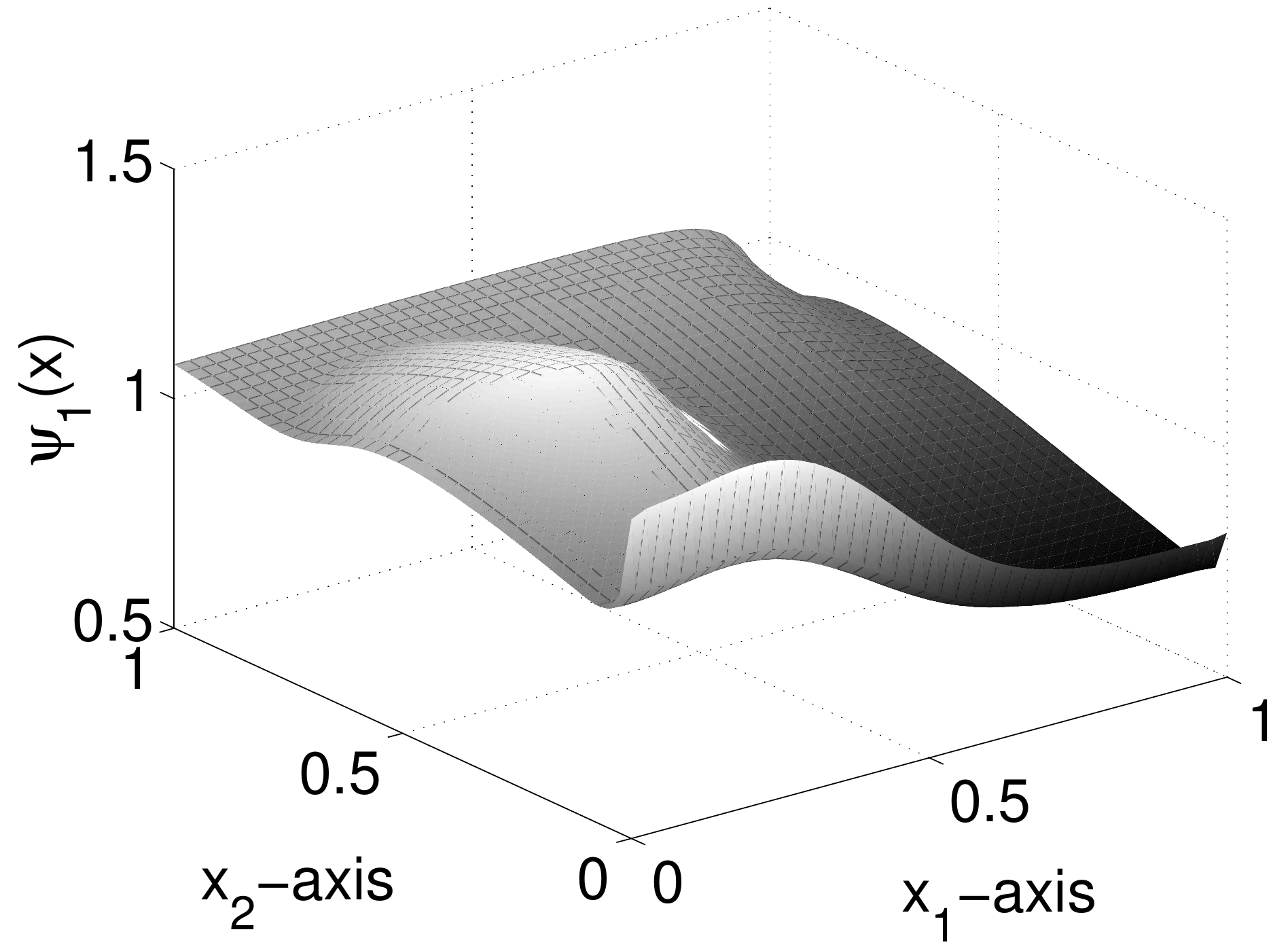}\hspace{2mm}
\includegraphics[height=32mm,width=33mm]{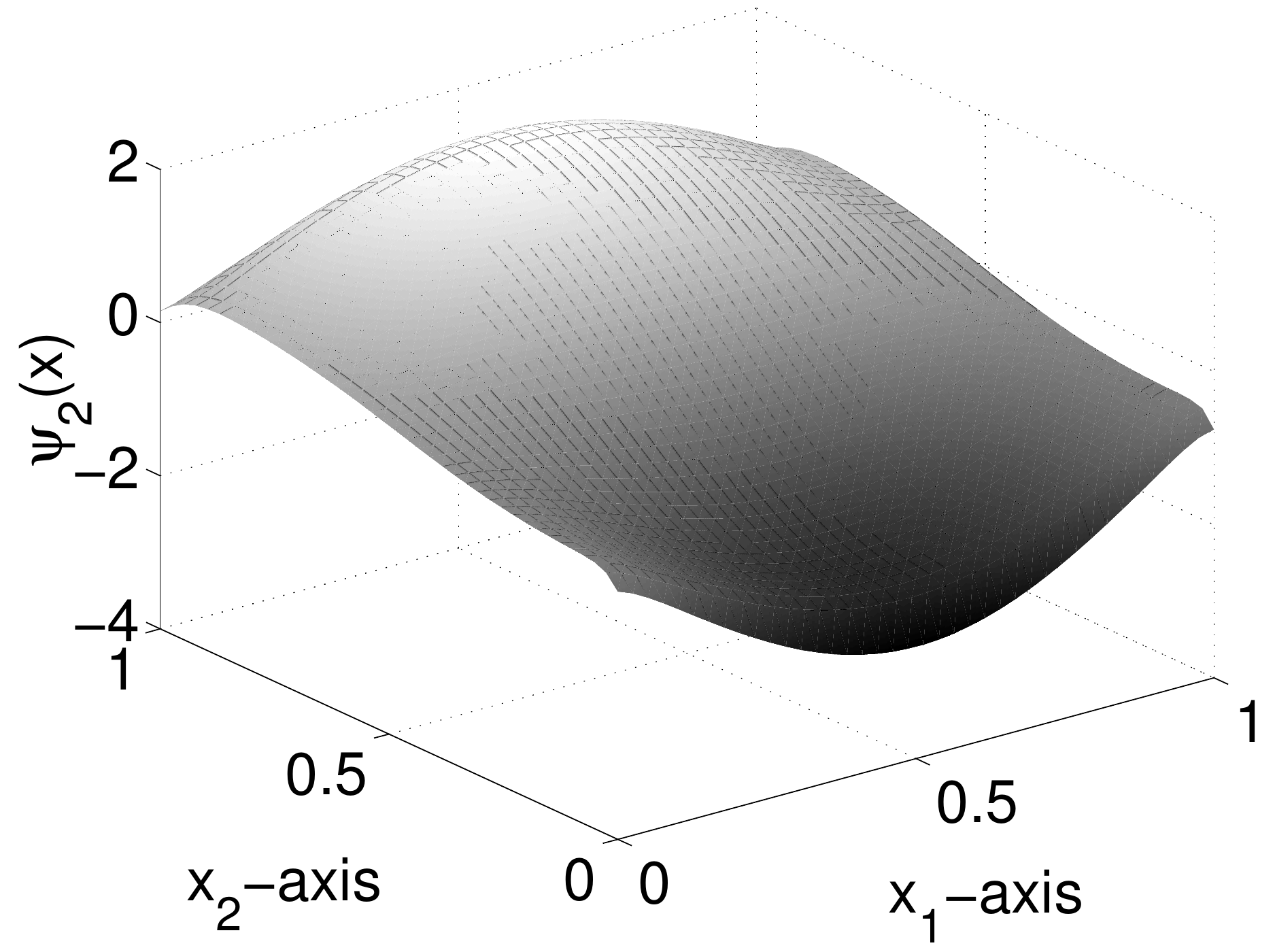}\hspace{2mm}
\includegraphics[height=32mm,width=33mm]{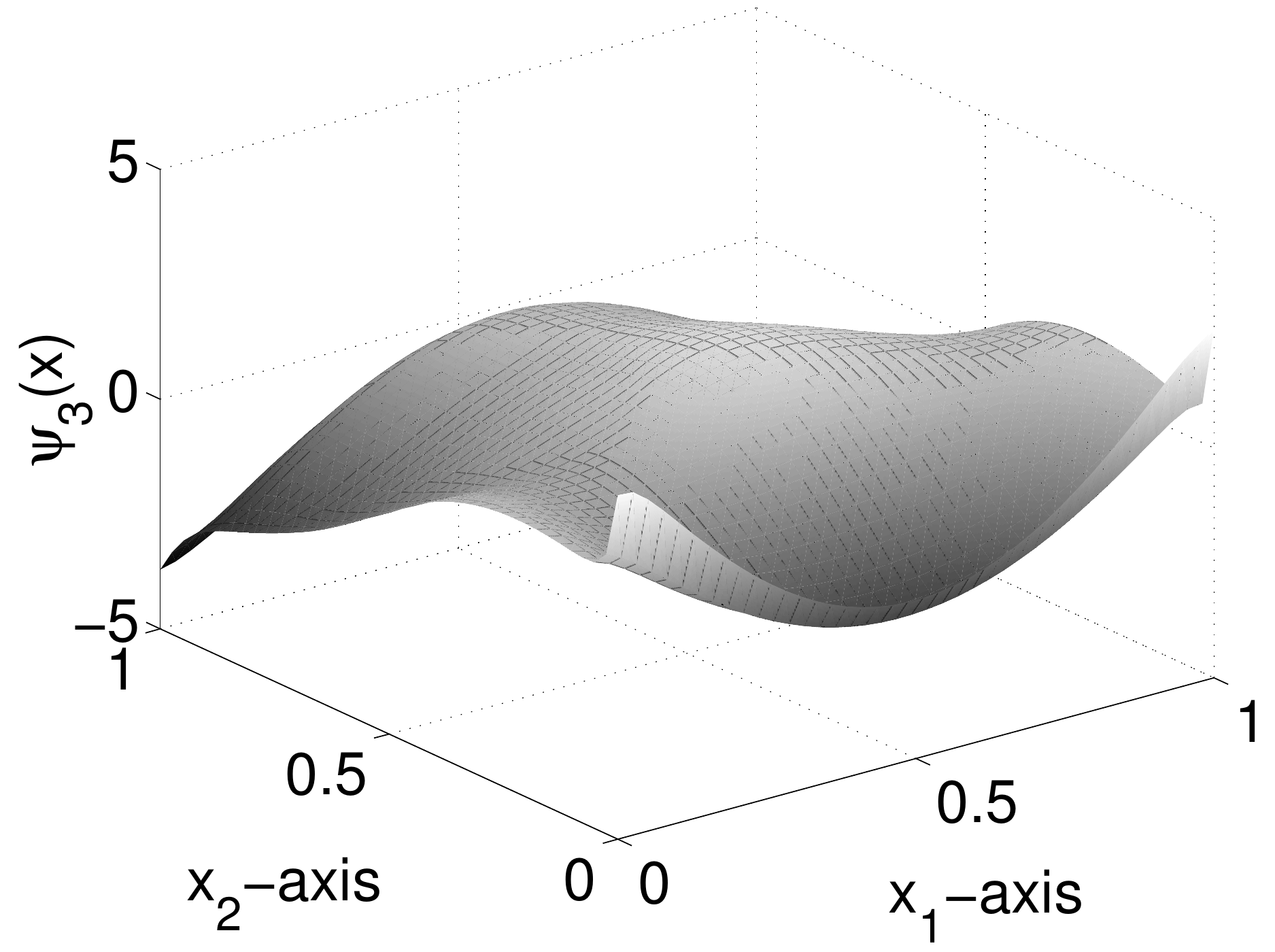}
\end{center}
\caption{Run~\ref{Run-KV1} (cf. \cite[Figures~5 and 7]{KV10}). POD basis $\Psi_1$, $\Psi_2$, $\Psi_3$ for the initial additional time instances $\tau^0\in\mathbb R^4$ (upper three plots) and for the optimal additional time instances $\bar\tau\in\mathbb R^4$ (lower three plots).}
\label{Figure-Run-KV1b}
\end{figure}
\end{run}

%%%%%%%%%%%%%%%%%%%%%%%%%%%%%%%%%%%%%%%%%%%%%%%%%%%
\section{Optimal control with POD surrogate models}
\label{POD-opt}
%%%%%%%%%%%%%%%%%%%%%%%%%%%%%%%%%%%%%%%%%%%%%%%%%%%

Reduced-order models are used in PDE-constrained optimization in various ways; see, e.g., \cite{HV05,SV10} for a survey. In optimal control problems it is sometimes necessary to compute a feedback control law instead of a fixed optimal control. In the implementation of these feedback laws models of reduced-order can play an important, and very useful role, see \cite{ABK01,GMPV19,KVX04,LV06,MV19,Rav00}. Another useful application is the use in optimization problems, where a PDE solver is part of the function evaluation. Obviously, thinking of a gradient evaluation or even a step-size rule in the optimization algorithm, an expensive function evaluation leads to an enormous amount of computing time. Here, the reduced-order model can replace the system given by a PDE in the objective function. It is quite common that a PDE can be replaced by a five- or ten-dimensional system of ordinary differential equations. This results computationally in a very fast method for optimization compared to the effort for the computation of a single solution of a PDE. There is a large amount of literature in engineering applications in this regard, we mention only the papers \cite{LT01,NAMTT03}. Recent applications can also be found in finance using the reduced models generated with the reduced basis (RB) method \cite{Pir06} and the POD model \cite{SS13,Sch12} in the context of calibration for models in option pricing.

We refer to the survey article \cite{GV17}, where a linear-quadratic optimal control problem in an abstract setting is considered. Error estimates for the POD Galerkin approximations of the optimal control are proved. This is achieved by combining techniques from \cite{DH02,DH04,Hin05} and \cite{KV01,KV02}. For nonlinear problems we refer the reader to \cite{HV05,Pin08,SV10}. However, unless the snapshots are generating a sufficiently rich state space or are computed from the exact (unknown) optimal controls, it is not a priorly clear how far the optimal solution of the POD problem is from the exact one. On the other hand, the POD method is a universal tool that is applicable also to problems with time-dependent coefficients or to nonlinear equations. Moreover, by generating snapshots from the real (large) model, a space is constructed that inhibits the main and relevant physical properties of the state system. This, and its ease of use makes POD very competitive in practical use, despite of a certain heuristic flavor. In this context results for a POD a posteriori analysis are important, see e.g., \cite{TV09} and \cite{GV13,KV12,KTV13,SV13,TUV11,Vol11,VV12}. Using a fairly standard perturbation method it is deduced how far the suboptimal control, computed on the basis of the POD model, is from the (unknown) exact one. This idea turned out to be very efficient in our examples. It is able to compensate for the lack of a priori analysis for POD methods. Let us also refer to the papers \cite{Ded10,GK13,NRMQ13}, where a posteriori error bounds are computed for linear-quadratic optimal control problems approximated by the reduced basis method.

Data- and/or simulation-based POD models depend on the data (e.g. initial values, right hand sides, boundary conditions, oberservations, etc.) which is used to generate the snapshots. If those models are used as surrogates in e.g. optimization problems with PDE constraints the algorithmical framework has to account for this fact with providing mechanisms for accordingly updating the surrogate model during the solution process. Strategies proposed in this context for optimal flow control can be found in e.g. \cite{AH99,AH01,AFS00,Fahl2001,BCB05}. One of the most mature methods developed in this context is Trust-Region POD proposed in \cite{AFS00}, which since then has successfully been applied in many applications. We also refer to the work \cite{GGMRV17}, where strategies for updating the POD bases are compared.

The quality of the surrogate model highly depends on its information basis, which for snapshot-based methods is given by the snapshot set, compare Section~\ref{POD-snap}. The location of snapshots and also the choice of the initial control in surrogate-based optimal control is discussed in \cite{AGH18}. There, techniques from time-adaptive schemes for optimality systems of parabolic optimal control problems are adjusted to compute optimal time locations for snapshots generation in POD surrogate modeling for parabolic optimal control problems.

Concepts for the construction and use of POD surrogate modeling in robust optimal control of electrical machines are presented in \cite{LaU16,AHLKU18}. Those problems are governed by nonlinear partial differential equations with uncertain parameters, so that robustness can be achieved by considering a worst case formulation. The resulting optimization problem then is of bilevel structure and POD reduced-order models in combination with a posteriori error estimators are used to speed up the numerical computations.

\section{Miscellaneous}
\label{POD-misc}
POD model order reduction can also be applied to provide surrogate models for high-fidelity components in networks. The general perspective is discussed in e.g. \cite{HM13}. Related research for MOR of electrical networks is reported in e.g. \cite{BHtM11,HKMV14,HKSS12}. The basic idea here consists in a decoupling of MOR approaches for the network and high-fidelity components which in general are modeled by PDE systems. For the latter, simulation-based POD MOR techniques are used to construct surrogate models which then are stamped back into the (reduced) electrical network. Details and performance tests are reported e.g. in \cite{HK12,HKSS12}. A short lecture series with related topics is presented under \href{https://slideslive.com/38894790/mathematical-aspects-of-proper-orthogonal-decomposition-pod-iii}{Hinze-Pilsen}\footnote{\url{https://slideslive.com/38894790/mathematical-aspects-of-proper-orthogonal-decomposition-pod-iii}}. Further contributions to this topic can be found in \cite{B17}. 

Recent trends in data-driven and nonlinear MOR methods are discussed within a YouTube lecture series under \href{https://www.youtube.com/watch?v=KOHxCIx04Dg}{Carlberg-YouTube}\footnote{\url{https://www.youtube.com/watch?v=KOHxCIx04Dg}}.

\bibliographystyle{plain}
%\bibliography{References}

\end{document}